\documentclass{article}
\usepackage{arxiv}

\usepackage[utf8]{inputenc}
\usepackage[T1]{fontenc}

\usepackage{amsmath,amssymb,mathtools,mathrsfs}
\usepackage{amsthm}
\usepackage{hyperref}
\usepackage{url}
\usepackage{xcolor}
\usepackage{microtype}
\usepackage{graphicx}
\usepackage{natbib}
\setcitestyle{authoryear,round,semicolon,aysep={,},yysep={,}}
\usepackage{subfigure}
\usepackage{booktabs}
\usepackage{array}
\usepackage{tabularx}
\usepackage[shortlabels]{enumitem}
\usepackage{placeins}
\usepackage{algorithm}
\usepackage{algpseudocode}
\usepackage{etoolbox}
\usepackage{authblk}

\newcolumntype{P}[1]{>{\raggedright\arraybackslash}p{#1}}
\newcolumntype{Y}{>{\raggedright\arraybackslash}X}
\allowdisplaybreaks

\newtheorem{theorem}{Theorem}[section]
\newtheorem{lemma}[theorem]{Lemma}
\newtheorem{proposition}[theorem]{Proposition}
\newtheorem{corollary}[theorem]{Corollary}
\theoremstyle{definition}
\newtheorem{definition}[theorem]{Definition}
\newtheorem{assumption}[theorem]{Assumption}

\theoremstyle{remark}
\newtheorem{remark}[theorem]{Remark}

\providecommand{\dist}{\mathop{\rm dist}\nolimits}

\title{Continuized Quasar-Convex Acceleration via an\\
Unbiased Finite-Mixing Network Oracle}

\author[1,2]{Ilya Kuruzov\thanks{\texttt{kuruzov.ia@phystech.edu}}}
\author[3]{Artem Agafonov\thanks{\texttt{artem.agafonov@mbzuai.ac.ae}}}
\author[1]{Vadim Abronin\thanks{\texttt{abronin.vd@phystech.edu}}}
\author[4,1,5]{Darina Dvinskikh\thanks{\texttt{dmdvinskikh@hse.ru}}}
\author[2,1,6]{Alexander Gasnikov\thanks{\texttt{gasnikov@yandex.ru}}}
\affil[1]{Moscow Independent Research Institute of Artificial Intelligence, Moscow, Russia}
\affil[2]{Innopolis University, Innopolis, Russia}
\affil[3]{Department of Machine Learning, Mohamed bin Zayed University of Artificial Intelligence, Abu Dhabi, United Arab Emirates}
\affil[4]{HSE University, Moscow, Russia}
\affil[5]{Trusted AI Research Center RAS, Moscow, Russia}
\affil[6]{Institute for Information Transmission Problems of the Russian Academy of Sciences, Moscow, Russia}
\date{}

\renewcommand{\shorttitle}{Continuized Quasar-Convex Acceleration}

\hypersetup{
  colorlinks=true, linkcolor=blue!55!black, citecolor=blue!55!black, urlcolor=blue!55!black,
  pdftitle={Continuized Quasar-Convex Acceleration via an Unbiased Finite-Mixing Network Oracle},
  pdfauthor={Ilya Kuruzov, Artem Agafonov, Vadim Abronin, Darina Dvinskikh, Alexander Gasnikov},
  pdfkeywords={quasar-convexity, continuized acceleration, decentralized optimization, finite mixing, gossip}
}

\begin{document}
\maketitle

\begin{abstract}
We study decentralized minimization of a smooth network average that is
globally quasar-convex with respect to a fixed minimizer; the node functions
may be nonconvex and need not share a minimizer.  A unit-rate Poisson clock
drives an exactly simulable continuized process: between events the state
follows a closed-form accelerated flow, and each event uses one oracle vector.
One master theorem covers oracles with a predictable bias and a state-dependent
variance and keeps the inverse-quadratic term
$\mathcal O(LR_0^2/(\gamma^2T^2))$; one selected-row lemma supplies the
oracle moments of three network implementations.  CQA-Net keeps a common
state, mixes local rows for finitely many rounds, and disseminates one
uniformly selected row; its oracle is unbiased at every mixing depth, so
incomplete averaging affects only the variance.  A walking token removes
dissemination, the shared clock, and the shared root, and an importance weight
keeps it unbiased on fixed graphs.  D-CQA keeps local states and mixes them;
under bounded heterogeneity it needs
$\widetilde{\mathcal O}(\sqrt\chi\sqrt LR_0/(\gamma\theta\sqrt\varepsilon))$
communication rounds and no dissemination.  Changing graphs enter through
certified operator norms.  Experiments on a certified nonconvex regression on
static and changing graphs compare the three implementations.
\end{abstract}


\section{Introduction}\label{sec:introduction}

We consider a decentralized optimization problem over a possibly time-varying
network:
$$\min_{x \in \mathbb{R}^d}  F(x):=\frac1m\sum\limits_{i=1}^m f_i(x),$$
where node $i$ stores its own nonlinear-regression samples and the network
objective is the average of the local ones.  Such problems model learning from
data that cannot be collected in one place and arise in federated learning, in sensor and multi-agent systems, and in
large-scale training, where decentralized communication can be cheaper than
routing everything through a central server
\citep{nedic2009,mcmahan2017,lian2017}.

Even for a realizable monotone-link model, this objective need not be convex:
its Hessian can be indefinite.  Nevertheless, if the link has slope at least
$\alpha>0$ and is $M$-Lipschitz, the objective satisfies the one-point
inequality $\gamma[F(x)-F(x^\star)]\leq\langle\nabla F(x),x-x^\star\rangle$
with the globally valid parameter $\gamma=\min\{1,2\alpha/M\}$; see
Proposition~\ref{prop:model-certificate}.  This quasar-convexity condition
holds globally for generalized linear models with monotone links
(Section~\ref{sec:experiments}) and for star-convex objectives, which are
$1$-quasar-convex \citep{hinder2020}.  Weaker versions restricted to a domain
of parameters arise elsewhere, for instance in learning linear dynamical systems
\citep{hardt2018}; such cases need a projection or an invariance argument and are
not covered by our global theorems.  We ask whether the accelerated
inverse-quadratic rate survives when the gradient is assembled by only finitely
many peer-to-peer communication rounds.

Acceleration for convex decentralized problems is well studied
\citep{scaman2017,hendrikx2019,xu2020,kovalev2021adom,kovalev2022sliding,li2024agt},
but a convex decentralized proof cannot simply be copied.  Quasar-convexity
gives a lower model only for the pair $(x,x^\star)$, not for two arbitrary query
points, and it does not imply monotonicity of the gradient, so a consensus
error cannot be inserted into the usual affine estimate-sequence inequalities.

Our approach separates the acceleration mechanism from the network.  We first
prove a continuized acceleration result for an abstract inexact oracle with a
predictable bias and a state-dependent variance
(Theorem~\ref{thm:master-inexact}): after $T$ events in expectation, the
objective gap is $\mathcal O(L\|x_0-x^\star\|^2/(\gamma^2T^2))$ plus variance
and bias terms.  We then show that finitely many gossip rounds followed by the
selection of one node's mixed row implement such an oracle without bias, at
every mixing depth: unfinished averaging enters only as variance, which the
quasar-convex analysis can absorb.  We call this method CQA-Net.  It is
synchronized, blocking, and common-state, and it relies on three global primitives: exact
dissemination of the selected row, a shared event clock, and a shared uniform
root.  Dissemination is natural in federated learning, but in a peer-to-peer
network, and especially on a changing graph, it is a strong requirement.  We
therefore give two further implementations.  A walking token carries the common
state along a random walk, so that the clock and the root travel with the token
and nothing is disseminated.  D-CQA keeps local states and mixes them; only a
shared clock seed and synchronous rounds remain, and the dissemination depth
disappears from the round count.  Table~\ref{tab:nearby-models} compares the
three designs with related methods.  We do not address physical asynchrony,
fault recovery, partial participation, wall-clock time, or communication
optimality.

\begin{table}[!htbp]
\caption{Comparison with nearby network-optimization models.  ``Common'' means
that all nodes hold identical optimization states after each event; ``local''
allows node disagreement.}
\label{tab:nearby-models}
\centering
\small
\begin{tabular}{@{}P{2.4cm}P{1.6cm}P{2.6cm}P{2.3cm}P{4.75cm}@{}}
\toprule
Method & State & Timing & Geometry & Mechanism and result \\
\midrule
\textbf{CQA-Net} & Common & Blocking event and broadcast & Global quasar-convexity & Unbiased finite mixing; $\mathcal O(LR_0^2/(\gamma^2T^2))$ plus a variance term \\
\textbf{Token} & One mobile state & Sequential walk & Global quasar-convexity & Importance-weighted walk endpoint; no global primitive \\
\textbf{D-CQA} & Local & Synchronous rounds & Global quasar-convexity, bounded heterogeneity & State mixing; $\widetilde{\mathcal O}(\sqrt\chi\sqrt LR_0/(\gamma\theta\sqrt\varepsilon))$ rounds \\
Gradient tracking \citep{nedic2017,li2024agt} & Local & Synchronous rounds & Convex or strongly convex & Primal and tracker consensus \\
DADAO \citep{nabli2025dadao} & Local & Asynchronous clocks & Strongly convex & Independent Poisson communication and computation \\
\bottomrule
\end{tabular}
\end{table}

Section~\ref{sec:experiments} tests the methods on instances whose constants
are certified in closed form: a nonconvex leaky-$\tanh$ regression on static and
changing graphs, the individual variance mechanisms of CQA-Net, and a
comparison of CQA-Net, the token, and D-CQA by communication rounds.  It also quantifies how
conservative the global certificates can be.

\paragraph{Main contributions.}
\begin{enumerate}\setlength{\itemsep}{.25em}
\item Theorem~\ref{thm:master-inexact} gives one Lyapunov bound for oracles
with a predictable bias and a state-dependent variance.  It retains the variance
integral and a gradient reserve that absorbs multiplicative noise, and it
handles a mean-square bias by a square-root Gronwall argument, without
projection or damping.
\item Lemma~\ref{lem:uniform-row} gives the conditional moments of any
selected-row oracle: a mixed matrix of local rows, a root with an arbitrary
law, and an optional importance weight.  It yields unbiased finite mixing at
every depth (Theorem~\ref{thm:static-network}, fixed and changing graphs) and
the walking-token oracle (Corollaries~\ref{cor:token-fixed}
and~\ref{cor:token-changing}).
\item Theorem~\ref{thm:dcqa} analyzes the local-state variant D-CQA, and
Corollary~\ref{cor:dcqa-complexity} bounds its communication by
$\widetilde{\mathcal O}(\sqrt\chi\sqrt LR_0/(\gamma\theta\sqrt\varepsilon))$
rounds, without a dissemination term.
\item For shared-landscape finite sums, exact-snapshot restarts give a geometric
rate with an accuracy-independent inner depth and batch size
(Remark~\ref{rem:vr} and Appendix~\ref{app:vr}).
\end{enumerate}

\subsection{Related work}\label{sec:related-work}

\paragraph{Acceleration beyond convexity.}
Smooth convex acceleration originates with Nesterov
\citep{nesterov1983,nesterov2004}.  Accelerated methods for weakly quasi-convex
and star-convex objectives use auxiliary one- or low-dimensional searches
\citep{guminov2023,kuruzov2021sso,kuruzov2024subspace,hinder2020}, which are
costly across a network because every trial point may need a global
comparison.  Continuized acceleration avoids this search: \citet{even2021}
introduced exactly simulable Poisson-time acceleration, and
\citet{wang2023} obtained the optimal deterministic dependence on the quasar
parameter.  Stochastic, robust, non-Euclidean, constrained,
zeroth-order, proximal, and distributional developments appear in
\citep{fu2023,hermant2024,gupta2024agnes,lezane2026,martinezrubio2025,gbaguidi2026,jin2020,debrito2025,ahookhosh2026,ding2024gqc};
related one-point structures arise, on restricted domains, in learning linear
dynamical systems \citep{hardt2018}.  Nonconvex stochastic lower bounds
\citep{saad2025lower} do not yield a matching lower bound for our model.

\paragraph{Acceleration in networks.}
Communication-optimal convex decentralized methods combine objective and
spectral acceleration \citep{scaman2017,hendrikx2019,xu2020,kovalev2021adom},
with sliding, similarity, and second-order variants
\citep{kovalev2022sliding,agafonov2021,takezawa2025spdo,tian2022,daneshmand2021}.
Gradient tracking and primal-dual methods maintain local states
\citep{nedic2017,qu2018,li2024agt,hong2017}, DADAO uses asynchronous edge and
gradient clocks \citep{nabli2025dadao}, DCatalyst wraps a decentralized solver in
an accelerated outer loop \citep{cao2026dcatalyst}, and time-varying and
finite-sum settings are treated in
\citep{rogozin2019,metelev2023,metelev2024,lobanov2023,liu2024finite,lu2021}.
CQA-Net keeps one common state, so its round count contains the broadcast depth
$D$ in addition to the Chebyshev depth; at $\gamma=1$ it is a non-optimal
common-state specialization.  D-CQA removes $D$ and, under bounded
heterogeneity, matches the convex local-state round scale up to a logarithm.
Random-walk methods process one node per step along a walk
\citep{johansson2009randomized,mao2020walkman}; our token uses this mechanism
only to deliver the oracle, and push-sum-type protocols \citep{kempe2003gossip}
are a natural route to message loss, which we do not treat.

\paragraph{Node sampling.}
At $q=0$ CQA-Net is uniform node sampling followed by dissemination, as in
client sampling, Local SGD, and SCAFFOLD-type control variates
\citep{karimireddy2020scaffold,gorbunov2021localsgd,koloskova2020unified} and
randomized broadcast gossip \citep{boyd2006randomized}.  The uniform-row identity
is elementary; we use it to supply the exact conditional moments that the
quasar-specific Lyapunov theorem needs.

\paragraph{Organization.}
Section~\ref{sec:problem} sets the geometry and the stochastic model,
Section~\ref{sec:algorithm-design} defines the algorithms,
Section~\ref{sec:main-results} contains the analysis, and
Section~\ref{sec:experiments} the numerical study.  Technical proofs and
reproducibility details are in the appendices.

\section{Problem statement and notation}\label{sec:problem}

We separate the optimization state from the network data used to construct one
oracle vector.  The continuized dynamics live in $\mathbb R^d$ and are written
with lowercase variables.  A node index appears only on local rows and on the
transient matrices used by the communication protocol.  In CQA-Net every node
stores the same optimization vectors at the end of each event;
Section~\ref{sec:dcqa-design} introduces a variant with local states.

\subsection{The one-vector problem}\label{sec:central-problem}

Let $F:\mathbb R^d\to\mathbb R$ be differentiable and let $x^\star$ be a
fixed global minimizer.  Continuous-time states and the oracle lie in
$\mathbb R^d$ and are denoted $x_t,z_t,g_t$; all processes are right-continuous
with left limits.  For a jump index $k\geq1$,
\[
x^{(k)}=x_{\tau_k},\qquad y^{(k)}=x_{\tau_k-},\qquad
z^{(k)}=z_{\tau_k},\qquad g^{(k)}=g_{\tau_k},
\]
so $y^{(k)}$ is the pre-jump query and $x^{(k)}$ the post-jump state, and
$\mathcal F_{\tau_k-}$ collects everything revealed strictly before $\tau_k$.
We set $x^{(0)}=x_0$, $z^{(0)}=z_0$, $\tau_0=0$, $E_{k-1}=\tau_k-\tau_{k-1}$,
and $R_0:=\|x_0-x^\star\|$.

\begin{definition}[Smoothness]\label{def:smooth}
The objective is $L$-smooth if
$\|\nabla F(x)-\nabla F(y)\|\leq L\|x-y\|$ for all $x,y\in\mathbb R^d$.
Then the descent inequality
\begin{equation}\label{eq:smoothness}
F(y)\leq F(x)+\langle\nabla F(x),y-x\rangle+\frac L2\|y-x\|^2.
\end{equation}
holds for all $x,y\in\mathbb R^d$.
\end{definition}

\begin{definition}[Quasar-convexity / weak quasi-convexity]\label{def:quasar}
For $\gamma\in(0,1]$, the objective is $\gamma$-quasar-convex with respect
to $x^\star$ if
\begin{equation}\label{eq:quasar}
\gamma\bigl(F(x)-F(x^\star)\bigr)
\leq\langle\nabla F(x),x-x^\star\rangle
\qquad\text{for every }x\in\mathbb R^d.
\end{equation}
The number $\gamma$ is any certified admissible value.  Every statement remains
valid after replacing it by a smaller positive value.
\end{definition}

\begin{definition}[Quadratic growth]\label{def:qg}
The objective satisfies fixed-minimizer quadratic growth with constant $\mu>0$
if
\begin{equation}\label{eq:qg}
F(x)-F(x^\star)\geq\frac\mu2\|x-x^\star\|^2
\qquad\text{for every }x\in\mathbb R^d.
\end{equation}
Condition~\eqref{eq:qg} uses the same distinguished minimizer as quasar-convexity.
\end{definition}

Convexity implies Definition~\ref{def:quasar} with $\gamma=1$, but the converse
is false.  The condition supplies a lower model only in the pair $(x,x^\star)$:
it gives neither an affine lower model at an arbitrary second point nor
monotonicity of the gradient.  It nevertheless implies that every stationary
point is globally optimal.  Quadratic growth is invoked only for the dissipative closure of
Corollary~\ref{cor:closed-static}(c).

\subsection{Certified quasar-convex models}\label{sec:quasar-overview}

A concrete globally nonconvex class is realizable nonlinear regression with a
strongly monotone link.  Let $\varphi:\mathbb R\to\mathbb R$ be a
differentiable link function; quantitative assumptions on $\varphi'$ are
stated in Proposition~\ref{prop:model-certificate}.  Given features
$a_r\in\mathbb R^d$ and a target $x^\star$, the responses are realizable,
$y_r=\varphi(a_r^\top x^\star)$, and the squared-loss objective is
\begin{equation}\label{eq:link-objective-main}
F(x)=\frac1{2N}\sum_{r=1}^N
\bigl(\varphi(a_r^\top x)-y_r\bigr)^2.
\end{equation}Proposition~\ref{prop:model-certificate} gives the constants used in the
numerical section.  Appendix~\ref{app:quasar-certificates} contains the proof.

\begin{proposition}[Monotone-link and leaky-$\tanh$ certificates]
\label{prop:model-certificate}
Assume that $\varphi$ is differentiable and
$0<\alpha\leq\varphi'(u)\leq M$ for every $u$.  Then
\eqref{eq:link-objective-main} is globally $\gamma$-quasar-convex with
$\gamma=\min\{1,2\alpha/M\}$.  If
$\Sigma=N^{-1}\sum_ra_ra_r^\top\succ0$, it satisfies quadratic growth with
$\mu=\alpha^2\lambda_{\min}(\Sigma)$.

For $\alpha\in(0,1]$ and
$\varphi_\alpha(u)=\alpha u+(1-\alpha)\tanh u$, let
$v_r=a_r^\top x^\star$ and
\[
B_\alpha(v)=1+\frac{4\alpha(1-\alpha)}e
 +\frac{4\alpha(1-\alpha)}{3\sqrt3}|v|
 +\frac{8(1-\alpha)^2}{3\sqrt3}.
\]
Then $\nabla F$ is globally Lipschitz with any
$L\geq N^{-1}\sum_rB_\alpha(v_r)\|a_r\|^2$.
\end{proposition}

The quasar part sharpens Lemma~3 of \citet{wang2023}, which gives the constant $2\alpha^2/M^2$ under the same assumptions, by the factor $M/\alpha$; the quadratic-growth constant matches their quasar-growth bound.  The certificate is global and permits negative curvature.  In the scalar case
$N=d=1$, $\alpha=1/2$, and $x^\star=-2$, direct interval evaluation gives
$-0.253<F''(1)<-0.250$.  Realizability is essential: response noise creates a
nonzero residual at $x^\star$ and requires a separate one-point argument.

\paragraph{Certified examples used below.}
The smooth convex case has $\gamma=1$ and serves as a calibration and ridge
benchmark.  The realizable leaky-$\tanh$ regression instance is globally
quasar-convex and satisfies quadratic growth despite being nonconvex; its
rounded constants are $\gamma=0.36$, $L=2.3631$, and
$\mu=1.457\cdot10^{-3}$.  It is used for the static- and changing-network
experiments.  The aligned multiplicative quadratic has a
convex network average with $\gamma=\mu=L=1$ but nonconvex local terms and
heterogeneity coefficient $h_1=1.28$.  Full-precision constants are stored in the accompanying metadata.

\subsection{Network average and local rows}\label{sec:network-problem}

Node $i\in\{1,\ldots,m\}$ stores a differentiable local function
$f_i:\mathbb R^d\to\mathbb R$, and
\begin{equation}\label{eq:network-average}
F(x)=\frac1m\sum_{i=1}^m f_i(x).
\end{equation}
Only the average $F$ is assumed smooth and quasar-convex; the $f_i$ may be
nonconvex, may fail the quasar condition, and may have different minimizers.
Weighted empirical partitions are encoded in the definitions of the $f_i$, so
\eqref{eq:network-average} remains exact when local sample counts differ.

At jump $k\geq1$, all nodes receive the same predictable query $y^{(k)}\in\mathbb
R^d$.  Node $i$ produces a random row $g_i^{(k)}(y^{(k)})\in\mathbb R^d$ with
\begin{equation}\label{eq:local-oracle}
\mathbb E[g_i^{(k)}(y^{(k)})\mid\mathcal F_{\tau_k-}]
 =\nabla f_i(y^{(k)}),\qquad
\mathbb E\!\left[\|g_i^{(k)}(y^{(k)})-\nabla f_i(y^{(k)})\|^2
 \mid\mathcal F_{\tau_k-}\right]\leq\frac{\sigma^2}{b}.
\end{equation}
Write $\mathsf K_{i,k}$ for a regular conditional kernel, given
$\mathcal F_{\tau_k-}$, of the row rule prescribed at that event (either the
ordinary local row above or a control-variate row introduced later).  Whenever
a root $R_k$ is used, it is conditionally uniform and conditionally
independent, given $\mathcal F_{\tau_k-}$, of a latent family of row
innovations having kernels
$(\mathsf K_{i,k})_{i=1}^m$.  Equivalently, in the lazy $q=0$ implementation,
after $R_k=i$ is revealed the only evaluated row is sampled from
$\mathsf K_{i,k}$ using fresh randomness independent of the root, so that
\begin{equation}\label{eq:lazy-root-row-kernel}
\mathcal L\!\left(g_i^{(k)}(y^{(k)})\mid
 \mathcal F_{\tau_k-},R_k=i\right)=\mathsf K_{i,k}
 =\mathcal L\!\left(g_i^{(k)}(y^{(k)})\mid
 \mathcal F_{\tau_k-}\right).
\end{equation}
This latent coupling is a statement about the conditional law; it does not
require unelected nodes to evaluate gradients.  Cross-node row innovations may
be dependent unless a later covariance bound explicitly requires otherwise.
The stack of rows is $G^{(k)}\in\mathbb R^{m\times d}$, and $e_i\in\mathbb
R^m$ denotes the $i$th standard basis vector.  The single vector
returned after finite mixing, row selection, and broadcast is denoted
$\widehat g^{(k)}\in\mathbb R^d$.  The squared heterogeneity at a common point
is
\begin{equation}\label{eq:heterogeneity}
\zeta^2(y)=\frac1m\sum_{i=1}^m
\|\nabla f_i(y)-\nabla F(y)\|^2.
\end{equation}
The quantity in~\eqref{eq:heterogeneity} measures disagreement of local mean rows, not disagreement of node states.

\subsection{Analytical and simulation filtrations}\label{sec:oracle-model}

Let $n_t$ be a unit-rate Poisson counting process with event times
$0<\tau_1<\tau_2<\cdots$.  The analysis is carried out in a filtration
$\mathbb F=(\mathcal F_t)_{t\geq0}$ relative to which $n_t$ has predictable
intensity one.  A random kernel $\mathsf K_t(\omega,dg)$ is
\emph{$\mathbb F$-predictable} if, for every bounded Borel function $\phi$,
the map $(t,\omega)\mapsto\int\phi(g)\mathsf K_t(\omega,dg)$ is measurable
with respect to the predictable sigma-field of $\mathbb F$.  An exact simulator
may instead pre-draw the exponential interarrivals; the enlarged filtration
then knows the next jump time, so the compensator identity is used only in
$\mathbb F$.  Lemma~\ref{lem:clock-law-equivalence} in
Appendix~\ref{app:generator-moments} shows that both constructions have the
same law provided every event parameter (graph, depth, batch size, operator) is
chosen from the revealed past before the corresponding interarrival is drawn.
An adaptive scheduler may therefore depend on all past states, gradients,
roots, and graph products, but not on the just-drawn interarrival.

At a jump the oracle value is $g_t\in\mathbb R^d$, queried at $x_{t-}$.

\begin{assumption}[Conditional oracle moments]\label{ass:central-oracle-moments}
There are predictable vectors $e_t\in\mathbb R^d$ and predictable scalars
$\varsigma_t\geq0$ such that
\begin{equation}\label{eq:central-oracle-moments}
\mathbb E[g_t\mid\mathcal F_{t-}]=\nabla F(x_{t-})+e_t,
\qquad
\mathbb E\!\left[
 \|g_t-\mathbb E[g_t\mid\mathcal F_{t-}]\|^2
 \mid\mathcal F_{t-}\right]\leq\varsigma_t^2.
\end{equation}
\end{assumption}

CQA-Net and the fixed-graph token have $e_t=0$; the changing-graph token,
and D-CQA have a predictable bias.

\subsection{Communication chronology and cost metrics}\label{sec:communication-model}

A communication subround exchanges one $d$-vector along every active edge.  At
jump $k\geq1$, $W_{k,q}\in\mathbb R^{m\times m}$ is the input-output operator of
the depth-$q$ mean-preserving protocol (Section~\ref{sec:network-oracle-design}).
The graph and protocol parameters for jump $k$ are fixed from the revealed past
before $E_{k-1}$ is exposed.  At $\tau_k$ the common pre-jump query is formed.
If $q>0$, fresh rows are generated and mixed first, and only then is a
conditionally uniform root, independent of their innovations, revealed.  If
$q=0$, the root may be revealed first, but only through the lazy
coupling~\eqref{eq:lazy-root-row-kernel}, which has the same law and avoids all
unelected gradient evaluations.  The selected row is then broadcast.
Pre-shared clock and root streams implement this chronology for an oblivious
schedule; an adaptive deployment needs separate primitives for clock and root
generation, whose cost is not included in the round ledger.

We report expected events, component-gradient evaluations per node and in
aggregate, critical-path neighbor rounds, and transmitted $d$-vectors.  On a
fixed undirected graph, $D$ is a broadcast depth, $M_{\rm br}$ bounds its
transmitted vectors, $\chi:=\lambda_m/\lambda_2$ is the ratio of the largest to
the smallest nonzero Laplacian eigenvalue, and one Laplacian multiplication
sends at most $2|E|$ directed vectors.  Changing-graph costs are expressed
through certified mixing and dissemination bounds.

\section{Algorithm design}\label{sec:algorithm-design}

\begin{figure}[t]
  \centering
  \includegraphics[width=.96\textwidth]{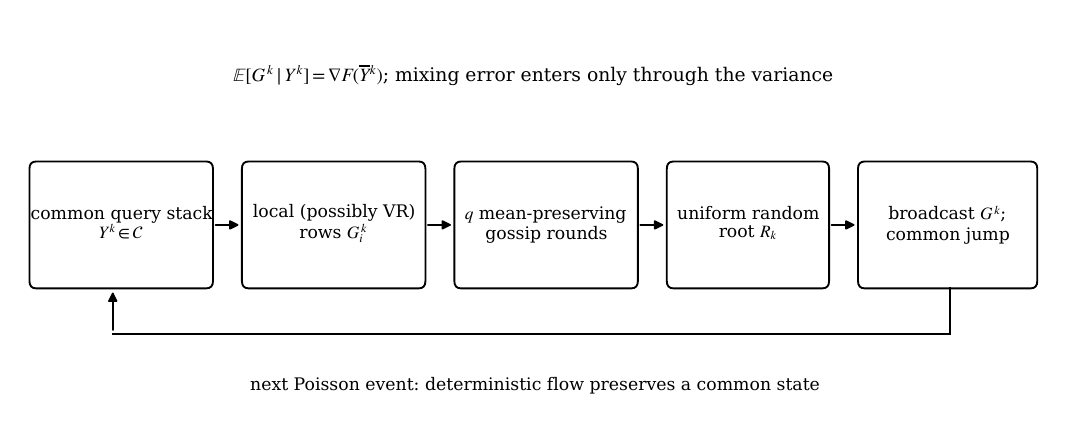}
  \caption{One event of CQA-Net.  All nodes share the query, form local rows,
  run a finite mean-preserving protocol, reveal one row conditionally uniformly,
  and use its broadcast value in the common accelerated jump.  At $q=0$, the
  law-equivalent lazy coupling evaluates only the selected row, using a family
  of conditional row kernels fixed independently of the root.}
  \label{fig:algorithm-pipeline}
\end{figure}

The method has three layers, each with a narrow input-output interface.
Figure~\ref{fig:algorithm-pipeline} shows one network event.  This section
describes how to run the algorithms; moment and convergence estimates are in
Section~\ref{sec:main-results}.

The layers are the following.
\begin{enumerate}
\item The exact flow and Poisson jumps in Algorithm~\ref{alg:base-process}
produce one oracle vector per event; Theorem~\ref{thm:master-inexact} analyzes
this one-vector interface.
\item Algorithm~\ref{alg:network-oracle} combines local rows, depth-$q$
mixing, a uniform root, and broadcast into one common vector targeted at the
average row; Lemma~\ref{lem:uniform-row} establishes the required moments.
\item Algorithm~\ref{alg:network-composition} substitutes the network oracle
into the flow and therefore gives every node the same trajectory $x_t$;
Theorem~\ref{thm:static-network} treats static and changing graphs.
\end{enumerate}
Layers~2 and~3 use three global primitives, which
Section~\ref{sec:global-primitives} makes explicit.
Sections~\ref{sec:token-design} and~\ref{sec:dcqa-design} replace these layers
by a walking token (Algorithms~\ref{alg:token-oracle}
and~\ref{alg:token-cqa}) and by local states (Algorithm~\ref{alg:dcqa}); the
first layer is unchanged.

\subsection{The basic continuized process}\label{sec:base-process-design}

The first layer replaces the auxiliary low-dimensional searches used by early
accelerated quasar-convex methods with an exactly solvable event-driven flow.
The point $x_t$ is the point at which the objective is controlled, while
$z_t$ is a momentum-like auxiliary state.  Between events, $z_t$ is fixed and
$x_t$ moves toward it; at an event, both states use the same oracle vector with
different weights.  The coefficients solve the two cancellation requirements of
the quasar Lyapunov balance.  The deterministic undamped specialization is the
continuized quasar dynamics of \citet{wang2023}; the parameter
$\theta\in(0,1]$ leaves a tunable dissipative reserve when $\theta<1$.

Starting from $x_0=z_0\in\mathbb R^d$, define
\begin{equation}\label{eq:base-process}
 d x_t=\frac{2}{\gamma t}(z_t-x_t)\,dt-\frac1L g_t\,dn_t,
 \qquad
 d z_t=-\frac{\gamma\theta^2t}{2L}g_t\,dn_t.
\end{equation}
The oracle is queried at $x_{t-}$.  If no event occurs on $(r,t]$, then
\begin{equation}\label{eq:base-flow}
 z_t=z_r,
 \qquad
 x_t=z_r+\left(\frac r t\right)^{2/\gamma}(x_r-z_r).
\end{equation}
The singularity at the origin is removable because $x_0=z_0$: before the first
event the state remains fixed.  Because of~\eqref{eq:base-flow}, the process can
be simulated without a line search.

For a predictable biased oracle one may project the auxiliary jump onto a
closed convex set $\mathcal Q\ni\{x_0,x^\star\}$.  This projected variant is
used only by the projected bound of Theorem~\ref{thm:master-inexact}(b) and its
application to the changing-graph token; CQA-Net, the fixed-graph token, and D-CQA never
project.

\begin{algorithm}[t]
\caption{CQA: exact simulation of the one-vector continuized process}
\label{alg:base-process}
\begin{algorithmic}[1]
\Require horizon $T>0$, $L>0$, $\gamma\in(0,1]$, $\theta\in(0,1]$, initial point $x^{(0)}$, oracle, optional $\mathcal Q$
\State $\tau_0\gets0$, $z^{(0)}\gets x^{(0)}$, $k\gets0$
\While{true}
  \State draw $E_k\sim\mathop{\rm Exp}\nolimits(1)$ and set $\tau_{k+1}\gets\tau_k+E_k$
  \If{$\tau_{k+1}>T$}
    \State \Return $z^{(k)}+(\tau_k/T)^{2/\gamma}(x^{(k)}-z^{(k)})$, with $(0/T)^{2/\gamma}:=0$
  \EndIf
  \State $y^{(k+1)}\gets z^{(k)}+(\tau_k/\tau_{k+1})^{2/\gamma}(x^{(k)}-z^{(k)})$
  \State $g^{(k+1)}\gets\mathop{\rm Oracle}\nolimits(y^{(k+1)})$
  \State $x^{(k+1)}\gets y^{(k+1)}-g^{(k+1)}/L$
  \State $z^{(k+1)}\gets z^{(k)}-\gamma\theta^2\tau_{k+1}g^{(k+1)}/(2L)$
  \If{the biased localized variant is used}
    \State $z^{(k+1)}\gets\Pi_{\mathcal Q}(z^{(k+1)})$
  \EndIf
  \State $k\gets k+1$
\EndWhile
\end{algorithmic}
\end{algorithm}

Algorithm~\ref{alg:base-process} is an exact simulator of
\eqref{eq:base-process}, not a time discretization.  The Poisson clock is chosen
because its compensator turns the expected weighted jump sum into the time
integral of the continuous drift; this alignment permits the cancellation proof
and gives $\mathbb E n_T=T$.  We do not claim physical asynchrony, and we do not
know whether the random clock is necessary; a deterministic grid would need a
separate discrete estimate-sequence analysis.

\subsection{Distributed gradient estimation}\label{sec:network-oracle-design}

A depth-$q$ mean-preserving mixer is a linear distributed protocol that is
implementable in $q$ neighbor rounds, may use finite local memory, and has an
input-output map $U\mapsto W_{k,q}U$ satisfying
$W_{k,q}\mathbf1=\mathbf1$ and
$\mathbf1^\top W_{k,q}=\mathbf1^\top$; set $W_{k,0}=I_m$.  Put
$J:=m^{-1}\mathbf1\mathbf1^\top$, the exact-averaging operator.  Mean
preservation is equivalent to $W_{k,q}J=JW_{k,q}=J$, so $W_{k,q}-J$ acts only
on the disagreement subspace.  The definition covers ordinary products of
one-round matrices and polynomial filters.

For a fixed connected undirected graph, a Chebyshev filter can be executed by a
three-term recurrence.  Let $\mathsf L$ be the Laplacian, let
$0<\underline\lambda\leq\lambda_2(\mathsf L)$ and
$\lambda_m(\mathsf L)\leq\overline\lambda$, and put
$c=(\overline\lambda+\underline\lambda)/2$,
$r=(\overline\lambda-\underline\lambda)/2$,
$\mathsf R=(cI-\mathsf L)/r$, and $\xi=c/r$.  For an input row matrix $U$, set
$V_0=U$, $V_1=\mathsf R U$, $a_0=1$, $a_1=\xi$, and iterate
\[
 V_{\ell+1}=2\mathsf R V_\ell-V_{\ell-1},
 \qquad
 a_{\ell+1}=2\xi a_\ell-a_{\ell-1}.
\]
The output is $W_qU=V_q/a_q$.  Each recurrence step uses one Laplacian
multiplication and stores two matrices locally.  When the certified interval
collapses, the one-round polynomial $I-\mathsf L/\lambda_2$ gives exact
averaging.  The spectral contraction supplied by this implementation is analyzed
after the oracle has been constructed.

Algorithm~\ref{alg:network-oracle} turns arbitrary local rows into one common
vector.  At positive depth, the ordering ``protocol and rows first, root
second'' prevents the selected index from influencing the row matrix whose
average is to be estimated.  At depth zero, a root-first implementation is
permitted only through the root-independent lazy conditional kernel
in~\eqref{eq:lazy-root-row-kernel}; it is distributionally identical to
selecting from a latent row family fixed independently of the root.

\begin{algorithm}[t]
\caption{Finite-mixing uniform-row network oracle}
\label{alg:network-oracle}
\begin{algorithmic}[1]
\Require jump index $k\geq1$, common query $y^{(k)}$, local row kernels,
batch size $b\in\mathbb N$ with $b\geq1$, mean-preserving $W_{k,q}$
\If{$q=0$}
  \State reveal $R_k$ conditionally uniformly, independently of the latent row innovations
  \State only node $R_k$ samples $\mathsf K_{R_k,k}$ using a fresh seed independent of $R_k$
  \State node $R_k$ computes $g_{R_k}^{(k)}(y^{(k)})$
  \State $\widehat g^{(k)}\gets g_{R_k}^{(k)}(y^{(k)})$
\Else
  \ForAll{nodes $i$ in parallel}
    \State compute a fresh row $g_i^{(k)}(y^{(k)})$
  \EndFor
  \State stack the rows in $G^{(k)}$ and run the $q$-round protocol to form $W_{k,q}G^{(k)}$
  \State reveal $R_k$ conditionally uniformly after $W_{k,q}$ and all fresh rows are fixed
  \State $\widehat g^{(k)}\gets e_{R_k}^\top W_{k,q}G^{(k)}$
\EndIf
\State node $R_k$ broadcasts $\widehat g^{(k)}$; every node returns that vector
\end{algorithmic}
\end{algorithm}

At $q=0$ this is uniform node sampling followed by broadcast.  Equation
\eqref{eq:lazy-root-row-kernel} is essential: marginal row unbiasedness alone
would not suffice if the post-selection row law could depend on the root.
Positive depth is not needed for unbiasedness.  It can reduce the selected-row
variance when the certified norm profile of $W_{k,q}$ contracts; mean
preservation by itself does not imply such a contraction.  The price is one
fresh row at every node and the additional neighbor rounds.

\subsection{Composition: CQA-Net}\label{sec:composition-design}

CQA requires one vector per event, and Algorithm~\ref{alg:network-oracle}
returns exactly one common vector.  Their composition therefore needs no new
state variable.  Every node stores identical copies of $x^{(k)}$ and
$z^{(k)}$; after a common flow and a broadcast jump, the copies remain
identical.  We call the resulting method \emph{CQA-Net}.

\begin{algorithm}[t]
\caption{CQA-Net: synchronized common-state network acceleration}
\label{alg:network-composition}
\begin{algorithmic}[1]
\Require the inputs of Algorithm~\ref{alg:base-process}, depth $q$, batch $b$, and an admissible event-indexed network schedule
\State every node keeps a copy of $(x,z,\tau)$, initialized to $(x^{(0)},x^{(0)},0)$, and runs Algorithm~\ref{alg:base-process} on it with two changes:
\State \quad before $E_k$ is drawn, the graph and $W_{k+1,q}$ are chosen from the revealed history, and $E_k$ is common to all nodes
\State \quad line~8 is replaced by Algorithm~\ref{alg:network-oracle} at the common query $y^{(k+1)}$; every node uses the broadcast $\widehat g^{(k+1)}$
\end{algorithmic}
\end{algorithm}

Because all rows of the optimization state are identical, the objective is
always evaluated at the common vector $x_t\in\mathbb R^d$; no matrix-valued
extension of $F$ is introduced.  The network matrix acts only on the transient
stack of local oracle rows.  This reduction works only because the last line
of Algorithm~\ref{alg:network-oracle} delivers the same vector to every node;
Section~\ref{sec:global-primitives} discusses the cost of that step and
Sections~\ref{sec:token-design}--\ref{sec:dcqa-design} remove it.

\subsection{Global primitives of CQA-Net}\label{sec:global-primitives}

Algorithms~\ref{alg:network-oracle} and~\ref{alg:network-composition} use
three operations that involve the whole network at every event: (G1)~exact
dissemination of $\widehat g^{(k)}$ from $R_k$ to every node (line~14 of
Algorithm~\ref{alg:network-oracle}); (G2)~a shared event clock; and (G3)~a
shared conditionally uniform root.  (G1) is essential: the common-state
invariant, and hence the reduction to Theorem~\ref{thm:master-inexact}, needs
every node to receive exactly the same vector.  Approximate averaging cannot
replace it, since finitely many gossip rounds give close but unequal vectors
and the copies of $(x,z)$ drift apart.  In federated learning (G1) is a server
broadcast.  In a peer-to-peer network it costs $D$ rounds per event, which is
exactly the term that separates the round count of CQA-Net from the
local-state scale (Section~\ref{sec:network-results}).  On a changing
graph Theorem~\ref{thm:static-network} takes a dissemination bound $B_{\rm br}$
as an input, and such a bound exists only if the schedule is known in advance or
every window of rounds is connected.

We call a method \emph{fully decentralized} if (R1)~every message travels along
one edge of the current graph; (R2)~all nodes run the same code, with no fixed or
moving coordinator; (R3)~every node keeps its own copy of the optimization
variables, brought together only by neighbor exchanges; and (R4)~nodes never
wait for an operation to finish over the whole network.
Sections~\ref{sec:token-design} and~\ref{sec:dcqa-design} give two
implementations without (G1); Table~\ref{tab:decentralization} lists which
requirements each design meets and what one event costs.  The token removes all
three global primitives but keeps a single mobile state.  D-CQA has local states
and keeps only a shared clock seed and synchronous rounds, as do synchronous
gossip methods such as decentralized gradient descent and gradient tracking.

\begin{table}[!htbp]
\caption{Decentralization requirements and the cost of one event in the three
implementations.}
\label{tab:decentralization}
\centering
\small
\begin{tabular}{@{}P{2.6cm}P{3.6cm}P{4.6cm}P{3.6cm}@{}}
\toprule
 & CQA-Net (Alg.~\ref{alg:network-composition}) & Token (Alg.~\ref{alg:token-cqa}) & D-CQA (Alg.~\ref{alg:dcqa}) \\
\midrule
(R1) local communication & no (dissemination) & yes & yes \\
(R2) no distinguished node & no (root $R_k$) & almost (current holder) & yes \\
(R3) local states & identical copies only & no (one mobile state) & yes \\
(R4) no per-event coordination & no (barrier, shared clock and root) & yes & almost (shared clock seed) \\
\midrule
Critical-path rounds & $q+D$ (fixed), $q+B_{\rm br}$ (changing) & $h_{\rm w}+2q$, sequential & $q$ \\
Row evaluations & $m$ if $q>0$, one if $q=0$ & $|B_q(R_k)|$; one if $q=0$ & $m$ \\
Transmitted data & $2q|E|+M_{\rm br}$ vectors in $\mathbb R^d$ & at most $h_{\rm w}$ token messages of $2d+2$ scalars, plus at most $2(q+1)|E(B_q(R_k))|$ vectors in $\mathbb R^d$ & $4|E|q$ vectors in $\mathbb R^d$ \\
One-time setup & none & fixed mode: node $i$ stores row $i$ of $P^{h_{\rm w}}$ & none \\
\bottomrule
\end{tabular}
\end{table}

\subsection{A walking token}\label{sec:token-design}

In the token variant the state $(x,z,\tau)$ is stored in a single message, the
\emph{token}, that moves along edges.  The node holding it is the
\emph{holder}.  The holder draws the next interarrival itself, so no shared
clock is needed.  The root is the node at which a random walk of the token
stops, so no shared root is needed.  The oracle vector is used where it is
computed, so nothing is disseminated.  Algorithm~\ref{alg:base-process} is
executed unchanged; only the place where its lines run moves through the
network.

\paragraph{Walk matrices.}
At event $k$ the token makes $h_{\rm w}\geq0$ hops with matrices
$P_{k,1},\ldots,P_{k,h_{\rm w}}$.  Each of them is entrywise nonnegative,
doubly stochastic, and supported on the self-loops and the active edges of its
hop; Metropolis weights of the current graph are an example.  The sequence is
fixed from the revealed event history before the interarrival that produces
$\tau_k$ is drawn, and it does not depend on the positions visited by the walk
itself.  We write $o_k$ for the holder at the start of event $k$, so
$o_k=R_{k-1}$ for $k\geq2$ and $o_1$ is the initial holder.

\paragraph{Local mixer.}
For $q\geq0$ the token oracle uses a mean-preserving, $\mathbb F$-predictable
$W_{k,q}$ with the \emph{ball-support property}
$(W_{k,q})_{vi}=0$ whenever $\dist(v,i)>q$; we set $W_{k,0}=I_m$.  On a fixed
graph every degree-$q$ Laplacian polynomial has this property, in particular
the normalized Chebyshev filter of Section~\ref{sec:network-oracle-design}.
Then the selected row $e_v^\top W_{k,q}G^{(k)}$ depends only on the rows of the
ball $B_q(v)=\{i:\dist(v,i)\leq q\}$.  To evaluate it, $v$ floods the query to
$B_q(v)$, every node of the ball forms its row, and the three-term recurrence
is run only where it is needed: step $\ell$ of the recurrence is computed at the
nodes with $\dist(\cdot,v)\leq q-\ell$.  This costs $2q$ subrounds.  On a
changing graph we use $q=0$, because a ball-local evaluation of a
time-varying product is not considered here.

\paragraph{Two modes.}
In the \emph{fixed} mode the graph and the walk matrix
$P_{k,s}=P=P^\top$ are fixed, and the oracle is reweighted by
$\omega_k=1/(m(P^{h_{\rm w}})_{o_kR_k})$, the inverse of $m$ times the
probability of the realized endpoint.  Node $i$ stores row~$i$ of
$P^{h_{\rm w}}$, which a one-time setup computes by $h_{\rm w}$ rounds of
gossip with messages of $m$ scalars.  In the \emph{changing} mode no weight is
used, $\omega_k=1$.  Optionally, the auxiliary jump is projected onto a closed
bounded convex set $\mathcal Q\ni x^{(0)},x^\star$, as in the localized variant
of Algorithm~\ref{alg:base-process}; the projection is needed only for the
second of the two bounds in Corollary~\ref{cor:token-changing}.  The fixed mode
needs $\theta<1$, because the walk creates multiplicative noise that the
reserve must absorb.

\begin{algorithm}[t]
\caption{Walking-token network oracle (replaces Algorithm~\ref{alg:network-oracle})}
\label{alg:token-oracle}
\begin{algorithmic}[1]
\Require event index $k\geq1$; token $(y^{(k)},z,\tau_k,o_k)$ at node $o_k$;
hop count $h_{\rm w}$; walk matrices $P_{k,1},\ldots,P_{k,h_{\rm w}}$; depth $q$
and $W_{k,q}$ as above; batch $b\geq1$; mode $\in\{\text{fixed},\text{changing}\}$
\State $v\gets o_k$
\For{$s=1,\ldots,h_{\rm w}$}
  \State node $v$ draws $j$ with probability $(P_{k,s})_{vj}$, independently of all row innovations and interarrivals
  \State if $j\neq v$, node $v$ sends the token to its neighbor $j$; \ $v\gets j$
\EndFor
\State $R_k\gets v$ \Comment{replaces the uniform root of Algorithm~\ref{alg:network-oracle}}
\If{$q=0$}
  \State node $R_k$ samples $\widetilde g^{(k)}\sim\mathsf K_{R_k,k}$ with a fresh seed
\Else
  \State $R_k$ floods $y^{(k)}$ to $B_q(R_k)$ \Comment{$q$ subrounds}
  \State every $i\in B_q(R_k)$ samples a fresh row $g_i^{(k)}(y^{(k)})\sim\mathsf K_{i,k}$
  \State run the recurrence of Section~\ref{sec:network-oracle-design}; step $\ell$ only at $\dist(\cdot,R_k)\leq q-\ell$ \Comment{$q$ subrounds}
  \State $\widetilde g^{(k)}\gets e_{R_k}^\top W_{k,q}G^{(k)}$ at node $R_k$
\EndIf
\If{mode $=$ fixed}
  \State $\omega_k\gets1/\bigl(m(P^{h_{\rm w}})_{o_kR_k}\bigr)$ \Comment{read from the stored row of $P^{h_{\rm w}}$}
\Else
  \State $\omega_k\gets1$
\EndIf
\State $\widehat g^{(k)}\gets\omega_k\widetilde g^{(k)}$; the vector stays at node $R_k$ \Comment{no dissemination}
\end{algorithmic}
\end{algorithm}

\begin{algorithm}[t]
\caption{CQA with a walking token (replaces Algorithm~\ref{alg:network-composition})}
\label{alg:token-cqa}
\begin{algorithmic}[1]
\Require the inputs of Algorithm~\ref{alg:base-process} with $\theta<1$ in fixed mode; $h_{\rm w}$, $q$, $b$; initial holder $o_1$; optionally, in changing mode, a closed bounded convex $\mathcal Q\ni x^{(0)},x^\star$
\State the token stores $(x,z,\tau)=(x^{(0)},x^{(0)},0)$ and the identity of its holder; the holder runs Algorithm~\ref{alg:base-process} with three changes:
\State \quad $h_{\rm w}$, $q$, $W_{k+1,q}$ and the walk matrices are fixed from the revealed history, and the holder alone draws $E_k$ \Comment{no shared clock}
\State \quad line~8 is replaced by Algorithm~\ref{alg:token-oracle}, which moves the token to $R_{k+1}$; lines~9--10 are executed there
\State \quad if $\mathcal Q$ is given, line~12 projects $z$ onto $\mathcal Q$
\end{algorithmic}
\end{algorithm}

\paragraph{Chronology, storage, and cost.}
The walk runs after the query is formed and before the row is sampled, so its
endpoint is part of the event mark, exactly as $R_k$ in
Section~\ref{sec:communication-model}; the walk randomness is independent of
the interarrivals, so Lemma~\ref{lem:clock-law-equivalence} applies verbatim.
The token carries $(x,z,\tau)$ between events and $(y,z,\tau',o)$ during the
walk, that is, $2d+2$ scalars.  Table~\ref{tab:decentralization} gives the cost
of one event: $h_{\rm w}$ sequential hops replace the $D$ rounds of
dissemination, and Remark~\ref{rem:token-length} shows that $h_{\rm w}$ is of
order $t_{\rm rel}\log m$, where $t_{\rm rel}=(1-\lambda_P)^{-1}$.

\subsection{Local states: D-CQA}\label{sec:dcqa-design}

The second implementation removes the common state.  Every node keeps its own
pair $(x_i,z_i)$, runs the closed-form flow, takes a local step with its own
row at its own query, and then $q$ rounds of a mean-preserving protocol mix the
\emph{states}, not the rows.  There is no root, no dissemination and no token.
We write $X,Z\in\mathbb R^{m\times d}$ for the stacked states, with rows
$x_i,z_i$, and $\bar x=m^{-1}X^\top\mathbf1$, $\bar z=m^{-1}Z^\top\mathbf1$.
The disagreement is $X_\perp:=(I-J)X$.  Between events the flowed states are
$X(t)=Z+(\tau_k/t)^{2/\gamma}(X-Z)$.

\begin{algorithm}[t]
\caption{D-CQA: decentralized CQA with local states (replaces Algorithms~\ref{alg:network-oracle} and~\ref{alg:network-composition})}
\label{alg:dcqa}
\begin{algorithmic}[1]
\Require $T,L>0$; $\gamma,\theta\in(0,1]$; batch $b\geq1$; common $x^{(0)}$; an admissible sequence of mean-preserving $q$-round mixers $W_k$
\State every node sets $x_i\gets x^{(0)}$, $z_i\gets x^{(0)}$; \ $\tau_0\gets0$, $k\gets0$
\While{true}
  \State fix $W_{k+1}$ from the revealed history \Comment{before $E_k$}
  \State every node computes the common $E_k\sim\mathop{\rm Exp}\nolimits(1)$ from a shared seed; \ $\tau_{k+1}\gets\tau_k+E_k$
  \If{$\tau_{k+1}>T$}
    \State every node returns $z_i+(\tau_k/T)^{2/\gamma}(x_i-z_i)$, with $(0/T)^{2/\gamma}:=0$; \textbf{stop}
  \EndIf
  \ForAll{nodes $i$ in parallel}
    \State $y_i\gets z_i+(\tau_k/\tau_{k+1})^{2/\gamma}(x_i-z_i)$
    \State $g_i\gets$ fresh local row at $y_i$ \Comment{own query, own data}
    \State $\widetilde x_i\gets y_i-g_i/L$, \quad $\widetilde z_i\gets z_i-\gamma\theta^2\tau_{k+1}g_i/(2L)$
  \EndFor
  \State run the $q$-round protocol on both stacks: $X\gets W_{k+1}\widetilde X$, \ $Z\gets W_{k+1}\widetilde Z$ \Comment{each message: two vectors in $\mathbb R^d$}
  \State $k\gets k+1$
\EndWhile
\end{algorithmic}
\end{algorithm}

Relative to Algorithms~\ref{alg:network-oracle}--\ref{alg:network-composition},
the rows are no longer mixed, the uniform root and the dissemination step
disappear, the queries $y_i$ differ across nodes, and every node returns its
own point.  The only global element left is the shared clock seed: every node
generates the same event times locally, without communication, and rounds are
synchronous.  Each event costs one row per node, $q$ rounds and $4|E|q$
transmitted vectors in $\mathbb R^d$.  Because the queries differ,
quasar-convexity can no longer be applied at the point where the gradients are
evaluated.  Section~\ref{sec:dcqa-results} therefore analyzes the network
average, for which the disagreement acts as a bias.

\section{Main results}\label{sec:main-results}

Two results carry the analysis.  Theorem~\ref{thm:master-inexact} bounds the continuized process for any oracle
with a predictable bias and a state-dependent variance.
Lemma~\ref{lem:uniform-row} computes the conditional moments of every
network oracle that returns one, possibly reweighted, row of a mixed matrix.
Every other statement substitutes moments into these two results, certifies the
network operators, or closes the resulting oracle inequality.  For CQA-Net,
Lemma~\ref{lem:uniform-row}(ii) and Theorem~\ref{thm:master-inexact}(a) give
Theorem~\ref{thm:static-network}, which Corollary~\ref{cor:closed-static}
closes, and Lemma~\ref{lem:chebyshev-filter} and
Corollary~\ref{cor:block-contractive-changing} certify the operators.  For the
token, Lemma~\ref{lem:walk-endpoint} and Lemma~\ref{lem:uniform-row}(ii)--(iii)
give Corollaries~\ref{cor:token-fixed} and~\ref{cor:token-changing}.  For D-CQA,
Lemmas~\ref{lem:dcqa-average} and~\ref{lem:dcqa-disagreement} and
Theorem~\ref{thm:master-inexact}(a) give Theorem~\ref{thm:dcqa}.  The quantities
to be certified are therefore the oracle bias and variance
envelope, the operator norms $\|W_{k,q}-J\|_2$ and
$m^{-1}\|W_{k,q}\|_{\mathrm F}^2$, and the walk contraction.

\subsection{The basic process and the one-vector oracle}\label{sec:base-results}

For $t>0$, put
\begin{equation}\label{eq:lyapunov-weights}
A_t=\frac{\gamma^2\theta^2t^2}{4L},
\qquad
c_t=\frac{\gamma\theta^2t}{2L},
\qquad
V_\theta=\frac{1+\theta^2}{2},
\end{equation}
and define
\begin{equation}\label{eq:master-potential}
\Phi_t=A_t\bigl(F(x_t)-F(x^\star)\bigr)
 +\frac12\|z_t-x^\star\|^2.
\end{equation}
The weights solve the two cancellation equations
$A_t'=\gamma c_t$ and $c_t^2=\theta^2A_t/L$.  Together with
$2A_t/(\gamma t)=c_t$ and $A_t/A_T=t^2/T^2$, these identities follow by direct
substitution.  Within the nonzero positive quadratic--linear family, the
linear equation alone leaves an overall scale free, but the quadratic equation
fixes that scale; with the process normalization above the two equations
determine the displayed weights exactly.  (The identically zero pair is the
irrelevant degenerate solution.)

For a predictable mark kernel $\mathsf K_t$ from
Section~\ref{sec:oracle-model}, define the predictable extended drift operator
on a smooth test function $\Psi(t,x,z)$ by
\begin{equation}\label{eq:predictable-drift-operator}
\begin{aligned}
\mathscr D_t\Psi(t,x,z)
={}&\partial_t\Psi(t,x,z)
 +\left\langle\nabla_x\Psi(t,x,z),\frac{2}{\gamma t}(z-x)\right\rangle\\
&+\int\!\bigl[\Psi(t,x-g/L,z-c_tg)-\Psi(t,x,z)\bigr]
  \mathsf K_t(dg),
\end{aligned}
\end{equation}
with $z-c_tg$ replaced by its projection in the localized variant.  Under
Markovian rules this is the usual time-inhomogeneous PDMP generator; under the
history-dependent model it is a predictable drift characteristic, and the
proof uses the It\^o compensator identity rather than the Markov property.

\begin{lemma}[Predictable Lyapunov drift]\label{lem:predictable-drift}
Suppose that $F$ is $L$-smooth and $\gamma$-quasar-convex, and that
Assumption~\ref{ass:central-oracle-moments} holds.  For the unprojected process,
\begin{align}
\mathscr D_t\Phi_t
\leq{}&A_t'\bigl(F(x_{t-})-F(x^\star)\bigr)
 -c_t\langle\nabla F(x_{t-}),x_{t-}-x^\star\rangle\notag\\
&-\frac{A_t(1-\theta^2)}{2L}\|\nabla F(x_{t-})\|^2
 +\frac{A_tV_\theta}{L}\varsigma_t^2\notag\\
&-c_t\langle z_{t-}-x^\star,e_t\rangle
 +\frac{A_t\theta^2}{L}\langle\nabla F(x_{t-}),e_t\rangle
 +\frac{A_tV_\theta}{L}\|e_t\|^2.
\label{eq:master-drift}
\end{align}
The first line is nonpositive because $A_t'=\gamma c_t$ and
\eqref{eq:quasar}.  The third term is the dissipative reserve.  If the
auxiliary jump is projected onto a bounded convex
$\mathcal Q\ni\{x_0,x^\star\}$, the term
$-c_t\langle z_{t-}-x^\star,e_t\rangle$ is at most
$c_tR_{\mathcal Q}\|e_t\|$, where
$R_{\mathcal Q}=\sup_{z\in\mathcal Q}\|z-x^\star\|$.
\end{lemma}

\begin{proof}
The deterministic drift and the two jump increments are computed before taking
conditional moments.  In the last part, the weights
\eqref{eq:lyapunov-weights} remove every occurrence of the auxiliary state and
leave the quasar inequality at the single point $x_{t-}$.

\smallskip\noindent\textit{Step 1: drift between jumps.}\quad
Between events, $z_t$ is constant and
$\dot x_t=2(z_t-x_t)/(\gamma t)$.  The time derivative of
\eqref{eq:master-potential} and the drift of $x_t$ give
\[
A_t'\bigl(F(x_{t-})-F(x^\star)\bigr)
 +\frac{2A_t}{\gamma t}
  \langle\nabla F(x_{t-}),z_{t-}-x_{t-}\rangle.
\]
The identity $2A_t/(\gamma t)=c_t$ converts the second coefficient to the
same $c_t$ that appears in the auxiliary jump.

\smallskip\noindent\textit{Step 2: jump contributions.}\quad
Smoothness gives
\[
F(x_{t-}-g_t/L)-F(x_{t-})
\leq-\frac1L\langle\nabla F(x_{t-}),g_t\rangle
 +\frac1{2L}\|g_t\|^2.
\]
The auxiliary distance changes exactly by
\[
\frac12\|z_{t-}-c_tg_t-x^\star\|^2
 -\frac12\|z_{t-}-x^\star\|^2
=-c_t\langle z_{t-}-x^\star,g_t\rangle
 +\frac{c_t^2}{2}\|g_t\|^2.
\]
Projection cannot increase this distance because $x^\star\in\mathcal Q$.

\smallskip\noindent\textit{Step 3: conditional moments and cancellation.}\quad
The conditional moment relations~\eqref{eq:central-oracle-moments} imply
\[
\mathbb E[g_t\mid\mathcal F_{t-}]=\nabla F(x_{t-})+e_t,
\qquad
\mathbb E[\|g_t\|^2\mid\mathcal F_{t-}]
\leq\|\nabla F(x_{t-})+e_t\|^2+\varsigma_t^2.
\]
The term $A_t'[F(x_{t-})-F(x^\star)]$ pairs with the quasar inner product
through $A_t'=\gamma c_t$; the $z_{t-}$ part of the drift term cancels the
linear auxiliary jump and leaves $-c_t\langle\nabla F(x_{t-}),x_{t-}-x^\star\rangle$;
the linear part of the $x_t$ jump, $-A_t\|\nabla F(x_{t-})\|^2/L$, is partly
cancelled by the two jump squares, $A_t(1+\theta^2)\|\nabla F(x_{t-})\|^2/(2L)$,
which leaves the reserve; and the centered fluctuation gives
$A_tV_\theta\varsigma_t^2/L$.  Here $c_t^2=\theta^2A_t/L$ and $V_\theta=(1+\theta^2)/2$ combine the quadratic
jump terms.  Expanding the mean as $\nabla F(x_{t-})+e_t$ gives the three bias
terms in the last two lines of~\eqref{eq:master-drift}.
\end{proof}

The next theorem integrates Lemma~\ref{lem:predictable-drift}.  It is the only
place where the Poisson clock, rather than merely the event-driven algorithm,
enters the convergence proof, and every rate in this paper is obtained from it
by substituting conditional oracle moments.  Two elementary integral
inequalities are used to close the estimates.

\begin{lemma}[Linear and square-root Gronwall inequalities]\label{lem:gronwall}
Let $u:[0,T]\to[0,\infty)$ be bounded and measurable, let $a\geq0$, and let
$p,h\geq0$ be integrable on $[0,T]$.
\begin{enumerate}[label=(\roman*),leftmargin=2.2em]
\item If $\lambda\geq0$ and $u(t)\leq a+\int_0^t[\lambda u(r)+h(r)]dr$ for all
$t$, then
\begin{equation}\label{eq:conditional-gronwall-main}
u(t)\leq e^{\lambda t}
\left[a+\int_0^te^{-\lambda r}h(r)dr\right]
\leq e^{\lambda t}\left[a+\int_0^th(r)dr\right].
\end{equation}
\item If $a>0$ and $u(t)\leq a+\int_0^t\bigl[p(r)\sqrt{u(r)}+h(r)\bigr]dr$
for all $t$, then
\begin{equation}\label{eq:sqrt-gronwall}
\sqrt{u(t)}\leq\sqrt{a+\int_0^th(r)\,dr}+\frac12\int_0^tp(r)\,dr .
\end{equation}
\end{enumerate}
\end{lemma}

\begin{proof}
Part~(i) is the classical Gronwall inequality.  For part~(ii), let $v(t)$ be
the right side of the assumed inequality and $H(t)=a+\int_0^th$.  Then $v$ is
absolutely continuous, $u\leq v$, $v\geq H>0$, and $v'\leq p\sqrt v+h$ almost
everywhere.  Therefore
$(\sqrt v)'=v'/(2\sqrt v)\leq p/2+h/(2\sqrt H)=p/2+(\sqrt H)'$, and
integration from $0$, where $v(0)=H(0)=a$, gives~\eqref{eq:sqrt-gronwall}.
\end{proof}

\begin{theorem}[Master continuized acceleration with an inexact oracle]
\label{thm:master-inexact}
Suppose that $F$ is $L$-smooth and $\gamma$-quasar-convex with respect to
$x^\star$, let $\theta\in(0,1]$, and let Algorithm~\ref{alg:base-process}
start from $x_0=z_0$ with $R_0=\|x_0-x^\star\|$.  Let
Assumption~\ref{ass:central-oracle-moments} hold, and suppose that for a
constant $\eta\geq0$ and a predictable $\bar\varsigma_t\geq0$,
\begin{equation}\label{eq:master-variance}
\varsigma_t^2\leq\eta\|\nabla F(x_{t-})\|^2+\bar\varsigma_t^2,
\qquad
V_\theta\eta\leq\frac{1-\theta^2}{2}.
\end{equation}
Assume the finite-horizon conditions
\begin{align*}
&\mathbb E\sup_{0\leq t\leq T}\Phi_t<\infty,
&&\mathbb E\sup_{0\leq t\leq T}(\|x_t\|+\|z_t\|)<\infty,\\
&\mathbb E\int_0^TA_t\|\nabla F(x_t)\|^2\,dt<\infty,
&&\mathbb E\int_0^TA_t\bar\varsigma_t^2\,dt<\infty.
\end{align*}
\begin{enumerate}[label=(\alph*),leftmargin=2.2em]
\item \emph{Mean-square bias, no projection.}  Let $\bar e:[0,T]\to[0,\infty)$
be deterministic, measurable, with $\int_0^Tt^2\bar e_t^2\,dt<\infty$, and
suppose $\mathbb E\|e_t\|^2\leq\bar e_t^2$ for $t\in[0,T]$.  Put
\[
\Gamma_T:=\frac{2LR_0^2}{\gamma^2\theta^2T^2}
 +\frac{V_\theta}{LT^2}\int_0^Tt^2\bigl(\mathbb E\bar\varsigma_t^2+\bar e_t^2\bigr)dt,
\qquad
B_T:=\frac{(1+\theta)\theta}{\sqrt{2L}\,T}\int_0^Tt\,\bar e_t\,dt .
\]
Then
\begin{equation}\label{eq:master-biased}
\mathbb E\bigl[F(x_T)-F(x^\star)\bigr]
 +\frac{(1-\theta^2)/2-V_\theta\eta}{LT^2}
   \int_0^Tt^2\mathbb E\|\nabla F(x_t)\|^2\,dt
\leq\bigl(\sqrt{\Gamma_T}+B_T\bigr)^2 .
\end{equation}
In particular, for an unbiased oracle ($e_t=0$, $\bar e=0$),
\begin{align}
&\mathbb E\bigl[F(x_T)-F(x^\star)\bigr]
 +\frac{(1-\theta^2)/2-V_\theta\eta}{LT^2}
   \int_0^Tt^2\mathbb E\|\nabla F(x_t)\|^2\,dt\notag\\
&\qquad\leq\frac{2LR_0^2}{\gamma^2\theta^2T^2}
 +\frac{V_\theta}{LT^2}\int_0^Tt^2\mathbb E\bar\varsigma_t^2\,dt.
\label{eq:master-unbiased}
\end{align}
\item \emph{Bias with projection.}  Let $\theta\in(0,1)$ and
$V_\theta\eta<(1-\theta^2)/2$, project every $z$-jump onto a closed convex
$\mathcal Q\ni x_0,x^\star$ with
$R_{\mathcal Q}:=\sup_{z\in\mathcal Q}\|z-x^\star\|<\infty$, and assume
$\mathbb E\int_0^Tt\|e_t\|\,dt<\infty$ and
$\mathbb E\int_0^Tt^2\|e_t\|^2dt<\infty$.  Put
$C_{\theta,\eta}=V_\theta+\theta^4/[2(1-\theta^2-2V_\theta\eta)]$, so that
$C_{\theta,0}=V_\theta+\theta^4/[2(1-\theta^2)]$.  Then
\begin{align}
\mathbb E[F(x_T)-F(x^\star)]
&\leq\frac{2LR_0^2}{\gamma^2\theta^2T^2}
 +\frac{2R_{\mathcal Q}}{\gamma T^2}
     \int_0^Tt\,\mathbb E\|e_t\|\,dt\notag\\
&\quad+\frac{C_{\theta,\eta}}{LT^2}
     \int_0^Tt^2\mathbb E\|e_t\|^2\,dt
 +\frac{V_\theta}{LT^2}
     \int_0^Tt^2\mathbb E\bar\varsigma_t^2\,dt.
\label{eq:biased-localized-rate}
\end{align}
\end{enumerate}
\end{theorem}

\begin{proof}[Proof of part~(a)]
\smallskip\noindent\textit{Step 1: drift.}\quad
Insert~\eqref{eq:master-variance} into~\eqref{eq:master-drift}.  The quasar
bracket is nonpositive, and the term
$A_tV_\theta\eta\|\nabla F(x_{t-})\|^2/L$ is absorbed by the dissipative
reserve, which leaves the coefficient
$\kappa_\eta:=(1-\theta^2)/2-V_\theta\eta\geq0$.  For the two linear bias
terms, $\frac12\|z_{t-}-x^\star\|^2\leq\Phi_{t-}$ gives
$-c_t\langle z_{t-}-x^\star,e_t\rangle\leq c_t\sqrt{2\Phi_{t-}}\,\|e_t\|$.
Because $x^\star$ is a global minimizer,
$\|\nabla F(x)\|^2\leq2L[F(x)-F(x^\star)]$, hence
$\|\nabla F(x_{t-})\|^2\leq2L\Phi_{t-}/A_t$; with $c_t=\theta\sqrt{A_t/L}$
this gives
$\frac{A_t\theta^2}L\langle\nabla F(x_{t-}),e_t\rangle\leq\theta c_t\sqrt{2\Phi_{t-}}\,\|e_t\|$.
Therefore
\[
\mathscr D_t\Phi_t
\leq-\frac{\kappa_\eta A_t}{L}\|\nabla F(x_{t-})\|^2
 +(1+\theta)\sqrt2\,c_t\sqrt{\Phi_{t-}}\,\|e_t\|
 +\frac{A_tV_\theta}{L}\bigl(\bar\varsigma_t^2+\|e_t\|^2\bigr).
\]
Neither a projection nor the reserve is spent on the bias.

\smallskip\noindent\textit{Step 2: localization.}\quad
Let $\tau_R=\inf\{t\geq0:\|x_t\|+\|z_t\|+\Phi_t\geq R\}\wedge T$, and put
$u_R(t)=\mathbb E\Phi_{t\wedge\tau_R}$ and
$w_R(t)=u_R(t)+\frac{\kappa_\eta}L\mathbb E\int_0^{t\wedge\tau_R}A_s\|\nabla F(x_{s-})\|^2ds$;
both are bounded by the finite-horizon conditions.  For fixed $s$,
$\Phi_{s-}=\Phi_s$ almost surely and
$\mathbf 1\{s\leq\tau_R\}\Phi_s\leq\Phi_{s\wedge\tau_R}$.  The stopped
compensator identity~\eqref{eq:ito-compensator-general} and Cauchy--Schwarz
therefore give
\[
w_R(t)\leq\Phi_0+\int_0^t\Bigl[(1+\theta)\sqrt2\,c_s\bar e_s\sqrt{u_R(s)}
 +\frac{A_sV_\theta}{L}\bigl(\mathbb E\bar\varsigma_s^2+\bar e_s^2\bigr)\Bigr]ds .
\]
Since $u_R\leq w_R$, the same inequality holds with $\sqrt{w_R}$ in place of
$\sqrt{u_R}$.

\smallskip\noindent\textit{Step 3: Gronwall and the limit $R\to\infty$.}\quad
Lemma~\ref{lem:gronwall}(ii) with $a=\Phi_0=R_0^2/2$ (if $R_0=0$, let
$a\downarrow0$) gives
\[
\sqrt{w_R(T)}\leq\sqrt{\Phi_0+\frac{V_\theta}{L}\int_0^TA_s\bigl(\mathbb E\bar\varsigma_s^2+\bar e_s^2\bigr)ds}
 +\frac{1+\theta}{\sqrt2}\int_0^Tc_s\bar e_s\,ds .
\]
Finite-horizon nonexplosion gives $\tau_R=T$ for all large $R$ almost surely.
Fatou's lemma for $\Phi_{T\wedge\tau_R}\geq0$ and monotone convergence for the
gradient integral let $R\to\infty$ on the left.

\smallskip\noindent\textit{Step 4: normalization.}\quad
The jump times are Lebesgue-null, so $x_{t-}$ may be replaced by $x_t$ inside
time integrals.  Divide by $A_T$ and use $\Phi_T\geq A_T[F(x_T)-F(x^\star)]$,
$A_t/A_T=t^2/T^2$, $\Phi_0/A_T=2LR_0^2/(\gamma^2\theta^2T^2)$, and
$A_T^{-1/2}c_t=\theta t/(\sqrt L\,T)$.  This gives~\eqref{eq:master-biased};
for $\bar e=0$, $B_T=0$ and~\eqref{eq:master-biased}
is~\eqref{eq:master-unbiased}.
\end{proof}

The proof of part~(b) spends the reserve on the bias and is given in
Appendix~\ref{app:projected-bound}.  A primitive sufficient condition for the
finite-horizon requirements is a linear conditional second-moment bound on the
oracle: Lemma~\ref{lem:finite-moments} shows that it implies all of them, and
every oracle below satisfies it whenever the local gradient second moments grow
at most linearly.

\begin{remark}[Which part is used where]\label{rem:master-parts}
Part~(a) with $e_t=0$ covers CQA-Net (Theorem~\ref{thm:static-network}) and the
fixed-graph token (Corollary~\ref{cor:token-fixed}); there $\eta=0$ or $\eta$
is the multiplicative noise created by the walk.  Part~(a) with a bias covers
the changing-graph token without projection (Corollary~\ref{cor:token-changing})
and D-CQA (Theorem~\ref{thm:dcqa}).  Part~(b) gives the alternative
projected bound for the changing-graph token.  The two parts complement each other.
For a constant bias $\bar e$ and vanishing variance, expanding the square
in~\eqref{eq:master-biased} gives the cross term $(1+\theta)R_0\bar e/\gamma$,
which has the form of the term $R_{\mathcal Q}\bar e/\gamma$
of~\eqref{eq:biased-localized-rate} with $R_0\leq R_{\mathcal Q}$; the
quadratic bias term grows like $T^2\bar e^2/L$ in part~(a) and like
$T\bar e^2/L$ in part~(b).  Part~(a) needs neither
$\mathcal Q$ nor $\theta<1$ and keeps the gradient integral on the left.  A
persistent nonzero bias leaves an accuracy floor in both parts.
\end{remark}

The exact-gradient specialization calibrates the constants and the
$\gamma$-dependence against the known continuized quasar dynamics.

\begin{corollary}[Exact-gradient acceleration]\label{cor:exact-base}
Under Theorem~\ref{thm:master-inexact}, if
$g_t=\nabla F(x_{t-})$, then
\begin{equation}\label{eq:exact-base-rate}
\mathbb E\bigl[F(x_T)-F(x^\star)\bigr]
\leq\frac{2L\|x_0-x^\star\|^2}{\gamma^2\theta^2T^2}.
\end{equation}
Estimate~\eqref{eq:exact-base-rate} shows that $T=\mathcal O(\sqrt L\|x_0-x^\star\|/(\gamma\theta\sqrt\varepsilon))$
is sufficient for expected error $\varepsilon$, and $\mathbb E n_T=T$.
\end{corollary}

At $\theta=1$, Corollary~\ref{cor:exact-base} is the deterministic
continuized result of \citet{wang2023}.  The network results
need the predictable, state-dependent variance term, the bias terms, and the
explicit reserve, which can absorb centered noise proportional to the squared
gradient norm.

\subsection{The selected-row oracle and CQA-Net}\label{sec:network-results}

Every network oracle in this paper returns one row of a mixed matrix of local
rows, possibly reweighted: CQA-Net selects a uniform row, and the walking token
selects the endpoint of a random walk.  One lemma gives the conditional
moments in all cases, even when individual row means are not local gradients.

\paragraph{Setting.}
Condition on $\mathcal F_{\tau_k-}$.  Let
$\widehat U^{(k)},U^{(k)}\in\mathbb R^{m\times d}$ have rows
$\widehat U_i^{(k)}$ and
$U_i^{(k)}=\mathbb E[\widehat U_i^{(k)}\mid\mathcal F_{\tau_k-}]$, put
$\bar u^{(k)}=m^{-1}(U^{(k)})^\top\mathbf1$ and
$D^{(k)}=U^{(k)}-\mathbf1(\bar u^{(k)})^\top$, and define the block covariances
\[
\Sigma_{ij}^{(k)}=
\mathbb E\!\left[(\widehat U_i^{(k)}-U_i^{(k)})
 (\widehat U_j^{(k)}-U_j^{(k)})^\top
 \mid\mathcal F_{\tau_k-}\right].
\]
Let $W\in\mathbb R^{m\times m}$ be fixed under this conditioning and mean
preserving, $W\mathbf1=\mathbf1$ and $\mathbf1^\top W=\mathbf1^\top$, and for
every node $v$ let $S(v)$ be a fixed set with $W_{vi}=0$ for $i\notin S(v)$
(the choice $S(v)=\{1,\ldots,m\}$ is always admissible).  Put
\[
\mathcal N^{(k)}:=\frac1m\sum_{i,j=1}^m(W^\top W)_{ij}
\mathop{\rm Tr}\nolimits\Sigma_{ij}^{(k)},
\qquad
\mathcal B^{(k)}:=\frac1m\|(W-J)D^{(k)}\|_{\mathrm F}^2+\mathcal N^{(k)} .
\]
The root $R_k$ has a conditional law $\pi$ and is independent of the latent row
innovations, and $\omega$ is a positive weight fixed under the conditioning.
The selected rows are produced in one of two chronologies: \emph{rows first},
where the latent rows are fixed before $R_k$ is revealed, or \emph{root first},
where after $R_k=v$ the rows indexed by $S(v)$ are sampled jointly from the
restriction of the latent law, with randomness independent of $R_k$.  The
estimator is
\begin{equation}\label{eq:uniform-row-estimator}
\widehat g^{(k)}=\omega(R_k)\,e_{R_k}^\top W\widehat U^{(k)},
\qquad
\Omega:=\max_v m\,\pi(v)\,\omega(v)^2 .
\end{equation}

\begin{lemma}[Selected-row oracle]\label{lem:uniform-row}
In the setting above the following statements hold.
\begin{enumerate}[label=(\roman*),leftmargin=2.2em]
\item \emph{Law.}  Given $R_k=v$, the selected row
$e_v^\top W\widehat U^{(k)}$ has the same conditional law in both
chronologies; root first evaluates only the $|S(v)|$ rows of $S(v)$.
\item \emph{Importance weight, $\pi\omega\equiv1/m$.}  Then
$\mathbb E[\widehat g^{(k)}\mid\mathcal F_{\tau_k-}]=\bar u^{(k)}$ and
\begin{equation}\label{eq:uniform-row-covariance}
\mathbb E\bigl[\|\widehat g^{(k)}-\bar u^{(k)}\|^2\mid\mathcal F_{\tau_k-}\bigr]
\leq\Omega\,\mathcal B^{(k)}+(\Omega-1)\|\bar u^{(k)}\|^2 ,
\end{equation}
with equality for a uniform root and $\omega\equiv1$, where $\Omega=1$.
\item \emph{No weight, $\omega\equiv1$.}  Then
$\mathbb E[\widehat g^{(k)}\mid\mathcal F_{\tau_k-}]=\bar u^{(k)}+e_k$ with
$e_k^\top=(\pi-\mathbf1/m)^\top(W-J)D^{(k)}$, so
$\|e_k\|\leq\|\pi-\mathbf1/m\|_2\|W-J\|_2\|D^{(k)}\|_{\mathrm F}$, and the
conditional centered variance is at most $\Omega\,\mathcal B^{(k)}$.
\end{enumerate}
Moreover $m^{-1}\|(W-J)D^{(k)}\|_{\mathrm F}^2\leq
\|W-J\|_2^2\,m^{-1}\|D^{(k)}\|_{\mathrm F}^2$, and for $W=I_m$ with
$S(v)=\{v\}$ (lazy depth zero) only the diagonal blocks $\Sigma_{ii}^{(k)}$
enter $\mathcal N^{(k)}$, so the bounds do not depend on how unevaluated
latent rows are coupled.
\end{lemma}

\begin{proof}
\smallskip\noindent\textit{Step 1: law.}\quad
Since $W_{vi}=0$ for $i\notin S(v)$, the selected row
$e_v^\top W\widehat U^{(k)}=\sum_{i\in S(v)}W_{vi}\widehat U_i^{(k)}$ is a fixed
linear function of the rows indexed by $S(v)$.  In both chronologies these rows
have, given $R_k=v$, the restricted latent law, because the root is independent
of the latent innovations and the root-first randomness is independent of the
root.  This proves~(i); all expectations below may therefore be computed under
the latent construction.

\smallskip\noindent\textit{Step 2: the average over all rows.}\quad
Write $E^{(k)}=\widehat U^{(k)}-U^{(k)}$.  Since $\mathbf1^\top D^{(k)}=0$ and
$WJ=JW=J$, $W\widehat U^{(k)}=\mathbf1(\bar u^{(k)})^\top+(W-J)D^{(k)}+WE^{(k)}$.
The three terms are orthogonal in expectation: the first two are
Frobenius-orthogonal because $\mathbf1^\top(W-J)=0$, and $E^{(k)}$ is
conditionally centered.  Hence
\begin{equation}\label{eq:all-rows-identity}
\frac1m\sum_v\mathbb E\|e_v^\top W\widehat U^{(k)}-\bar u^{(k)}\|^2=\mathcal B^{(k)},
\qquad
\frac1m\sum_v\mathbb E\|e_v^\top W\widehat U^{(k)}\|^2=\|\bar u^{(k)}\|^2+\mathcal B^{(k)},
\end{equation}
where $m^{-1}\mathbb E\|WE^{(k)}\|_{\mathrm F}^2=\mathcal N^{(k)}$ by expanding
the square.

\smallskip\noindent\textit{Step 3: means.}\quad
By Step~1,
$\mathbb E[\widehat g^{(k)}\mid\mathcal F_{\tau_k-}]^\top
=\sum_v\pi(v)\omega(v)\,e_v^\top WU^{(k)}$.  If $\pi\omega\equiv1/m$, this is
$m^{-1}\mathbf1^\top WU^{(k)}=(\bar u^{(k)})^\top$.  If $\omega\equiv1$, write
$WU^{(k)}=\mathbf1(\bar u^{(k)})^\top+(W-J)D^{(k)}$ and use
$(\pi-\mathbf1/m)^\top\mathbf1=0$; this gives the formula for $e_k$, and the
norm bound follows from
$\|D^{(k)}\|_2\leq\|D^{(k)}\|_{\mathrm F}$.

\smallskip\noindent\textit{Step 4: variances.}\quad
If $\pi\omega\equiv1/m$, then $\pi\omega^2=\omega/m$ and
$\Omega=\max_v\omega(v)$, so by~\eqref{eq:all-rows-identity}
\[
\mathbb E\|\widehat g^{(k)}\|^2
=\frac1m\sum_v\omega(v)\,\mathbb E\|e_v^\top W\widehat U^{(k)}\|^2
\leq\Omega\bigl(\|\bar u^{(k)}\|^2+\mathcal B^{(k)}\bigr),
\]
with equality when $\omega\equiv1$.  Subtracting the squared mean
$\|\bar u^{(k)}\|^2$ gives~\eqref{eq:uniform-row-covariance}.  If
$\omega\equiv1$, the centered variance is at most the second moment about
$\bar u^{(k)}$, which is
$\sum_v\pi(v)\,\mathbb E\|e_v^\top W\widehat U^{(k)}-\bar u^{(k)}\|^2
\leq\Omega\,\mathcal B^{(k)}$ by~\eqref{eq:all-rows-identity}.  The final
claims follow from $\|(W-J)D\|_{\mathrm F}\leq\|W-J\|_2\|D\|_{\mathrm F}$ and
from $(W^\top W)_{ij}=\delta_{ij}$ when $W=I_m$.
\end{proof}

We keep the noise term in block form because pairwise nonpositive
cross-covariances do not automatically reduce it for a signed polynomial
filter: the entries of $W^\top W$ need not be nonnegative.  The
simpler scalar term is used only under the direct domination
\begin{equation}\label{eq:row-noise-trace-domination}
\frac1m\sum_{i,j=1}^m(W^\top W)_{ij}
\mathop{\rm Tr}\nolimits\Sigma_{ij}^{(k)}
\leq\frac{\|W\|_{\mathrm F}^2}{m}\frac{\sigma^2}{b}.
\end{equation}
Conditional cross-node uncorrelatedness is sufficient.  More generally,
$\Sigma^{(k)}\preceq(\sigma^2/(bd))I_{md}$ is a checkable correlated-noise
example that implies~\eqref{eq:row-noise-trace-domination}.  For ordinary
rows $U_i^{(k)}=\nabla f_i(y^{(k)})$, so $\bar u^{(k)}=\nabla F(y^{(k)})$ and
$m^{-1}\|D^{(k)}\|_{\mathrm F}^2=\zeta^2(y^{(k)})$; then
\begin{equation}\label{eq:bracket-ordinary}
\mathcal B^{(k)}\leq\|W-J\|_2^2\,\zeta^2(y^{(k)})
 +\frac{\|W\|_{\mathrm F}^2}{m}\frac{\sigma^2}{b}.
\end{equation}

The network enters only through $\|W-J\|_2$ and $m^{-1}\|W\|_{\mathrm F}^2$.
On a fixed graph the Chebyshev recurrence of
Section~\ref{sec:network-oracle-design} controls both.

\begin{lemma}[Chebyshev disagreement bound]\label{lem:chebyshev-filter}
Let a fixed connected undirected graph have Laplacian eigenvalues
$0=\lambda_1<\lambda_2\leq\cdots\leq\lambda_m$ and
$\chi=\lambda_m/\lambda_2>1$.  The degree-$q$ normalized Chebyshev protocol
has a mean-preserving input-output operator $W_q$ satisfying
\begin{equation}\label{eq:chebyshev-profile}
\|W_q-J\|_2
\leq2\left(\frac{\sqrt\chi-1}{\sqrt\chi+1}\right)^q
\leq2\exp\!\left(-\frac{2q}{\sqrt\chi+1}\right),
\qquad
\frac1m\|W_q\|_{\mathrm F}^2
\leq\frac1m+\|W_q-J\|_2^2.
\end{equation}
Certified spectral endpoints may replace the exact eigenvalues.  If
$\chi=1$, one round gives $W_1=J$.
\end{lemma}

The proof is in Appendix~\ref{app:network-proofs}.  The
$\mathcal O(\sqrt\chi\log(1/\delta))$ depth concerns only pre-mixing; every
CQA-Net event also pays the broadcast depth $D$.  On changing graphs an
instantaneous spectral gap cannot be inserted into a noncommuting product, so
we condition on the actual input-output operator and require bounds on that
operator itself; Corollary~\ref{cor:block-contractive-changing} below gives one
explicit class.

\begin{theorem}[Network oracle inequality for CQA-Net]
\label{thm:static-network}
Suppose that $F$ is $L$-smooth and $\gamma$-quasar-convex, the local rows
satisfy~\eqref{eq:local-oracle}, and CQA-Net uses the root--row chronology of
Section~\ref{sec:communication-model}, including the lazy
kernel~\eqref{eq:lazy-root-row-kernel} when $q=0$, over a horizon $T>0$ with
$\theta\in(0,1]$.  At jump $k\geq1$ an event-indexed scheduler chooses a
mean-preserving $W_{k,q}$ from the revealed past before $E_{k-1}$ is exposed,
and the centered row noise satisfies~\eqref{eq:row-noise-trace-domination} for
the selected $W_{k,q}$.  Let $\bar w_q\geq\|W_{k,q}-J\|_2$ and
$\bar\nu_q\geq m^{-1}\|W_{k,q}\|_{\mathrm F}^2$ be deterministic uniform bounds
(or essential suprema over all revealed histories), and assume the
finite-horizon conditions of Theorem~\ref{thm:master-inexact}.  Then
\begin{align}
&\mathbb E\bigl[F(x_T)-F(x^\star)\bigr]
 +\frac{1-\theta^2}{2LT^2}
  \int_0^Tt^2\mathbb E\|\nabla F(x_t)\|^2\,dt\notag\\
&\quad\leq\frac{2LR_0^2}{\gamma^2\theta^2T^2}
 +\frac{V_\theta\bar w_q^2}{LT^2}
  \int_0^Tt^2\mathbb E\zeta^2(x_t)\,dt
 +\frac{V_\theta\bar\nu_q\sigma^2T}{3Lb}.
\label{eq:static-network-rate}
\end{align}
On a fixed graph with $W_{k,q}=W_q$ one may take $\bar w_q=\|W_q-J\|_2$ and
$\bar\nu_q=m^{-1}\|W_q\|_{\mathrm F}^2$.  The expected event count is $T$.  For
$q>0$, the direct implementation uses $bT$ expected samples per node and $mbT$
in aggregate; for $q=0$, root-first sampling uses $bT/m$ per node and $bT$ in
aggregate.  On a fixed graph, expected critical-path rounds and transmitted
vectors are at most $T(q+D)$ and $T(2q|E|+M_{\rm br})$.  On a changing graph
where the selected row can be disseminated within $B_{\rm br}$ subrounds and
$M_{\rm br}^{\rm tv}$ transmitted vectors and mixing uses at most
$M_{\rm mix}^{\rm tv}(q)$ vectors, they are at most $T(q+B_{\rm br})$ and
$T[M_{\rm mix}^{\rm tv}(q)+M_{\rm br}^{\rm tv}]$.
\end{theorem}

\begin{proof}
Fix jump $k\geq1$.  The operator $W_{k,q}$ is measurable with respect to the
history through jump $k-1$, and after the interarrival has produced the
predictable query, $W_{k,q}$ and the query are fixed under
$\mathcal F_{\tau_k-}$ while the rows and the root remain fresh.
Lemma~\ref{lem:uniform-row}(ii) with a uniform root, $\omega\equiv1$ and
$S(v)=\{1,\ldots,m\}$ (or $W=I_m$, $S(v)=\{v\}$ at $q=0$) gives conditional
mean $\nabla F(x_{\tau_k-})$, and~\eqref{eq:bracket-ordinary} gives the
centered variance
\[
\varsigma_t^2\leq\bar w_q^2\,\zeta^2(x_{t-})
 +\bar\nu_q\frac{\sigma^2}{b}.
\]
Insert this envelope into~\eqref{eq:master-unbiased} with $\eta=0$.  The
heterogeneity remains inside the weighted integral because it depends on the
trajectory; the mini-batch term is constant and integrates using
$T^{-2}\int_0^Tt^2dt=T/3$.  This step uses the common-state property:
Theorem~\ref{thm:master-inexact} is applied to $x_t\in\mathbb R^d$, not to a
matrix.  Finally, $\mathbb E n_T=T$, and the two branches of
Algorithm~\ref{alg:network-oracle} together with the certified per-event
mixing and dissemination costs give the ledgers.
\end{proof}

Theorem~\ref{thm:static-network} reduces the network problem to externally
certified products and does not prove a new consensus result.  It is also an
oracle inequality: the integral of $\zeta^2(x_t)$ depends on the trajectory.  There are
three useful closures.  On changing graphs each of them holds with
$\|W_q-J\|_2$ and $m^{-1}\|W_q\|_{\mathrm F}^2$ replaced by $\bar w_q$ and
$\bar\nu_q$.

\begin{corollary}[Closed network rates]\label{cor:closed-static}
Under Theorem~\ref{thm:static-network}, the following statements hold.

\emph{(a) Gronwall closure.}  If $\zeta^2(x)\leq h_0+h_1[F(x)-F(x^\star)]$
globally, then
\begin{equation}\label{eq:affine-closed-rate}
\mathbb E\bigl[F(x_T)-F(x^\star)\bigr]
\leq e^{V_\theta\|W_q-J\|_2^2h_1T/L}
\Bigl[\frac{2LR_0^2}{\gamma^2\theta^2T^2}
+\frac{V_\theta T}{3L}\Bigl(\|W_q-J\|_2^2h_0
 +\frac{\|W_q\|_{\mathrm F}^2}{m}\frac{\sigma^2}{b}\Bigr)\Bigr].
\end{equation}
It is at most $\varepsilon$ if
$T\geq2\sqrt LR_0/(\gamma\theta\sqrt\varepsilon)$, the exponent is at most
$\log(4/3)$, and the round bracket is at most $3L\varepsilon/(4V_\theta T)$.

\emph{(b) Global envelope.}  If $\zeta^2(x)\leq\bar\zeta^2$ globally,
then~\eqref{eq:affine-closed-rate} holds with $h_0=\bar\zeta^2$, $h_1=0$ and
no exponential factor.  It is at most $\varepsilon$ if
$T\geq2\sqrt LR_0/(\gamma\theta\sqrt\varepsilon)$ and each of the two variance
terms in the round bracket is at most $3L\varepsilon/(4V_\theta T)$; for the
Chebyshev protocol of Lemma~\ref{lem:chebyshev-filter} the disagreement
condition is implied by
$q\geq\frac{\sqrt\chi+1}{4}[\log(16V_\theta\bar\zeta^2T/(3L\varepsilon))]_+$,
where $[a]_+:=\max\{a,0\}$.

\emph{(c) Dissipative closure.}  Suppose quadratic growth holds,
$\zeta^2(x)\leq h_0+h_1[F(x)-F(x^\star)]$, $\theta\in(0,1)$, and
\begin{equation}\label{eq:absorption-condition}
\frac{2V_\theta h_1\|W_q-J\|_2^2}{\gamma^2\mu}
\leq\frac{1-\theta^2}{2}.
\end{equation}
Then~\eqref{eq:affine-closed-rate} holds without the exponential factor.
\end{corollary}

\begin{proof}
(a) Insert the affine bound into the variance envelope of
Theorem~\ref{thm:static-network}.  Since
$A_t[F(x_{t-})-F(x^\star)]\leq\Phi_{t-}$, the gap-dependent term contributes at
most $(V_\theta\|W_q-J\|_2^2h_1/L)\Phi_{t-}$ to the stopped drift; apply
Lemma~\ref{lem:gronwall}(i), remove the stopping as in
Theorem~\ref{thm:master-inexact}, and use $A_T^{-1}\int_0^TA_tdt=T/3$.  The
stated conditions bound the bracket by $3\varepsilon/4$, and the exponential
factor is at most $4/3$.  (b) Substitute $\zeta^2\leq\bar\zeta^2$ directly
into~\eqref{eq:static-network-rate}.  (c) Quasar-convexity, quadratic growth, and
Cauchy--Schwarz give $\|\nabla F(x)\|^2\geq\gamma^2\mu[F(x)-F(x^\star)]/2$, so the
gap-dependent variance is at most $\eta\|\nabla F(x_{t-})\|^2$ with
$\eta=2h_1\|W_q-J\|_2^2/(\gamma^2\mu)$, and~\eqref{eq:absorption-condition} is
$V_\theta\eta\leq(1-\theta^2)/2$; apply Theorem~\ref{thm:master-inexact} with
this $\eta$ and drop the gradient integral.
\end{proof}

Part~(a) is general but may be vacuous on long horizons; part~(c) is
horizon-free but requires quadratic growth, damping, and enough mixing.  Their
common premise holds with $h_0=2\zeta^2(x^\star)$ and $h_1=4\bar\delta^2/\mu$
under quadratic growth if
$\|[\nabla f_i(x)-\nabla F(x)]-[\nabla f_i(y)-\nabla F(y)]\|\leq\delta_i\|x-y\|$
and $\bar\delta^2=m^{-1}\sum_i\delta_i^2$ (Appendix~\ref{app:network-proofs}).
Mixing as accurately as possible is not generally optimal: if $T_q$ is the
smallest horizon satisfying a closed bound at depth $q$, the communication
choice is $q_{\rm comm}\in\arg\min_{q\geq0}(q+D)T_q$, so positive finite depth
trades computation against communication and does not always help.

\begin{remark}[Variance reduction for shared-landscape finite sums]\label{rem:vr}
The additive variance in~\eqref{eq:static-network-rate} grows with the horizon.
For finite sums $f_i=n_i^{-1}\sum_jf_{ij}$ under quadratic growth, restarting
CQA-Net every $T_{\rm st}=4\sqrt{L/\mu}/(\gamma\theta)$ time units from an exact
SVRG snapshot, computed by one all-reduce per stage, makes the inner variance
proportional to the current and snapshot gaps (Algorithm~\ref{alg:restarted-vr}).
In the shared-landscape case $f_{ij}=a_{ij}F+\langle c_{ij},x\rangle+d_{ij}$, where
node and component functions may have different minimizers and $F$ may have
negative curvature, put $a_i=n_i^{-1}\sum_ja_{ij}$,
$v_{\rm node}=m^{-1}\sum_i(a_i-1)^2$, and
$v_{\rm comp}=\max_in_i^{-1}\sum_j(a_{ij}-a_i)^2$.  If
\[
v_{\rm node}\|W_q-J\|_2^2+\frac{v_{\rm comp}}{b}\frac{\|W_q\|_{\mathrm F}^2}{m}
\leq\frac{\gamma\theta\log(4/3)}{16V_\theta}\sqrt{\frac\mu L},
\]
every stage contracts the expected gap by the factor
$r_{\rm st}=\frac43(\frac14+\frac13\log\frac43)\approx0.461$, so
$\mathcal O(\log(1/\varepsilon))$ stages suffice with an inner depth $q$ and a
batch $b$ that do not depend on $\varepsilon$ (Theorem~\ref{thm:vr-main} and
Corollary~\ref{cor:vr-shared}).  The exact snapshot is again a global
primitive, and generic distance-based similarity certificates lose a factor
$1/\mu$; Appendix~\ref{app:vr} gives the statements, the proofs, and a numerical
check on a nonconvex instance where the threshold holds.
\end{remark}

\smallskip\noindent\textit{Complexity to a target accuracy.}\quad
Set
\[
T_\varepsilon=\Theta\!\left(
  \frac{\sqrt L R_0}{\gamma\theta\sqrt\varepsilon}
\right).
\]
Let $q_\varepsilon,b_\varepsilon$ satisfy a chosen closure above.  A CQA-Net
event costs $q_\varepsilon+D$ critical-path rounds and $2q_\varepsilon|E|+M_{\rm br}$
transmitted vectors, and uses $mb_\varepsilon$ row samples for $q_\varepsilon>0$
and $b_\varepsilon$ for $q=0$, so the total round count is
$(q_\varepsilon+D)T_\varepsilon$.

At $\gamma=1$ the optimal convex local-state scale is
$\widetilde{\mathcal O}(\sqrt\chi\sqrt LR_0/\sqrt\varepsilon)$ rounds
\citep{scaman2017}.  Up to logarithmic and numerical factors, the ratio of the
common-state round bound to this scale is roughly
$(D+\sqrt\chi\log(1/\varepsilon))/\sqrt\chi$.  On a bounded-degree expander
family,
$\sqrt\chi=\mathcal O(1)$ whereas a common broadcast still has
$D=\Theta(\log m)$.  In the stochastic regime,
$q_\varepsilon$ generally grows logarithmically with $1/\varepsilon$, while
the sufficient batch size scales as $\varepsilon^{-3/2}$ because
$T_\varepsilon=\Theta(\varepsilon^{-1/2})$ and the stochastic term must be at
most $\varepsilon$.  The resulting aggregate stochastic sample scale remains
$\Theta(\varepsilon^{-2})$; acceleration acts on the deterministic term.
D-CQA has no dissemination term: its round count
$q_\varepsilon T_\varepsilon$ (Corollary~\ref{cor:dcqa-complexity}) matches the local-state scale up to a logarithm,
at the price of the assumptions of Section~\ref{sec:dcqa-results}.

Block contraction gives one explicit class of certificates for
Theorem~\ref{thm:static-network} on changing graphs; an unfinished block must
be controlled separately.

\begin{corollary}[Uniformly block-contractive changing graphs]
\label{cor:block-contractive-changing}
Suppose each complete block of $B_{\rm blk}$ communication rounds implements a
mean-preserving $Q_\ell$ with
$\|Q_\ell-J\|_2\leq\eta_{\rm blk}<1$.  For each remainder
$s\in\{0,\ldots,B_{\rm blk}-1\}$, suppose every admissible partial product
$R_s$ is mean preserving and
$\|R_s-J\|_2\leq\kappa_s$, with $\kappa_0=1$.  If
$q=rB_{\rm blk}+s$, then
\begin{equation}\label{eq:block-profile-corrected}
\|W_{k,q}-J\|_2\leq\kappa_s\eta_{\rm blk}^r,
\qquad
\frac1m\|W_{k,q}\|_{\mathrm F}^2
\leq\frac1m+\kappa_s^2\eta_{\rm blk}^{2r}.
\end{equation}
For $q=rB_{\rm blk}$ the remainder factor is absent.  Products of symmetric
lazy doubly stochastic maps are nonexpansive on the disagreement subspace, so
one may take $\kappa_s\leq1$.
\end{corollary}

\begin{proof}
Factorization controls the disagreement operator.  The Frobenius estimate then
matches the quantity in the row-noise term.

\smallskip\noindent\textit{Step 1: product factorization.}\quad
Mean preservation gives $AJ=JA=J$ for every factor.  Therefore
$(A-J)(B-J)=AB-J$, and induction yields
\[
W_{k,q}-J=(R_s-J)(Q_r-J)\cdots(Q_1-J).
\]
Submultiplicativity gives the first bound in
\eqref{eq:block-profile-corrected}.

\smallskip\noindent\textit{Step 2: row-noise factor.}\quad
The orthogonal decomposition
$\|W_{k,q}\|_{\mathrm F}^2=\|J\|_{\mathrm F}^2+
\|W_{k,q}-J\|_{\mathrm F}^2$, together with
$\|J\|_{\mathrm F}^2=1$ and
$m^{-1}\|A\|_{\mathrm F}^2\leq\|A\|_2^2$, gives the second.
\end{proof}

Corollary~\ref{cor:block-contractive-changing} is one explicit way to supply the product norms $\bar w_q,\bar\nu_q$ required by Theorem~\ref{thm:static-network}.  Exact averaging sets $W_{k,q}=J$ and removes row-mean disagreement, but
the stochastic-row upper bound still contains $\sigma^2/(mb)$.  The exact-gradient
corollary is recovered only when the local rows are deterministic as well.

\subsection{The walking-token oracle}\label{sec:token-results}

This subsection analyzes Algorithms~\ref{alg:token-oracle}
and~\ref{alg:token-cqa}.  At event $k$ put
\[
\mathcal P_k:=P_{k,1}P_{k,2}\cdots P_{k,h_{\rm w}},
\qquad
\pi_k:=\mathcal P_k^\top e_{o_k},
\qquad
\varrho_k:=\|\mathcal P_k-J\|_2 .
\]
Given $\mathcal F_{\tau_k-}$, the vector $\pi_k$ is the conditional law of the
walk endpoint $R_k$, which is independent of the latent row innovations.  On a
fixed graph $\mathcal P_k=P^{h_{\rm w}}$ with $P=P^\top$; we put
$\lambda_P:=\max\{|\lambda_2(P)|,|\lambda_m(P)|\}<1$ and
$\rho_h:=m\lambda_P^{h_{\rm w}}$.  The token oracle is the selected-row oracle
of Lemma~\ref{lem:uniform-row} with root law $\pi_k$, root-first chronology,
$S(v)=B_q(v)$, and weight $\omega_k=1/(m\pi_k)$ in fixed mode or
$\omega_k\equiv1$ in changing mode.  Only the walk law is new.

\begin{lemma}[Walk endpoint law]\label{lem:walk-endpoint}
If every $P_{k,s}$ is doubly stochastic, then
$\|\pi_k-\mathbf1/m\|_2\leq\varrho_k$ and $\max_im\pi_k(i)\leq1+m\varrho_k$.
On a fixed graph with $P=P^\top$, moreover $\varrho_k=\lambda_P^{h_{\rm w}}$,
$|m\pi_k(i)-1|\leq\rho_h$ for every $i$, and
$\omega_k=1/(m(P^{h_{\rm w}})_{o_kR_k})=1/(m\pi_k(R_k))$.  If $\rho_h\leq1/2$,
every node is reached with positive probability and
$\max_i\omega_k(i)\leq(1-\rho_h)^{-1}\leq1+2\rho_h$.
\end{lemma}

\begin{proof}
Since $\mathcal P_k^\top\mathbf1=\mathbf1$ and $J(e_{o_k}-\mathbf1/m)=0$,
$\pi_k-\mathbf1/m=(\mathcal P_k^\top-J)(e_{o_k}-\mathbf1/m)$, and
$\|e_{o_k}-\mathbf1/m\|_2^2=(m-1)/m\leq1$.  This gives the first two claims.
For symmetric $P$, $\|P^{h_{\rm w}}-J\|_2=\lambda_P^{h_{\rm w}}$, each
coordinate deviation is at most the Euclidean norm, and
$\pi_k(v)=(P^{h_{\rm w}})_{o_kv}$.
\end{proof}

\begin{corollary}[Fixed-graph token rate]\label{cor:token-fixed}
Assume the hypotheses of Theorem~\ref{thm:static-network} on $F$ and on the
rows, a fixed graph, $\theta\in(0,1)$, a ball-supported $W_q$, and
\begin{equation}\label{eq:token-hops}
h_{\rm w}\geq\frac{\ln\bigl(4mV_\theta/(1-\theta^2)\bigr)}{\ln(1/\lambda_P)},
\qquad\text{so that}\qquad
\rho_h\leq\frac{1-\theta^2}{4V_\theta}\ \Bigl(\leq\frac12\Bigr).
\end{equation}
Then the oracle of Algorithm~\ref{alg:token-oracle} in fixed mode is unbiased,
its centered variance is at most
$(1+2\rho_h)\bigl[\|W_q-J\|_2^2\zeta^2(y^{(k)})+\|W_q\|_{\mathrm F}^2\sigma^2/(mb)\bigr]
+2\rho_h\|\nabla F(y^{(k)})\|^2$, and Algorithm~\ref{alg:token-cqa} satisfies
\begin{align*}
&\mathbb E\bigl[F(x_T)-F(x^\star)\bigr]
 +\frac{(1-\theta^2)/2-2V_\theta\rho_h}{LT^2}
  \int_0^Tt^2\mathbb E\|\nabla F(x_t)\|^2\,dt\\
&\quad\leq\frac{2LR_0^2}{\gamma^2\theta^2T^2}
 +(1+2\rho_h)\frac{V_\theta\|W_q-J\|_2^2}{LT^2}
  \int_0^Tt^2\mathbb E\zeta^2(x_t)\,dt
 +(1+2\rho_h)\frac{V_\theta\|W_q\|_{\mathrm F}^2\sigma^2T}{3Lmb}.
\end{align*}
Corollary~\ref{cor:closed-static}(a) and~(b) hold with their variance terms
multiplied by $1+2\rho_h$.  Corollary~\ref{cor:closed-static}(c) holds if
condition~\eqref{eq:absorption-condition} is replaced by
\[
\frac{2V_\theta(1+2\rho_h)h_1\|W_q-J\|_2^2}{\gamma^2\mu}
\leq\frac{1-\theta^2}{2}-2V_\theta\rho_h.
\]
\end{corollary}

\begin{proof}
By Lemma~\ref{lem:walk-endpoint}, $\pi_k\omega_k\equiv1/m$ and
$\Omega=\max\omega_k\leq1+2\rho_h$.  Lemma~\ref{lem:uniform-row}(ii)
and~\eqref{eq:bracket-ordinary} give unbiasedness and the displayed variance,
which has the form~\eqref{eq:master-variance} with $\eta=2\rho_h$ and
$\bar\varsigma_t^2=(1+2\rho_h)[\|W_q-J\|_2^2\zeta^2(x_{t-})
+\|W_q\|_{\mathrm F}^2\sigma^2/(mb)]$.  The condition
$V_\theta\eta\leq(1-\theta^2)/2$ is exactly
$\rho_h\leq(1-\theta^2)/(4V_\theta)$.  Since $\omega_k\leq1+2\rho_h$, the
linear second-moment bound~\eqref{eq:linear-growth-oracle} is inherited from the
rows, and Lemma~\ref{lem:finite-moments} supplies the finite-horizon
conditions.  Theorem~\ref{thm:master-inexact}, in the
form~\eqref{eq:master-unbiased}, gives the rate.  For the closures, the affine
heterogeneity bound enters with the factor $1+2\rho_h$; in part~(c) the
multiplicative coefficient becomes
$\eta=2\rho_h+2(1+2\rho_h)h_1\|W_q-J\|_2^2/(\gamma^2\mu)$, and
$V_\theta\eta\leq(1-\theta^2)/2$ is the displayed condition.
\end{proof}

On a changing graph the endpoint law $\pi_k$ cannot be computed locally, so no
weight is used and the walk introduces a predictable bias.

\begin{corollary}[Changing-graph token rate]\label{cor:token-changing}
Let $F$ be $L$-smooth and $\gamma$-quasar-convex, let $\zeta(x)\leq\bar\zeta$
globally, let the rows satisfy~\eqref{eq:local-oracle}
and~\eqref{eq:row-noise-trace-domination}, and assume the uniform certificates
$\varrho_k\leq\bar\varrho$, $\|W_{k,q}-J\|_2\leq\bar w_q$ and
$m^{-1}\|W_{k,q}\|_{\mathrm F}^2\leq\bar\nu_q$ (for $q=0$,
$\bar w_0=\bar\nu_0=1$).  Put
$\bar e:=\sqrt m\,\bar\varrho\,\bar w_q\bar\zeta$ and
$\bar\sigma_{\rm w}^2:=(1+m\bar\varrho)[\bar w_q^2\bar\zeta^2+\bar\nu_q\sigma^2/b]$.  Then
the oracle of Algorithm~\ref{alg:token-oracle} in changing mode has a
predictable bias $\|e_k\|\leq\bar e$ and centered variance at most $\bar\sigma_{\rm w}^2$,
and the following hold under the corresponding finite-horizon conditions of
Theorem~\ref{thm:master-inexact}.
\begin{enumerate}[label=(\alph*),leftmargin=2.2em]
\item Without projection and for $\theta\in(0,1]$,
\[
\mathbb E[F(x_T)-F(x^\star)]
\leq\left(\sqrt{\frac{2LR_0^2}{\gamma^2\theta^2T^2}
 +\frac{V_\theta T}{3L}\bigl(\bar\sigma_{\rm w}^2+\bar e^2\bigr)}
 +\frac{(1+\theta)\theta T\bar e}{2\sqrt{2L}}\right)^{\!2}.
\]
\item If every $z$-jump is projected onto $\mathcal Q$ and $\theta\in(0,1)$,
\[
\mathbb E[F(x_T)-F(x^\star)]\leq\frac{2LR_0^2}{\gamma^2\theta^2T^2}
 +\frac{R_{\mathcal Q}}{\gamma}\bar e
 +\frac{C_{\theta,0}T}{3L}\bar e^2
 +\frac{V_\theta T}{3L}\bar\sigma_{\rm w}^2 .
\]
\end{enumerate}
If every block of $B_{\rm blk}$ hops contracts as in
Corollary~\ref{cor:block-contractive-changing}, then
$\bar\varrho\leq\kappa_s\eta_{\rm blk}^r$ for $h_{\rm w}=rB_{\rm blk}+s$; in
both bounds the bias terms are then below any prescribed $\varepsilon/4$ once
$h_{\rm w}$ is logarithmic in $1/\varepsilon$, and the variance term behaves as
in Corollary~\ref{cor:closed-static}(b) with the extra factor
$1+m\bar\varrho\leq2$ once $m\bar\varrho\leq1$.
\end{corollary}

\begin{proof}
Lemma~\ref{lem:uniform-row}(iii) with $\omega\equiv1$, $\Omega=\max_im\pi_k(i)$
and~\eqref{eq:bracket-ordinary}, together with Lemma~\ref{lem:walk-endpoint},
gives the bias $\|e_k\|\leq\varrho_k\bar w_q\sqrt m\,\zeta(y^{(k)})\leq\bar e$
and the variance bound $\bar\sigma_{\rm w}^2$.  The bias is predictable because it depends
only on $o_k$, the walk matrices, $W_{k,q}$, and the query.  Part~(a) is
Theorem~\ref{thm:master-inexact}(a) with $\eta=0$, constant $\bar e_t=\bar e$
and $\bar\varsigma_t=\bar\sigma_{\rm w}$, after dropping the gradient integral; part~(b) is
Theorem~\ref{thm:master-inexact}(b) with $\eta=0$, using
$\int_0^Tt\,dt=T^2/2$ and $\int_0^Tt^2dt=T^3/3$.  The block bound is
Corollary~\ref{cor:block-contractive-changing} applied to the walk matrices,
which are mean preserving, and the logarithmic choice of $h_{\rm w}$ follows
because all bias terms are increasing in $\bar\varrho$ and vanish as
$\bar\varrho\to0$.
\end{proof}

\begin{remark}[Length of the walk]\label{rem:token-length}
On the fixed ring with $m=16$ and Metropolis weights, $\lambda_P=0.94925$, and
for $\theta^2=1/2$ condition~\eqref{eq:token-hops} gives $h_{\rm w}=88$ hops
per event against $D=8$ for dissemination.  On bounded-degree graphs
$h_{\rm w}$ scales like $\chi\log m$, whereas
$D=\mathcal O(\sqrt\chi\log m)$ \citep{chung1989diameters}; nonreversible lifted
walks can shorten the walk, for instance from $\mathcal O(m^2)$ to
$\mathcal O(m)$ on the cycle \citep{chen1999lifting}, with
Lemma~\ref{lem:walk-endpoint} replaced by the mixing bound of the lifted chain.
Parallel work beyond $B_q(R_k)$, several tokens, loss of the token, and a token
version of finite-sum variance reduction are not covered here.
\end{remark}

\subsection{Local states: analysis of D-CQA}\label{sec:dcqa-results}

D-CQA needs two assumptions that the common-state analysis does not use.

\begin{assumption}[Local-state model]\label{ass:dcqa}
$F$ is $L$-smooth and $\gamma$-quasar-convex with respect to $x^\star$, and
$R_0=\|x^{(0)}-x^\star\|$.  In addition:
\begin{enumerate}[label=(\roman*),leftmargin=2.2em]
\item every $f_i$ is $L_{\rm loc}$-smooth; without loss of generality
$L_{\rm loc}\geq L$, and we set $\kappa_{\rm loc}:=L_{\rm loc}/L\geq1$;
\item $\zeta(y)\leq\bar\zeta$ for every $y\in\mathbb R^d$, with $\zeta$ as
in~\eqref{eq:heterogeneity};
\item at event $k$, conditionally on $\mathcal F_{\tau_k-}$, node $i$ returns
$g_i^{(k)}$ with mean $\nabla f_i(y_i^{(k)})$ and
$\mathbb E[\|g_i^{(k)}-\nabla f_i(y_i^{(k)})\|^2\mid\mathcal F_{\tau_k-}]\leq\sigma^2/b$;
noises of different nodes are conditionally uncorrelated, and the innovations
at event $k$ are conditionally independent of the interarrivals
$(E_j)_{j\geq k}$;
\item every $W_k$ is mean preserving, is chosen from the revealed history
before the interarrival that produces $\tau_k$, and satisfies
$\|W_k-J\|_2\leq\delta_{\rm mix}$.
\end{enumerate}
\end{assumption}

Items~(i) and~(ii) are new.  Item~(ii) holds, for example, for
shared-landscape models $f_i=F+\langle c_i,x\rangle$ with $\sum_ic_i=0$, where
$\zeta\equiv(m^{-1}\sum_i\|c_i\|^2)^{1/2}$; there the $f_i$ have different
minimizers and fail the quasar condition at $x^\star$.  Item~(ii) excludes the
partitioned-data regression of Section~\ref{sec:experiments}, where $\zeta$
grows with $\|y-x^\star\|$.  Put $A_t,c_t,V_\theta$ as
in~\eqref{eq:lyapunov-weights} and
\[
\Phi_t=A_t\bigl[F(\bar x_t)-F(x^\star)\bigr]+\frac12\|\bar z_t-x^\star\|^2 .
\]

\begin{lemma}[The network average is a biased CQA process]\label{lem:dcqa-average}
Under Assumption~\ref{ass:dcqa}, the averages $(\bar x_t,\bar z_t)$ follow
the flow~\eqref{eq:base-flow} between events.  At event $k$ they satisfy
\[
\bar x^{(k)}=\bar y^{(k)}-\bar g^{(k)}/L,\qquad
\bar z^{(k)}=\bar z^{(k-1)}-c_{\tau_k}\bar g^{(k)},\qquad
\bar g^{(k)}:=\frac1m\sum_ig_i^{(k)},
\]
with $\bar y^{(k)}=\bar x_{\tau_k-}$.  Hence they form the
process~\eqref{eq:base-process} driven by the oracle $\bar g$, whose
conditional moments are
\[
\mathbb E[\bar g^{(k)}\mid\mathcal F_{\tau_k-}]=\nabla F(\bar y^{(k)})+e_k,
\qquad
e_k=\frac1m\sum_i\bigl(\nabla f_i(y_i^{(k)})-\nabla f_i(\bar y^{(k)})\bigr),
\]
with predictable $e_k$, $\|e_k\|\leq L_{\rm loc}\|Y_\perp^{(k)}\|_{\mathrm F}/\sqrt m$,
and centered variance at most $\sigma^2/(mb)$.
\end{lemma}

\begin{proof}
Since $\mathbf1^\top W_k=\mathbf1^\top$, the averages of $W_k\widetilde X$ and
$W_k\widetilde Z$ are those of $\widetilde X$ and $\widetilde Z$.  The flow is
linear and identical at all nodes, so it commutes with averaging.  For the
bias, local smoothness and Cauchy--Schwarz give
$\|e_k\|\leq\frac{L_{\rm loc}}m\sum_i\|y_i-\bar y\|
\leq L_{\rm loc}\|Y_\perp\|_{\mathrm F}/\sqrt m$.  The queries are
$\mathcal F_{\tau_k-}$-measurable, so $e_k$ is predictable, and uncorrelated
node noises give the variance bound.
\end{proof}

\begin{lemma}[Disagreement bound]\label{lem:dcqa-disagreement}
Let Assumption~\ref{ass:dcqa} hold with a common initialization.  Put
$c_T:=\gamma\theta^2T/(2L)$ and $\bar\zeta_b:=\bar\zeta+\sigma/\sqrt b$, and
assume
\begin{equation}\label{eq:dcqa-psi}
\psi:=\delta_{\rm mix}\bigl(2+\kappa_{\rm loc}+\kappa_{\rm loc}\gamma\theta^2T/2\bigr)\leq\frac12 .
\end{equation}
Let $\mathcal E$ be the sigma-field generated by all interarrivals and
$s_k^2:=\mathbb E\bigl[\|X_\perp^{(k)}\|_{\mathrm F}^2+\|Z_\perp^{(k)}\|_{\mathrm F}^2\mid\mathcal E\bigr]/m$.
Then $s_k\leq\bar s:=2\delta_{\rm mix}(1/L+c_T)\bar\zeta_b$ for every event with
$\tau_k\leq T$.  Moreover $\mathbb E\|X_\perp(t)\|_{\mathrm F}^2/m\leq\bar s^2$
for every $t\in[0,T]$, and the bias
$e_t:=m^{-1}\sum_i[\nabla f_i(x_i(t-))-\nabla f_i(\bar x_{t-})]$ satisfies
\[
\mathbb E\|e_t\|^2\leq\bar e^2,\qquad
\bar e:=L_{\rm loc}\bar s
=2\delta_{\rm mix}\kappa_{\rm loc}\bigl(1+\gamma\theta^2T/2\bigr)\bar\zeta_b .
\]
\end{lemma}

\begin{proof}
\smallskip\noindent\textit{Step 0: conditioning on the clock.}\quad
Condition on $\mathcal E$; event times and the number of events in $[0,T]$
become deterministic.  By Assumption~\ref{ass:dcqa}(iii) and the tower
property, the conditional second moment of the noise at event $k$ given
$\mathcal E$ is still at most $\sigma^2/b$ per node, because $\mathcal E$ is
generated by $\mathcal F_{\tau_k-}$ and the future interarrivals.  Write
$\mathrm{rms}(A):=(\mathbb E[\|A\|_{\mathrm F}^2\mid\mathcal E]/m)^{1/2}$; it
satisfies Minkowski's inequality.

\smallskip\noindent\textit{Step 1: flow.}\quad
Between events $Z_\perp$ is constant and
$X_\perp(t)=\alpha X_\perp+(1-\alpha)Z_\perp$ with $\alpha\in[0,1]$.  By
convexity of $\|\cdot\|_{\mathrm F}^2$,
$\|X_\perp(t)\|_{\mathrm F}^2\leq\|X_\perp^{(k)}\|_{\mathrm F}^2+\|Z_\perp^{(k)}\|_{\mathrm F}^2$.
In particular the next query satisfies $\mathrm{rms}(Y_\perp)\leq s_k$.

\smallskip\noindent\textit{Step 2: one event.}\quad
Since $(I-J)W=W-J=(W-J)(I-J)$,
\[
X_\perp^{(k+1)}=(W-J)\bigl(Y_\perp-G_\perp/L\bigr),\qquad
Z_\perp^{(k+1)}=(W-J)\bigl(Z_\perp^{(k)}-c\,G_\perp\bigr),
\]
with $W=W_{k+1}$, $G_\perp=(I-J)G$, and $c=c_{\tau_{k+1}}\leq c_T$.

\smallskip\noindent\textit{Step 3: gradient disagreement.}\quad
Write
$G=[\nabla f_i(y_i)-\nabla f_i(\bar y)]_i+[\nabla f_i(\bar y)]_i+\Xi$, where
$\Xi$ is the noise.  The projection $I-J$ does not increase the Frobenius norm,
and $(I-J)[\nabla f_i(\bar y)]_i$ has Frobenius norm
$\sqrt m\,\zeta(\bar y)\leq\sqrt m\,\bar\zeta$.  Hence
$\mathrm{rms}(G_\perp)\leq L_{\rm loc}s_k+\bar\zeta+\sigma/\sqrt b$.

\smallskip\noindent\textit{Step 4: recursion.}\quad
Minkowski's inequality and $\|(W-J)A\|_{\mathrm F}\leq\delta_{\rm mix}\|A\|_{\mathrm F}$
give
\[
s_{k+1}\leq\mathrm{rms}(X_\perp^{(k+1)})+\mathrm{rms}(Z_\perp^{(k+1)})
\leq\delta_{\rm mix}\bigl[(2+\kappa_{\rm loc}+c_TL_{\rm loc})s_k+(1/L+c_T)\bar\zeta_b\bigr]
=\psi s_k+\frac{\bar s}2,
\]
where $c_TL_{\rm loc}=\kappa_{\rm loc}\gamma\theta^2T/2$.  Since $s_0=0$ and
$\psi\leq1/2$, induction gives $s_k\leq\bar s$.

\smallskip\noindent\textit{Step 5: fixed times.}\quad
For $t\leq T$ the index $n_t$ is $\mathcal E$-measurable, and Step~1 gives
$\mathbb E[\|X_\perp(t)\|_{\mathrm F}^2\mid\mathcal E]/m\leq s_{n_t}^2\leq\bar s^2$.
Taking expectations and applying the bias bound of
Lemma~\ref{lem:dcqa-average} at $x_i(t-)$ completes the proof.
\end{proof}

The two lemmas place D-CQA inside Theorem~\ref{thm:master-inexact}(a): the
network average is a CQA process whose oracle has a mean-square bias bounded
by $\bar e$, and the square-root Gronwall argument absorbs the linear bias
terms without a projection.

\begin{theorem}[Fully decentralized accelerated bound]\label{thm:dcqa}
Let Assumption~\ref{ass:dcqa} hold with a common initialization,
$\theta\in(0,1]$, and~\eqref{eq:dcqa-psi}.  Assume the finite-horizon
integrability conditions of Theorem~\ref{thm:master-inexact} for the stacked
process; they hold by Lemma~\ref{lem:finite-moments}, because the local
gradients grow at most linearly.  Then Algorithm~\ref{alg:dcqa} satisfies
\begin{equation}\label{eq:dcqa-rate}
\mathbb E[F(\bar x_T)-F(x^\star)]
\leq\left(\sqrt{\frac{2LR_0^2}{\gamma^2\theta^2T^2}
 +\frac{V_\theta T}{3L}\Bigl(\bar e^2+\frac{\sigma^2}{mb}\Bigr)}
 +\frac{(1+\theta)\,\theta\,T\,\bar e}{2\sqrt{2L}}\right)^{\!2},
\end{equation}
and the node outputs satisfy
$\frac1m\sum_i\mathbb E[F(x_i(T))-F(x^\star)]
\leq\text{(right side of~\eqref{eq:dcqa-rate})}+\bar e^2/(2\kappa_{\rm loc}^2L)$.
\end{theorem}

\begin{proof}
By Lemma~\ref{lem:dcqa-average}, $(\bar x_t,\bar z_t)$ is the
process~\eqref{eq:base-process} driven by an oracle with predictable bias
$e_t$ and centered variance at most $\sigma^2/(mb)$, so
Theorem~\ref{thm:master-inexact}(a) applies with $\eta=0$,
$\bar\varsigma_t^2=\sigma^2/(mb)$, and, by Lemma~\ref{lem:dcqa-disagreement},
the constant $\bar e_t=\bar e$.  Then
$\Gamma_T=2LR_0^2/(\gamma^2\theta^2T^2)+V_\theta T(\bar e^2+\sigma^2/(mb))/(3L)$
and $B_T=(1+\theta)\theta T\bar e/(2\sqrt{2L})$; dropping the nonnegative
gradient integral gives~\eqref{eq:dcqa-rate}.  For the node outputs,
smoothness of $F$ and $\sum_i(x_i-\bar x)=0$ give
$\frac1m\sum_iF(x_i)\leq F(\bar x)+\frac L{2m}\|X_\perp\|_{\mathrm F}^2$, and
Lemma~\ref{lem:dcqa-disagreement} bounds the expectation of the last term by
$L\bar s^2/2=\bar e^2/(2\kappa_{\rm loc}^2L)$.
\end{proof}

\begin{corollary}[Complexity of D-CQA]\label{cor:dcqa-complexity}
Let $\varepsilon>0$ and suppose that
\begin{enumerate}[label=(\roman*),leftmargin=2.2em]
\item $T\geq5\sqrt L R_0/(\gamma\theta\sqrt\varepsilon)$;
\item $b\geq4V_\theta T\sigma^2/(Lm\varepsilon)$ (void if $\sigma=0$);
\item the mixing accuracy satisfies $\delta_{\rm mix}\leq\delta_\varepsilon$, where
\[
\delta_\varepsilon:=\min\left\{\frac{1}{2(2+\kappa_{\rm loc}+\kappa_{\rm loc}\gamma\theta^2T/2)},\
\frac{\min\bigl\{\sqrt{L\varepsilon/(4V_\theta T)},\ \sqrt{2L\varepsilon}/((1+\theta)\theta T)\bigr\}}
 {2\kappa_{\rm loc}(1+\gamma\theta^2T/2)\,\bar\zeta_b}\right\}.
\]
\end{enumerate}
Then $\mathbb E[F(\bar x_T)-F(x^\star)]\leq\varepsilon$, and the average of the
node outputs is at most $(1+\frac1{4T})\varepsilon$.  On a fixed graph,
Chebyshev filtering with $q\geq\frac{\sqrt\chi+1}{2}\ln(2/\delta_\varepsilon)$
suffices.  For small $\varepsilon$ the binding term of $\delta_\varepsilon$ is
of order
$\gamma\varepsilon^{3/2}/(\theta\kappa_{\rm loc}\sqrt LR_0^2\bar\zeta_b)$, so
\[
q_\varepsilon=\mathcal O\Bigl(\sqrt\chi\,
 \log\frac{\kappa_{\rm loc}\sqrt LR_0^2\bar\zeta_b}{\gamma\varepsilon}\Bigr),
\qquad
\text{rounds}=Tq_\varepsilon=\widetilde{\mathcal O}\Bigl(
 \frac{\sqrt\chi\sqrt LR_0}{\gamma\theta\sqrt\varepsilon}\Bigr).
\]
Per event, each node evaluates one row and the network transmits $4|E|q$
vectors in $\mathbb R^d$.  On changing graphs any certificate of the
form~\eqref{eq:block-profile-corrected} applies: with $\kappa_s\leq1$ it
suffices to run $r\geq\ln(1/\delta_\varepsilon)/\ln(1/\eta_{\rm blk})$ blocks of
$B_{\rm blk}$ rounds per event.  No dissemination bound is needed.
\end{corollary}

\begin{proof}
Let $a^2$ be the expression under the square root in~\eqref{eq:dcqa-rate} and
$\beta_e$ the second summand.  Condition~(i) gives
$2LR_0^2/(\gamma^2\theta^2T^2)\leq2\varepsilon/25\leq\varepsilon/12$, and~(ii)
gives $V_\theta T\sigma^2/(3Lmb)\leq\varepsilon/12$.  Condition~(iii) gives
$\psi\leq1/2$, $V_\theta T\bar e^2/(3L)\leq\varepsilon/12$, and
$\beta_e^2\leq\varepsilon/4$.  Hence
$(a+\beta_e)^2\leq2a^2+2\beta_e^2\leq\varepsilon$.  The node term is
$\bar e^2/(2\kappa_{\rm loc}^2L)\leq\varepsilon/(8\kappa_{\rm loc}^2V_\theta T)\leq\varepsilon/(4T)$.
The Chebyshev depth follows from Lemma~\ref{lem:chebyshev-filter}.  For the
order of $\delta_\varepsilon$, substitute $T$ from~(i): after division by
$1+\gamma\theta^2T/2$, the second entry of the inner minimum scales like
$\varepsilon^{3/2}$, the first like $\varepsilon^{5/4}$, and the first entry of
the outer minimum like $\varepsilon^{1/2}$, so the second entry of the inner
minimum is binding.
\end{proof}

\begin{remark}[Why the mixing accuracy must shrink like $T^{-2}$]\label{rem:dcqa-T2}
The auxiliary step $c_t=\gamma\theta^2t/(2L)$ grows linearly in time.  Each
event therefore injects disagreement $c_t\zeta$ into $Z$, and the bias grows
like $\delta_{\rm mix}c_T\bar\zeta$.  In~\eqref{eq:dcqa-rate} this bias is
multiplied once more by $T$, through $\int_0^Tc_t\,dt$.  In CQA-Net the same
heterogeneity enters as variance,
$V_\theta\|W_q-J\|_2^2\bar\zeta^2T/(3L)$, and only
$\|W_q-J\|_2\lesssim\varepsilon^{3/4}$ is needed.  Both requirements give
$\mathcal O(\sqrt\chi\log(1/\varepsilon))$ rounds per event and differ only in
the constant inside the logarithm; CQA-Net additionally pays the
dissemination depth $D$.
\end{remark}

\subsection{Comparison of the three implementations}\label{sec:implementations-comparison}

Table~\ref{tab:decentralization} compares the cost of one event.  CQA-Net keeps the strongest
guarantee, an unbiased oracle at every mixing depth with no additional
assumptions, and pays the global primitives (G1)--(G3).  The token removes all
global primitives with no new assumption on a fixed graph, but its walk is
sequential and long.  D-CQA meets the usual synchronous notion of
decentralization and removes the dissemination depth from the round count, at
the price of local smoothness, bounded heterogeneity, and a mixing accuracy
that must improve with the horizon ($\delta_{\rm mix}\lesssim\varepsilon^{3/2}$,
against $\|W_q-J\|_2\lesssim\varepsilon^{3/4}$ for CQA-Net).  The main results
are Theorem~\ref{thm:static-network} for CQA-Net,
Corollaries~\ref{cor:token-fixed} and~\ref{cor:token-changing} for the token,
and Theorem~\ref{thm:dcqa} with Corollary~\ref{cor:dcqa-complexity} for D-CQA.

\section{Numerical experiments}\label{sec:experiments}

The experiments examine individual parts of the analysis under constants that
can be checked from the supplied data.  We record both deterministic horizons
and realized oracle work, and report computation separately from
communication.  Static and two-window-connected graphs are tested, and a
separate experiment isolates persistent stochastic noise.
Section~\ref{sec:experiment-broadcast-free} runs the walking token and D-CQA of
Sections~\ref{sec:token-design}--\ref{sec:dcqa-design}; all other experiments
use CQA-Net.  Appendix~\ref{app:reproducibility}
describes the code, runwise data, full-precision constants, and validation
checks.

\subsection{Protocol and cost accounting}\label{sec:experiment-protocol}

We simulate the closed-form flow exactly.  For every event index, the graph and
mixer are fixed from the revealed past before the current exponential
interarrival is generated; local rows and the fresh uniform root are then
sampled in the order required by Lemma~\ref{lem:clock-law-equivalence}.  Unless
stated otherwise, CQA-Net uses $\theta=1$.  Curves compared with the theory are
sample means over $24$ independent runs with pointwise normal-approximation
bands $\text{mean}\pm1.96\,\mathrm{SE}$; medians and interquartile intervals
are used for typical-run cost plots.  Values below $10^{-16}$ are clipped only
for logarithmic display and values below roughly $10^{-14}$ are interpreted as
numerical zero.

On a ring of diameter $D$, a finite-mixing event is charged $q+D$
critical-path rounds and exact reduction plus dissemination is charged $2D$.
At $q=0$, the root is selected before local computation, so only one node
forms a row.  For $q>0$, every node forms a row before mixing.  Aggregate
row/full-gradient or component work, critical-path rounds, and
transmitted-vector volume are kept separate, as labeled.  In the
regression ledger, one exact local full-gradient row counts as one evaluation.  The
baselines are exact-average gradient descent with step $1/L$,
classical convex Nesterov acceleration, and synchronous gradient tracking with
a lazy Metropolis matrix.  The tracking step was selected from the declared
grid $\{1,1/2,1/4,1/8,1/16,1/32\}/L$; $1/L$ minimized the final gap.  The
plotted tracker average is an offline diagnostic, and nodewise and consensus
errors are reported below.

The principal instances and their certified geometry are described in
Section~\ref{sec:quasar-overview}.  All manuscript constants are rounded; the
machine-readable metadata retain full precision.

\subsection{Certified nonconvex regression on fixed and changing graphs}
\label{sec:experiment-certified}

Each of $m=16$ nodes stores $120$ samples in dimension $18$.  Responses are
realizable under the link
$\varphi_{0.18}(u)=0.18u+0.82\tanh u$.  Proposition~\ref{prop:model-certificate}
certifies $\gamma=0.36$, $L=2.3631$, and
$\mu=1.457\cdot10^{-3}$; the initial gap is $8.543\cdot10^{-2}$.  On the ring,
the measured pairs
$(\|W_q-J\|_2,m^{-1}\|W_q\|_{\mathrm F}^2)$ are
$(1,1)$, $(0.55887,0.22483)$, $(0.12524,0.06913)$, and
$(0.011736,0.062572)$ for $q=0,3,7,13$.  Thus unbiasedness is present already
at $q=0$, while finite mixing reduces both row disagreement and, toward its
$1/m$ floor, independent row noise.

\begin{figure}[!htbp]
\centering
\includegraphics[width=.6\textwidth]{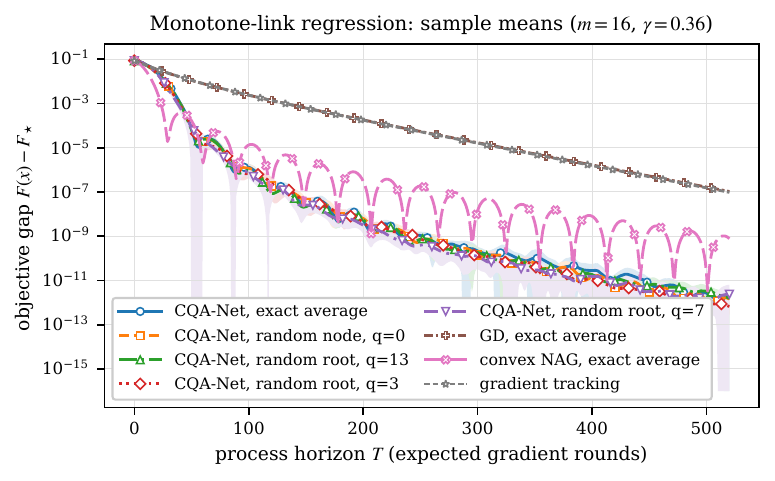}
\caption{Globally certified leaky-$\tanh$ regression on the static $16$-node
ring: sample-mean objective gap versus the process horizon.}
\label{fig:network-regression}
\end{figure}

Table~\ref{tab:static-ledger} merges accuracy and realized-work summaries, and
Figure~\ref{fig:network-regression} shows the sample-mean progress.  At gap $10^{-8}$, every CQA-Net run
reaches the target.

\begin{table}[!htbp]
\caption{Static-ring performance.  The first two columns use pointwise median
horizons; the remaining columns are medians of runwise first-passage work to
$10^{-8}$.}
\label{tab:static-ledger}
\centering
\scriptsize
\setlength{\tabcolsep}{3pt}
\begin{tabular}{@{}P{4.5cm}rrrrr@{}}
\toprule
Method & $T(10^{-6})$ & $T(10^{-8})$ & \shortstack{Actual events /\ iterations} & Rounds & \shortstack{Aggregate\ gradients} \\
\midrule
GD, exact average & 422 & not reached & not reached & not reached & not reached \\
Gradient tracking & 419 & not reached & not reached & not reached & not reached \\
CQA-Net, exact average & 85 & 166 & 119 & 1904 & 1904 \\
CQA-Net, $q=0$ & 83 & 170 & 147.5 & 1180 & 147.5 \\
CQA-Net, $q=3$ & 87 & 151 & 117 & 1287 & 1872 \\
CQA-Net, $q=7$ & 81 & 143 & 114 & 1710 & 1824 \\
CQA-Net, $q=13$ & 87 & 143 & 118 & 2478 & 1888 \\
Convex NAG, exact average & 88 & 177 & 177 & 2832 & 2832 \\
\bottomrule
\end{tabular}
\end{table}

At horizon $520$, realized Poisson counts range from $485$ to $570$, with
median $526$, so deterministic horizon and actual work are close only in
aggregate.  Along the exact representative trajectory the empirical quasar
ratio remains above $1.9917$, much larger than the global certificate $0.36$;
a separate near-tight one-dimensional diagnostic is reported below.

For a changing network, every communication subround is a disconnected perfect
matching.  Even and odd matchings alternate, so the union of two consecutive
subrounds is the ring.  The measured pairs are
\[
\begin{aligned}
q=1:&\quad (\|W_{k,q}-J\|_2,m^{-1}\|W_{k,q}\|_{\mathrm F}^2)=(1,0.5),\\
q=4:&\quad (\|W_{k,q}-J\|_2,m^{-1}\|W_{k,q}\|_{\mathrm F}^2)=(0.78858,0.15625),\\
q=10:&\quad (\|W_{k,q}-J\|_2,m^{-1}\|W_{k,q}\|_{\mathrm F}^2)=(0.49039,0.092804).
\end{aligned}
\]
Root dissemination is charged separately as $8$ matching subrounds, the exact
worst-case broadcast depth of the alternating $16$-node ring.  A full exact
reduction followed by broadcast costs $16$ matching subrounds.

The pointwise-median horizons to gap $10^{-10}$ are $287$, $294$, and $244$ for
$q=1,4,10$, compared with $267$ under exact averaging.  These values show that the optimization dynamics use the actual product norms supplied by the
external schedule, as stated in Theorem~\ref{thm:static-network}.

\subsection{Variance mechanisms: mixing, damping, and noise}
\label{sec:experiment-mechanisms}

\paragraph{Finite pre-mixing and the scope of the closed certificates.}
On the ring, consider $f_i(x)=a_ix^2/2$ with mean curvature one and local
curvature heterogeneity aligned with the fixed Fourier mode
$v_i=\sqrt{2}\cos(\pi i/2+\pi/2)$, $i=0,\ldots,15$, so that
$a_i=1+0.8v_i$.  This explicit phase avoids any dependence on a numerical
eigenbasis.  The network objective is $F(x)=x^2/2$ and
$\zeta^2(x)=1.28[F(x)-F(x^\star)]$.  With $q=3$, the selected mode is attenuated
by $0.06620$ even though the worst-case norm is
$\|W_3-J\|_2=0.55887$.  Empirically, the $q=0$ sample mean is dominated by rare
unstable paths and does not reach $10^{-6}$ by horizon $120$, whereas $q=3$
reaches it at horizon $26$, using $416$ aggregate gradients and $286$
communication rounds.  Exact averaging reaches it sooner, at horizon
$16$ with $256$ gradients and $256$ rounds, so the test contrasts mixing with the unstable unmixed sample mean and
makes no claim that it beats exact averaging.

The global worst-case norm certificate does not cover this example.  The
Gronwall small-exponent condition permits only $T\leq0.72$ at
$q=3$, while the factor at $T=26$ is about $3.27\cdot10^4$.  For
$\theta^2=1/2$, the dissipative condition
\eqref{eq:absorption-condition} has left-hand side $0.5997>0.25$ at $q=3$.
The favorable $q=3$ curves come from mode alignment; the global
$\|W_q-J\|_2$ certificate does not explain them.

At $q=7$, the global dissipative certificate does apply:
$2V_\theta h_1\|W_7-J\|_2^2=0.03012<0.25$ for $\theta^2=1/2$.  The resulting
ablation gives sample-mean hitting horizons $18$, $17$, and $24$ for
$\theta=1$, $\theta^2=1/2$, and $\theta=1/2$, respectively: moderate damping
absorbs the certified multiplicative noise, while stronger damping slows the
method down.

\paragraph{Persistent mini-batch noise.}
Balanced affine component perturbations create nonzero mini-batch variance
without changing any full local gradient or the global objective.  At horizon
$400$, the sample-mean full-gradient gap is $1.62\cdot10^{-9}$, whereas the
means for $b=16,8,4$ are $1.24\cdot10^{-3}$,
$2.06\cdot10^{-3}$, and $6.04\cdot10^{-3}$.  This is the noise accumulation
isolated by the $t^2$ variance integral.

\FloatBarrier

\subsection{Broadcast-free implementations}\label{sec:experiment-broadcast-free}

We run the walking token (Algorithm~\ref{alg:token-cqa}) and D-CQA
(Algorithm~\ref{alg:dcqa}) next to CQA-Net on the certified regression of
Section~\ref{sec:experiment-certified}, with $24$ independent runs to horizon
$520$.  All methods use $\theta^2=1/2$, because the token in fixed mode
requires $\theta<1$.  Costs per event follow
Table~\ref{tab:decentralization}: CQA-Net is charged $q+D$ rounds, the token
$h_{\rm w}+2q$ sequential rounds, and D-CQA $q$ rounds, in each of which every
node sends the pair $(x_i,z_i)$ to each neighbor.

On the fixed ring the weighted token walks with lazy Metropolis weights.  We
use $h_{\rm w}=48$, the smallest length for which the computed endpoint law
satisfies $\max_i|m\pi_k(i)-1|\leq1/6$; Lemma~\ref{lem:walk-endpoint} and
Corollary~\ref{cor:token-fixed} hold with this computed deviation in place of
$\rho_h$, whereas the spectral condition~\eqref{eq:token-hops} would require
$h_{\rm w}=88$.  We also run a lifted walk that carries a direction and
reverses it with probability $4/m$ per hop; it reaches the same deviation with
$h_{\rm w}=11$, and its exact endpoint law gives the weight, so
Lemma~\ref{lem:uniform-row}(ii) still makes the oracle unbiased.  D-CQA mixes
the states with the Chebyshev filter $W_q$, $q\in\{1,3,7\}$, and we report the
gap of the network average $\bar x$.  On the alternating matchings of
Section~\ref{sec:experiment-certified}, the token is unweighted, and D-CQA mixes
with the product of $q\in\{2,8,32\}$ matching subrounds.

\begin{figure}[!htbp]
\centering
\subfigure[Fixed ring: sample mean versus horizon.]{%
  \includegraphics[width=.48\textwidth]{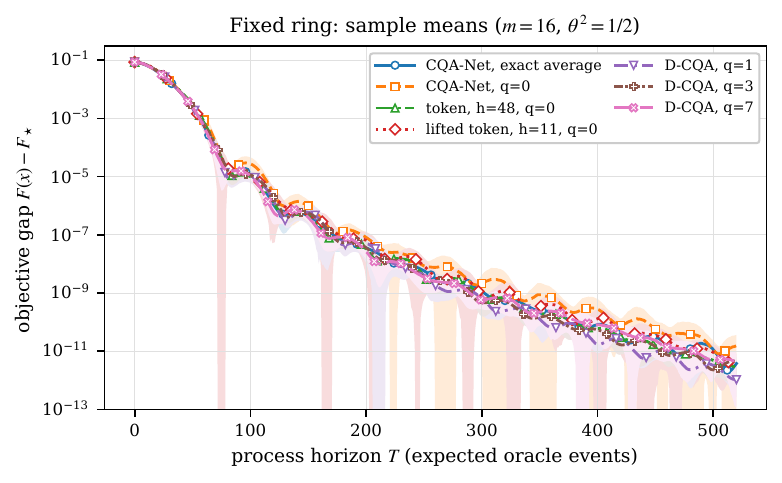}}\hfill
\subfigure[Fixed ring: median versus rounds.]{%
  \includegraphics[width=.48\textwidth]{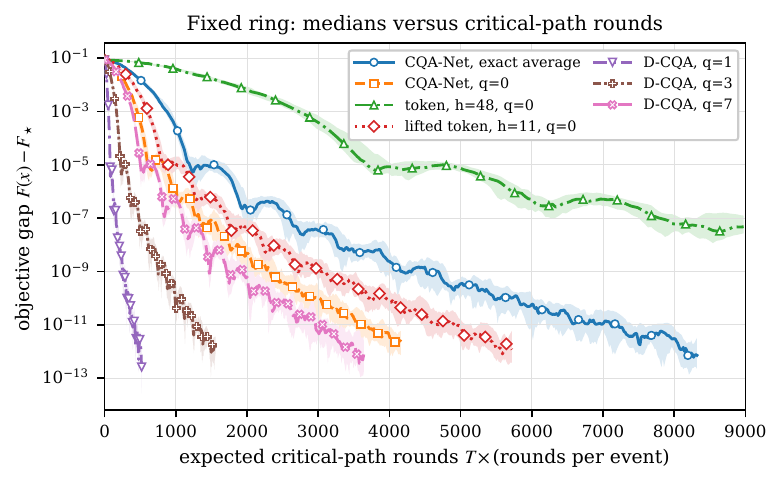}}

\medskip
\subfigure[Alternating matchings: median versus subrounds.]{%
  \includegraphics[width=.48\textwidth]{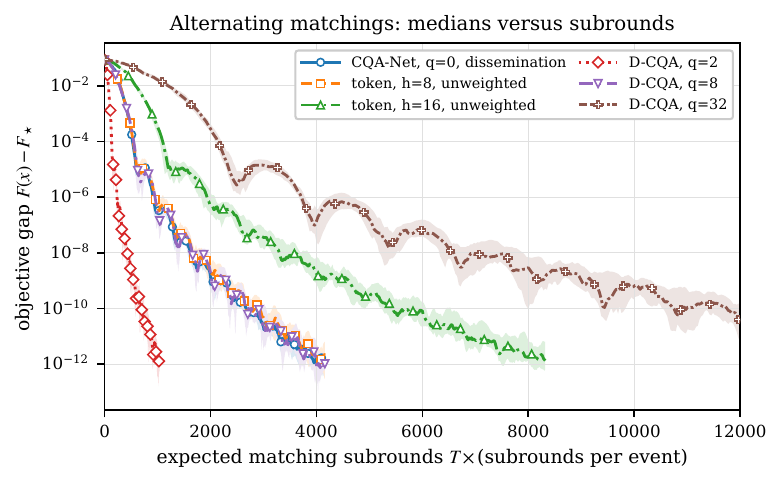}}
\caption{Broadcast-free implementations on the certified regression instance
($\theta^2=1/2$, $24$ runs).  (a)--(b)~CQA-Net, the weighted token (Metropolis
walk with $h_{\rm w}=48$, lifted walk with $h_{\rm w}=11$), and D-CQA with
Chebyshev mixing (network average $\bar x$).  (c)~Alternating matchings:
CQA-Net with dissemination, the unweighted token, and D-CQA with $q$ matching
subrounds per event.}
\label{fig:impl-main}
\end{figure}

Figure~\ref{fig:impl-main}(a) shows that per event all methods progress alike:
every run reaches $10^{-8}$, at median horizons between $202$ and $219$, and the
D-CQA node outputs reach it at the same horizons as $\bar x$.  The methods
differ in what one event costs (Figure~\ref{fig:impl-main}(b)).  To reach
$10^{-8}$, CQA-Net with $q=0$ needs $1{,}464$ rounds.  The Metropolis token
needs $8{,}520$, which is the latency price of a sequential walk, while the
lifted token needs only $1{,}991$.  D-CQA with one round per event needs $154$
rounds, $9.5$ times fewer than CQA-Net; it pays in computation, $2{,}464$ row
evaluations against $183$, and in data, $1.8\cdot10^5$ transmitted scalars
against $4.9\cdot10^4$.  On the alternating matchings
(Figure~\ref{fig:impl-main}(c)), D-CQA with $q=2$, one connected window per
event, needs $308$ subrounds, against $1{,}400$ for CQA-Net with dissemination
and $1{,}552$ for the token with $h_{\rm w}=8$; the token transmits the least
data, $2.9\cdot10^4$ scalars against $4.7\cdot10^4$ and $1.8\cdot10^5$.

These runs lie partly outside the hypotheses of the corresponding results.  The
partitioned regression violates item~(ii) of Assumption~\ref{ass:dcqa}, and the
unweighted token with $h_{\rm w}=8$ has $m\varrho_k=9.2$, far from the regime
$m\bar\varrho\leq1$ of Corollary~\ref{cor:token-changing}.  Both methods
nevertheless converge because the instance is realizable: every local gradient
vanishes at $x^\star$, so the walk bias and the disagreement injected by D-CQA
vanish there as well.  The counts are rounds and messages; they do not measure
the latency of sequential hops or of synchronous rounds.

\FloatBarrier

\subsection{Diagnostics and limits of interpretation}\label{sec:experiment-diagnostics}

\begin{figure}[!htbp]
\centering
\includegraphics[width=.6\textwidth]{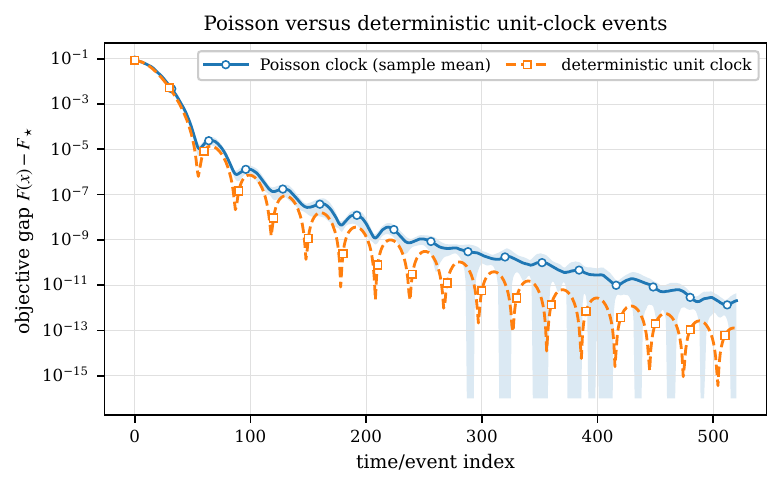}
\caption{Poisson versus unit-clock events on the certified regression.  The
deterministic-clock recurrence is not covered by the compensator proof.}
\label{fig:clock-tracking}
\end{figure}

The gradient-tracking comparator reaches an offline mean gap
$9.58\cdot10^{-8}$ at iteration $520$; the node-1 and worst-node gaps are
$9.60\cdot10^{-8}$ and $9.69\cdot10^{-8}$, with mean squared consensus error
$5.56\cdot10^{-9}$.  Replacing Poisson interarrivals by deterministic unit
intervals gives gap $1.20\cdot10^{-13}$ at index $520$, while the
Poisson-process sample mean is $2.05\cdot10^{-12}$ with standard error
$1.16\cdot10^{-12}$ (see Figure~\ref{fig:clock-tracking}).  On this instance the
random clock is needed for the proof and the exact simulator, but not
for convergence in practice.

Chebyshev filtering reaches $\|W_q-J\|_2\leq0.1$ on rings of size
$m=8,16,32,64$ at depths $4,8,16,31$, compared with
$11,45,179,717$ for plain lazy-Metropolis powers.  A separate smooth monotone
link has certified $\gamma=0.2$ and an observed one-point ratio $0.205405$ at
an explicit probe, only $2.70\%$ above the certificate.

The static, changing, and aligned-mode tests show how finite mixing changes the
variance; they do not show that it always beats $q=0$ or exact averaging.  The
certified $q=7$, $\theta^2=1/2$ ablation shows that damping absorbs
multiplicative noise, but it does not optimize $\theta$.  Agreement between
the continuized formulas, sample means, and realized event counts says nothing
about wall-clock time or about whether random clocks are necessary, and the
experiments include no lower bounds, cluster timings, fault tolerance, or
comparison with every local-state accelerated method.
\FloatBarrier

\section{Conclusion}\label{sec:conclusion}

Continuized quasar-convex acceleration survives finite network averaging.
One Lyapunov theorem handles any oracle with a predictable bias and a
state-dependent variance, and one lemma gives the moments of every
selected-row network oracle.  CQA-Net mixes local rows and disseminates one
uniformly selected row: the oracle is unbiased at every depth, incomplete
averaging enters only as variance, and the round count contains the
dissemination depth $D$.  The walking token removes every global primitive at
the price of a sequential walk of length of order $t_{\rm rel}\log m$; on a
fixed graph an importance weight keeps it unbiased, and on a changing graph
its bias is controlled by walk products.  D-CQA mixes local states: the network
average is a biased CQA process, the mixing accuracy must improve like
$T^{-2}$, and the round count
$\widetilde{\mathcal O}(\sqrt\chi\sqrt LR_0/(\gamma\theta\sqrt\varepsilon))$ has no
dissemination term, at the price of local smoothness and bounded
heterogeneity.  For shared-landscape finite sums, exact-snapshot restarts give a
geometric rate with fixed inner depth and batch size.

The results do not cover physical asynchrony, fault tolerance, or wall-clock
performance, and the broadcast-free implementations are simulated only on a
small instance ($m=16$, $d=18$).  Open questions include a lower bound on the per-event
communication of common-state methods with an unbiased oracle, gradient
tracking for D-CQA, which could remove the logarithm and the bounded
heterogeneity assumption, variance reduction for finite sums with mixed rather
than exact snapshots, and asynchronous clocks and message loss.

\section*{Acknowledgments}

The research was supported by Russian Science Foundation (project No. 23-11-00229), https://rscf.ru/en/project/23-11-00229/.

\bibliographystyle{plainnat}
\bibliography{references}

@article{chung1989diameters,
  author  = {Fan R. K. Chung},
  title   = {Diameters and eigenvalues},
  journal = {Journal of the American Mathematical Society},
  volume  = {2},
  number  = {2},
  pages   = {187--196},
  year    = {1989}
}

@inproceedings{chen1999lifting,
  author    = {Fang Chen and L{\'a}szl{\'o} Lov{\'a}sz and Igor Pak},
  title     = {Lifting {M}arkov chains to speed up mixing},
  booktitle = {Proceedings of the 31st Annual ACM Symposium on Theory of Computing},
  pages     = {275--281},
  publisher = {ACM},
  year      = {1999}
}

@inproceedings{kempe2003gossip,
  author    = {David Kempe and Alin Dobra and Johannes Gehrke},
  title     = {Gossip-based computation of aggregate information},
  booktitle = {Proceedings of the 44th Annual IEEE Symposium on Foundations of Computer Science},
  pages     = {482--491},
  publisher = {IEEE},
  year      = {2003}
}

@article{johansson2009randomized,
  author  = {Bj{\"o}rn Johansson and Maben Rabi and Mikael Johansson},
  title   = {A randomized incremental subgradient method for distributed optimization in networked systems},
  journal = {SIAM Journal on Optimization},
  volume  = {20},
  number  = {3},
  pages   = {1157--1170},
  year    = {2009}
}

@article{mao2020walkman,
  author  = {Xianghui Mao and Kun Yuan and Yubin Hu and Yuantao Gu and Ali H. Sayed and Wotao Yin},
  title   = {Walkman: A communication-efficient random-walk algorithm for decentralized optimization},
  journal = {IEEE Transactions on Signal Processing},
  volume  = {68},
  pages   = {2513--2528},
  year    = {2020}
}

@article{nedic2009,
  author  = {Nedi\'c, Angelia and Ozdaglar, Asuman},
  title   = {Distributed subgradient methods for multi-agent optimization},
  journal = {IEEE Transactions on Automatic Control},
  volume  = {54},
  number  = {1},
  pages   = {48--61},
  year    = {2009}
}

@inproceedings{mcmahan2017,
  author    = {McMahan, Brendan and Moore, Eider and Ramage, Daniel and Hampson, Seth and Ag{\"u}era y Arcas, Blaise},
  title     = {Communication-efficient learning of deep networks from decentralized data},
  booktitle = {Proceedings of the 20th International Conference on Artificial Intelligence and Statistics (AISTATS)},
  pages     = {1273--1282},
  year      = {2017}
}

@inproceedings{lian2017,
  author    = {Lian, Xiangru and Zhang, Ce and Zhang, Huan and Hsieh, Cho-Jui and Zhang, Wei and Liu, Ji},
  title     = {Can decentralized algorithms outperform centralized algorithms? {A} case study for decentralized parallel stochastic gradient descent},
  booktitle = {Advances in Neural Information Processing Systems 30 (NeurIPS)},
  pages     = {5330--5340},
  year      = {2017}
}

@article{nesterov1983,
  author  = {Nesterov, Yurii E.},
  title   = {A method for solving the convex programming problem with convergence rate $O(1/k^2)$},
  journal = {Doklady Akademii Nauk SSSR},
  volume  = {269},
  number  = {3},
  pages   = {543--547},
  year    = {1983}
}

@book{nesterov2004,
  author    = {Nesterov, Yurii},
  title     = {Introductory Lectures on Convex Optimization: A Basic Course},
  publisher = {Kluwer Academic Publishers},
  address   = {Boston},
  year      = {2004},
  doi       = {10.1007/978-1-4419-8853-9}
}

@article{guminov2023,
  author  = {Guminov, Sergey and Gasnikov, Alexander and Kuruzov, Ilya},
  title   = {Accelerated methods for weakly-quasi-convex optimization problems},
  journal = {Computational Management Science},
  volume  = {20},
  number  = {1},
  pages   = {36},
  year    = {2023},
  doi     = {10.1007/s10287-023-00468-w},
  note    = {Earlier version: arXiv:1710.00797}
}

@inproceedings{hinder2020,
  author    = {Hinder, Oliver and Sidford, Aaron and Sohoni, Nimit},
  title     = {Near-optimal methods for minimizing star-convex functions and beyond},
  booktitle = {Proceedings of the 33rd Conference on Learning Theory},
  series    = {Proceedings of Machine Learning Research},
  volume    = {125},
  pages     = {1894--1938},
  year      = {2020},
  publisher = {PMLR}
}

@inproceedings{wang2023,
  author    = {Wang, Jun-Kun and Wibisono, Andre},
  title     = {Continuized acceleration for quasar convex functions in non-convex optimization},
  booktitle = {International Conference on Learning Representations},
  year      = {2023},
  note      = {arXiv:2302.07851}
}

@inproceedings{fu2023,
  author    = {Fu, Qiang and Xu, Dongchu and Wilson, Ashia Camage},
  title     = {Accelerated stochastic optimization methods under quasar-convexity},
  booktitle = {Proceedings of the 40th International Conference on Machine Learning},
  series    = {Proceedings of Machine Learning Research},
  volume    = {202},
  pages     = {10431--10460},
  year      = {2023},
  publisher = {PMLR}
}

@article{hermant2024,
  author  = {Hermant, Julien and Aujol, Jean-Fran\c{c}ois and Dossal, Charles and Rondepierre, Aude},
  title   = {Study of the behaviour of {Nesterov} accelerated gradient in a non-convex setting: the strongly quasar-convex case},
  journal = {arXiv preprint arXiv:2405.19809},
  year    = {2024}
}

@inproceedings{gupta2024agnes,
  author    = {Gupta, Kanan and Siegel, Jonathan W. and Wojtowytsch, Stephan},
  title     = {Nesterov Acceleration despite Very Noisy Gradients},
  booktitle = {Advances in Neural Information Processing Systems},
  volume    = {37},
  pages     = {20694--20744},
  year      = {2024},
  doi       = {10.52202/079017-0653}
}

@inproceedings{lezane2026,
  author    = {Lezane, Clement and Langer, Sophie and Koolen, Wouter M.},
  title     = {Accelerated mirror descent for non-{Euclidean} star-convex functions},
  booktitle = {Proceedings of the 37th International Conference on Algorithmic Learning Theory},
  series    = {Proceedings of Machine Learning Research},
  volume    = {313},
  pages     = {1--41},
  year      = {2026},
  publisher = {PMLR}
}

@article{martinezrubio2025,
  author  = {Mart\'inez-Rubio, David},
  title   = {Smooth Quasar-Convex Optimization with Constraints},
  journal = {arXiv preprint arXiv:2510.01943},
  year    = {2025},
  note    = {Revised July 2026}
}

@article{gbaguidi2026,
  author  = {Gbaguidi, Em\'eric and Hermant, Julien},
  title   = {Accelerated stochastic zeroth-order quasar-convex optimization},
  journal = {arXiv preprint arXiv:2607.19965},
  year    = {2026}
}

@inproceedings{even2021,
  author    = {Even, Mathieu and Berthier, Rapha\"el and Bach, Francis and Flammarion, Nicolas and Hendrikx, Hadrien and Gaillard, Pierre and Massouli\'e, Laurent and Taylor, Adrien},
  title     = {Continuized accelerations of deterministic and stochastic gradient descents, and of gossip algorithms},
  booktitle = {Advances in Neural Information Processing Systems},
  volume    = {34},
  pages     = {28054--28066},
  year      = {2021}
}

@article{hardt2018,
  author  = {Hardt, Moritz and Ma, Tengyu and Recht, Benjamin},
  title   = {Gradient descent learns linear dynamical systems},
  journal = {Journal of Machine Learning Research},
  volume  = {19},
  number  = {29},
  pages   = {1--44},
  year    = {2018}
}

@inproceedings{scaman2017,
  author    = {Scaman, Kevin and Bach, Francis and Bubeck, S\'ebastien and Lee, Yin Tat and Massouli\'e, Laurent},
  title     = {Optimal algorithms for smooth and strongly convex distributed optimization in networks},
  booktitle = {Proceedings of the 34th International Conference on Machine Learning},
  series    = {Proceedings of Machine Learning Research},
  volume    = {70},
  pages     = {3027--3036},
  year      = {2017},
  publisher = {PMLR}
}

@inproceedings{hendrikx2019,
  author    = {Hendrikx, Hadrien and Bach, Francis and Massouli\'e, Laurent},
  title     = {Accelerated decentralized optimization with local updates for smooth and strongly convex objectives},
  booktitle = {Proceedings of the 22nd International Conference on Artificial Intelligence and Statistics},
  series    = {Proceedings of Machine Learning Research},
  volume    = {89},
  pages     = {897--906},
  year      = {2019},
  publisher = {PMLR}
}

@inproceedings{xu2020,
  author    = {Xu, Jinming and Tian, Ye and Sun, Ying and Scutari, Gesualdo},
  title     = {Accelerated primal-dual algorithms for distributed smooth convex optimization over networks},
  booktitle = {Proceedings of the 23rd International Conference on Artificial Intelligence and Statistics},
  series    = {Proceedings of Machine Learning Research},
  volume    = {108},
  pages     = {2381--2391},
  year      = {2020},
  publisher = {PMLR}
}

@inproceedings{kovalev2021adom,
  author    = {Kovalev, Dmitry and Shulgin, Egor and Richt\'arik, Peter and Rogozin, Alexander V. and Gasnikov, Alexander},
  title     = {{ADOM}: Accelerated decentralized optimization method for time-varying networks},
  booktitle = {Proceedings of the 38th International Conference on Machine Learning},
  series    = {Proceedings of Machine Learning Research},
  volume    = {139},
  pages     = {5784--5793},
  year      = {2021},
  publisher = {PMLR}
}

@inproceedings{metelev2023,
  author    = {Metelev, Dmitry and Rogozin, Alexander and Kovalev, Dmitry and Gasnikov, Alexander},
  title     = {Is consensus acceleration possible in decentralized optimization over slowly time-varying networks?},
  booktitle = {Proceedings of the 40th International Conference on Machine Learning},
  series    = {Proceedings of Machine Learning Research},
  volume    = {202},
  pages     = {24532--24554},
  year      = {2023},
  publisher = {PMLR}
}

@article{metelev2024,
  author  = {Metelev, Dmitry and Chezhegov, Savelii and Rogozin, Alexander and Beznosikov, Aleksandr and Sholokhov, Alexander and Gasnikov, Alexander and Kovalev, Dmitry},
  title   = {Decentralized finite-sum optimization over time-varying networks},
  journal = {arXiv preprint arXiv:2402.02490},
  year    = {2024}
}

@inproceedings{tian2022,
  author    = {Tian, Ye and Scutari, Gesualdo and Cao, Tianyu and Gasnikov, Alexander},
  title     = {Acceleration in distributed optimization under similarity},
  booktitle = {Proceedings of the 25th International Conference on Artificial Intelligence and Statistics},
  series    = {Proceedings of Machine Learning Research},
  volume    = {151},
  pages     = {5721--5756},
  year      = {2022},
  publisher = {PMLR}
}

@inproceedings{daneshmand2021,
  author    = {Daneshmand, Amir and Scutari, Gesualdo and Dvurechensky, Pavel and Gasnikov, Alexander},
  title     = {Newton method over networks is fast up to the statistical precision},
  booktitle = {Proceedings of the 38th International Conference on Machine Learning},
  series    = {Proceedings of Machine Learning Research},
  volume    = {139},
  pages     = {2398--2409},
  year      = {2021},
  publisher = {PMLR}
}

@article{qu2018,
  author  = {Qu, Guannan and Li, Na},
  title   = {Harnessing smoothness to accelerate distributed optimization},
  journal = {IEEE Transactions on Control of Network Systems},
  volume  = {5},
  number  = {3},
  pages   = {1245--1260},
  year    = {2018},
  doi     = {10.1109/TCNS.2017.2698261}
}

@article{nedic2017,
  author  = {Nedi\'c, Angelia and Olshevsky, Alex and Shi, Wei},
  title   = {Achieving geometric convergence for distributed optimization over time-varying graphs},
  journal = {SIAM Journal on Optimization},
  volume  = {27},
  number  = {4},
  pages   = {2597--2633},
  year    = {2017},
  doi     = {10.1137/16M1084316}
}

@article{lobanov2023,
  author  = {Lobanov, Aleksandr and Veprikov, Andrew and Konin, Georgiy and Beznosikov, Aleksandr and Gasnikov, Alexander and Kovalev, Dmitry},
  title   = {Non-smooth setting of stochastic decentralized convex optimization problem over time-varying graphs},
  journal = {Computational Management Science},
  volume  = {20},
  number  = {1},
  pages   = {48},
  year    = {2023},
  doi     = {10.1007/s10287-023-00479-7}
}

@inproceedings{hong2017,
  author    = {Hong, Mingyi and Hajinezhad, Davood and Zhao, Ming-Min},
  title     = {Prox-{PDA}: The proximal primal-dual algorithm for fast distributed nonconvex optimization and learning over networks},
  booktitle = {Proceedings of the 34th International Conference on Machine Learning},
  series    = {Proceedings of Machine Learning Research},
  volume    = {70},
  pages     = {1529--1538},
  year      = {2017},
  publisher = {PMLR}
}

@inproceedings{lu2021,
  author    = {Lu, Yucheng and De Sa, Christopher},
  title     = {Optimal complexity in decentralized training},
  booktitle = {Proceedings of the 38th International Conference on Machine Learning},
  series    = {Proceedings of Machine Learning Research},
  volume    = {139},
  pages     = {7111--7123},
  year      = {2021},
  publisher = {PMLR}
}

@article{agafonov2021,
  author  = {Agafonov, Artem and Dvurechensky, Pavel and Scutari, Gesualdo and Gasnikov, Alexander and Kamzolov, Dmitry and Lukashevich, Aleksandr and Daneshmand, Amir},
  title   = {An accelerated second-order method for distributed stochastic optimization},
  journal = {arXiv preprint arXiv:2103.14392},
  year    = {2021}
}

@inproceedings{kovalev2022sliding,
  author    = {Kovalev, Dmitry and Beznosikov, Aleksandr and Borodich, Ekaterina and Gasnikov, Alexander and Scutari, Gesualdo},
  title     = {Optimal Gradient Sliding and Its Application to Optimal Distributed Optimization under Similarity},
  booktitle = {Advances in Neural Information Processing Systems},
  volume    = {35},
  pages     = {33494--33507},
  year      = {2022}
}

@article{rogozin2019,
  author  = {Rogozin, Alexander and Uribe, C. A. and Gasnikov, Alexander and Malkovsky, Nikolay and Nedi\'c, Angelia},
  title   = {Optimal distributed optimization over slowly time-varying networks},
  journal = {IEEE Transactions on Control of Network Systems},
  volume  = {7},
  number  = {4},
  pages   = {1903--1915},
  year    = {2020}
}

@article{cao2026dcatalyst,
  author  = {Cao, Tianyu and Chen, Xiaokai and Scutari, Gesualdo},
  title   = {{DCatalyst}: A Unified Accelerated Framework for Decentralized Optimization},
  journal = {Journal of Machine Learning Research},
  volume  = {27},
  number  = {47},
  pages   = {1--57},
  year    = {2026}
}

@article{nabli2025dadao,
  author  = {Nabli, Adel and Oyallon, Edouard},
  title   = {Decentralized Asynchronous Optimization with {DADAO} Allows Decoupling and Acceleration},
  journal = {Journal of Machine Learning Research},
  volume  = {26},
  number  = {207},
  pages   = {1--48},
  year    = {2025}
}

@article{li2024agt,
  author  = {Li, Huan and Lin, Zhouchen},
  title   = {Accelerated Gradient Tracking over Time-Varying Graphs for Decentralized Optimization},
  journal = {Journal of Machine Learning Research},
  volume  = {25},
  number  = {274},
  pages   = {1--52},
  year    = {2024}
}

@inproceedings{takezawa2025spdo,
  author    = {Takezawa, Yuki and Jiang, Xiaowen and Rodomanov, Anton and Stich, Sebastian U.},
  title     = {Exploiting Similarity for Computation and Communication-Efficient Decentralized Optimization},
  booktitle = {Proceedings of the 42nd International Conference on Machine Learning},
  series    = {Proceedings of Machine Learning Research},
  volume    = {267},
  pages     = {58359--58405},
  year      = {2025},
  publisher = {PMLR}
}

@inproceedings{saad2025lower,
  author    = {Saad, El Mehdi and Lee, Wei-Cheng and Orabona, Francesco},
  title     = {New Lower Bounds for Non-Convex Stochastic Optimization through Divergence Decomposition},
  booktitle = {Proceedings of the 38th Conference on Learning Theory},
  series    = {Proceedings of Machine Learning Research},
  volume    = {291},
  pages     = {5054--5107},
  year      = {2025},
  publisher = {PMLR}
}

@inproceedings{ding2024gqc,
  author    = {Ding, Shihong and Yang, Long and Luo, Luo and Fang, Cong},
  title     = {Optimizing over Multiple Distributions under Generalized Quasar-Convexity Condition},
  booktitle = {Advances in Neural Information Processing Systems},
  volume    = {37},
  pages     = {4718--4764},
  year      = {2024},
  doi       = {10.52202/079017-0154}
}

@article{jin2020,
  author       = {Jikai Jin},
  title        = {On the Convergence of First Order Methods for Quasar-Convex Optimization},
  journal      = {arXiv preprint arXiv:2010.04937},
  year         = {2020}
}

@article{debrito2025,
  author       = {Jos{\'e} de Brito and Felipe Lara and Di Liu},
  title        = {Extending Linear Convergence of the Proximal Point Algorithm: The Quasar-Convex Case},
  journal      = {arXiv preprint arXiv:2509.04375},
  year         = {2025}
}

@article{ahookhosh2026,
  author       = {Masoud Ahookhosh and Jos{\'e} M. M. de Brito and Alireza Kabgani and Felipe Lara and Jinyun Yuan},
  title        = {Quasar-Convex Optimization: Fundamental Properties and High-Order Proximal-Point Methods},
  journal      = {arXiv preprint arXiv:2604.26735},
  year         = {2026}
}

@article{boyd2006randomized,
  author       = {Stephen Boyd and Arpita Ghosh and Balaji Prabhakar and Devavrat Shah},
  title        = {Randomized Gossip Algorithms},
  journal      = {IEEE Transactions on Information Theory},
  volume       = {52},
  number       = {6},
  pages        = {2508--2530},
  year         = {2006}
}

@inproceedings{karimireddy2020scaffold,
  author    = {Karimireddy, Sai Praneeth and Kale, Satyen and Mohri, Mehryar and Reddi, Sashank and Stich, Sebastian U. and Suresh, Ananda Theertha},
  title     = {{SCAFFOLD}: Stochastic Controlled Averaging for Federated Learning},
  booktitle = {Proceedings of the 37th International Conference on Machine Learning},
  series    = {Proceedings of Machine Learning Research},
  volume    = {119},
  pages     = {5132--5143},
  year      = {2020},
  publisher = {PMLR}
}

@inproceedings{gorbunov2021localsgd,
  author    = {Gorbunov, Eduard and Hanzely, Filip and Richt{\'a}rik, Peter},
  title     = {Local {SGD}: Unified Theory and New Efficient Methods},
  booktitle = {Proceedings of the 24th International Conference on Artificial Intelligence and Statistics},
  series    = {Proceedings of Machine Learning Research},
  volume    = {130},
  pages     = {3556--3564},
  year      = {2021},
  publisher = {PMLR}
}

@inproceedings{koloskova2020unified,
  author    = {Koloskova, Anastasia and Loizou, Nicolas and Boreiri, Sadra and Jaggi, Martin and Stich, Sebastian U.},
  title     = {A Unified Theory of Decentralized {SGD} with Changing Topology and Local Updates},
  booktitle = {Proceedings of the 37th International Conference on Machine Learning},
  series    = {Proceedings of Machine Learning Research},
  volume    = {119},
  pages     = {5381--5393},
  year      = {2020},
  publisher = {PMLR}
}

@inproceedings{liu2024finite,
  author    = {Liu, Yuxing and Chen, Lesi and Luo, Luo},
  title     = {Decentralized Convex Finite-Sum Optimization with Better Dependence on Condition Numbers},
  booktitle = {Proceedings of the 41st International Conference on Machine Learning},
  series    = {Proceedings of Machine Learning Research},
  volume    = {235},
  pages     = {30807--30841},
  year      = {2024},
  publisher = {PMLR}
}

@inproceedings{kuruzov2021sso,
  author    = {Kuruzov, Ilya A. and Stonyakin, Fedor S.},
  title     = {Sequential Subspace Optimization for Quasar-Convex Optimization Problems with Inexact Gradient},
  booktitle = {Advances in Optimization and Applications. OPTIMA 2021},
  series    = {Communications in Computer and Information Science},
  volume    = {1514},
  pages     = {19--33},
  publisher = {Springer, Cham},
  year      = {2021},
  doi       = {10.1007/978-3-030-92711-0_2},
  note      = {Earlier version: arXiv:2108.06297}
}

@article{kuruzov2024subspace,
  author        = {Kuruzov, Ilya and Stonyakin, Fedor},
  title         = {On Some Versions of Subspace Optimization Methods with Inexact Gradient Information},
  journal       = {arXiv preprint arXiv:2407.01051},
  year          = {2024},
  eprint        = {2407.01051},
  archivePrefix = {arXiv},
  primaryClass  = {math.OC}
}

\appendix

\section{Compensator calculus and localization}\label{app:generator-moments}

This appendix supplies the compensator, finite-horizon moment, and localization
facts used in Theorem~\ref{thm:master-inexact}.  None of them is needed to
implement the algorithms or to check the finite-mixing algebra.

\subsection{Predictable drift and the simulation law}

Let $\mathsf K_t(dg)$ be an $\mathbb F$-predictable version of the conditional
law of the event oracle.  For a smooth test function, the predictable extended
drift operator is the operator $\mathscr D_t$ defined in
\eqref{eq:predictable-drift-operator}.  If the oracle and mixer rules depend
only on a finite-dimensional current state, this is the extended generator of
a time-inhomogeneous piecewise-deterministic Markov process.  Under the more
general predictable history dependence used in the paper, it is a predictable
local characteristic; the proof needs no Markov property.

For a bounded stopping time $\rho\leq T$ and a test function of suitable
smoothness and growth, the stopped It\^o compensator identity reads
\begin{equation}\label{eq:ito-compensator-general}
\mathbb E\phi(\rho,x_\rho,z_\rho)=\mathbb E\phi(0,x_0,z_0)
 +\mathbb E\int_0^\rho(\mathscr D_t\phi)(t,x_{t-},z_{t-})\,dt.
\end{equation}
Identity~\eqref{eq:ito-compensator-general} follows from the unit predictable intensity of $n_t$ and the
usual jump It\^o formula and optional stopping of the compensated jump
integral.  It is applied below with $\rho=T\wedge\tau_R$.  Smooth truncation
first makes the stopped local martingale integrable; the stated moment
conditions then remove the truncation for the Lyapunov function.

\begin{lemma}[Law equivalence of on-line and pre-drawn clocks]
\label{lem:clock-law-equivalence}
Suppose the initial state is independent of i.i.d.\ $\mathop{\rm Exp}\nolimits(1)$
interarrivals.  After jump $k$, suppose the graph, the depth, the batch size,
and the mean-preserving input-output operator for jump $k+1$ are chosen from
the revealed past before $E_k$ is revealed.  Suppose also that the oracle/root
mark at $\tau_{k+1}$ is generated from an $\mathbb F$-predictable kernel, with
the root--row compatibility~\eqref{eq:lazy-root-row-kernel} when a lazy root is
used, and is conditionally independent of future interarrivals.  Then the
event-driven simulation using the pre-drawn
stream has the same finite-dimensional law as the jump process constructed in
$\mathbb F$ with unit compensator $dt$ and the same mark kernel.
\end{lemma}

\begin{proof}[Proof of Lemma~\ref{lem:clock-law-equivalence}]
The proof compares one conditional event transition and then iterates that
identity over the event index.

\medskip\noindent\textit{Step 1: one-event conditional law.}\quad
Let $\mathcal H_k$ contain the state and all revealed innovations after jump
$k$.  In the analytical construction, the graph, mixer, batch size, and all
other parameters for jump $k+1$ are selected as $\mathcal H_k$-measurable
functions.  An independent exponential variable $E_k$ is then drawn, the
deterministic flow forms the pre-jump query at $\tau_{k+1}$, and the fresh row
and root innovations are sampled from their predictable conditional kernels.
In the simulation construction the sequence
$(E_k)_{k\geq0}$ may have been generated in advance, but the same event
parameters are required to be functions of $\mathcal H_k$ alone and hence do
not use the unrevealed value of $E_k$.

Conditionally on $\mathcal H_k$, both constructions therefore use the same
iterated kernel: first the $\mathop{\rm Exp}\nolimits(1)$ law for $E_k$, then
the predictable mark kernel evaluated at the resulting pre-jump history.  This
need not be a product law because a state-dependent oracle kernel can depend on
the flowed query.  Within the mark kernel, positive-depth rows precede an
independent uniform root; at depth zero,
\eqref{eq:lazy-root-row-kernel} gives the equivalent product coupling of a
uniform root and latent row seeds.  The closed-form flow and jump map are
deterministic functions of the sampled quantities.

\medskip\noindent\textit{Step 2: iteration and compensator filtration.}\quad
Induction on $k$ gives equality of the joint laws of every finite event prefix.
Since i.i.d.\ positive exponential interarrivals are nonexplosive, these prefix
laws determine the finite-horizon path law.  The analytical filtration reveals
$E_k$ only through the occurrence of the next event, so $n_t$ has compensator
$t$ there.  The enlarged simulation filtration is useful for implementation
but is not the filtration in which the compensator calculation is performed.
\end{proof}

The coefficient $2/(\gamma t)$ has a removable singularity at the origin.
Between two events, $z_t=z_r$ and
$x_t-z_r=(r/t)^{2/\gamma}(x_r-z_r)$ for $t\geq r>0$.  Since $x_0=z_0$, the
continuous extension is constant before the first event.  This also explains
why the process can be simulated without a time discretization.

\subsection{Finite-horizon moments}

Assume that for constants $0\leq a_0,a_1<\infty$ and all event times $t\leq T$,
\begin{equation}\label{eq:linear-growth-oracle}
\mathbb E[\|g_t\|^2\mid\mathcal F_{t-}]
\leq a_0+a_1(\|x_{t-}\|^2+\|z_{t-}\|^2).
\end{equation}

\begin{lemma}[Finite-horizon state moments]\label{lem:finite-moments}
Under \eqref{eq:linear-growth-oracle}, for every finite $T$ there is
$C_T<\infty$ such that
\[
\mathbb E\sup_{0\leq t\leq T}
 (1+\|x_t\|^2+\|z_t\|^2)
\leq C_T(1+\|x_0\|^2+\|z_0\|^2).
\]
Global smoothness and the linear-growth oracle therefore imply the
finite-horizon integrability used in Theorem~\ref{thm:master-inexact}.
\end{lemma}

\begin{proof}
The exact flow controls the state between events.  A stopped compensator
estimate then handles the jump increments.

\medskip\noindent\textit{Step 1: between-event control.}\quad
Put $S_t=1+\|x_t\|^2+\|z_t\|^2$ and
$M_t=\sup_{0\leq r\leq t}S_r$.  Between events,
$x_t=(r/t)^{2/\gamma}x_r+[1-(r/t)^{2/\gamma}]z_r$; this is a convex
combination because $(r/t)^{2/\gamma}\in[0,1]$.  Hence
$\|x_t\|\leq\max\{\|x_r\|,\|z_r\|\}$, so the previously implicit
between-event constant is exactly one.

\medskip\noindent\textit{Step 2: stopped jump estimate.}\quad
At a jump, the two state increments are linear in $g_t$, with coefficients
$1/L$ and $c_t\leq\gamma\theta^2T/(2L)$.  The elementary square inequality
therefore gives a constant $C=C(L,\gamma,\theta,T)$ such that, pathwise,
\[
M_{T\wedge\tau_R}
\leq C M_0+C\int_0^{T\wedge\tau_R}
 (1+M_{t-}+\|g_t\|^2)\,dn_t.
\]
The marked compensator identity and~\eqref{eq:linear-growth-oracle} then give
\[
\mathbb E M_{T\wedge\tau_R}
\leq C M_0+C\int_0^T
 \bigl[(1+a_0)+(1+a_1)\mathbb E M_{t\wedge\tau_R}\bigr]dt.
\]
Gronwall gives a bound independent of $R$, now with constants depending also
on $a_0,a_1$, and Fatou yields the claim.  Global smoothness and the existence
of the fixed minimizer give the corresponding linear gradient-growth bounds;
the oracle bias and centered variance are controlled by its conditional second
moment.
\end{proof}

\subsection{Localization and exit times}

For local smoothness or local quasar geometry on a region $\mathcal B$, let
$\tau_{\mathcal B}$ be the first time at which a pre-jump query, a post-jump
state, or a line segment used by the smoothness model exits $\mathcal B$.
The same calculation proves the estimate only up to
$T\wedge\tau_{\mathcal B}$.  A global conclusion then requires a separate
invariance proof or an exit-probability bound.

\section{Network and communication lemmas}\label{app:network-proofs}

The spectral estimate proves Lemma~\ref{lem:chebyshev-filter}, the affine
heterogeneity bound closes Corollary~\ref{cor:closed-static}, and a Poisson
tail bound supports the work interpretation in Section~\ref{sec:experiments}.

\subsection{Chebyshev spectral estimate}

\begin{proof}[Proof of Lemma~\ref{lem:chebyshev-filter}]
The polynomial is controlled on the disagreement spectrum.  The resulting
spectral estimate gives the Frobenius row-noise factor.

\medskip\noindent\textit{Step 1: spectral contraction.}\quad
When $\chi=1$, all nonzero Laplacian eigenvalues coincide and the degree-one
polynomial $1-\lambda/\lambda_2$ gives exact averaging.  Assume $\chi>1$.
The normalized polynomial equals one at the consensus eigenvalue.  For every
$\lambda\in[\lambda_2,\lambda_m]$, its numerator argument lies in $[-1,1]$,
so the absolute numerator is at most one.  Hence
\[
\|W_q-J\|_2
\leq\frac1{\cosh(q\mathop{\rm arcosh}\nolimits((\chi+1)/(\chi-1)))}.
\]
Now $1/\cosh u\leq2e^{-u}$ and
$\mathop{\rm arcosh}\nolimits((\chi+1)/(\chi-1))
=\log((\sqrt\chi+1)/(\sqrt\chi-1))
\geq2/(\sqrt\chi+1)$, which gives
\eqref{eq:chebyshev-profile}.

\medskip\noindent\textit{Step 2: Frobenius factor.}\quad
Finally, mean preservation implies
$W_qJ=JW_q=J$ and Frobenius orthogonality of $J$ and $W_q-J$; therefore
$m^{-1}\|W_q\|_{\mathrm F}^2
=m^{-1}+m^{-1}\|W_q-J\|_{\mathrm F}^2
\leq m^{-1}+\|W_q-J\|_2^2$.
\end{proof}

The argument concerns a fixed symmetric graph or a certified common spectral
interval.  It is not a substitute for controlling noncommuting products on an
arbitrary changing graph.

\subsection{Affine heterogeneity from primitive local regularity}

Suppose
$\|[\nabla f_i(x)-\nabla F(x)]-[\nabla f_i(x^\star)-\nabla F(x^\star)]\|
\leq\delta_i\|x-x^\star\|$ and let
$\bar\delta^2=m^{-1}\sum_i\delta_i^2$.  Writing
$d_i(x)=\nabla f_i(x)-\nabla F(x)$ gives
$\|d_i(x)\|^2\leq2\|d_i(x^\star)\|^2
+2\delta_i^2\|x-x^\star\|^2$.  Averaging and using quadratic growth yields
\[
\zeta^2(x)
\leq2\zeta^2(x^\star)
 +\frac{4\bar\delta^2}{\mu}[F(x)-F(x^\star)],
\]
which is the certificate quoted after Corollary~\ref{cor:closed-static}.

\subsection{Actual versus expected event work}

The event count $n_T$ is Poisson with mean $T$, so the standard Bernstein bound
gives $|n_T-T|\leq\sqrt{2T\log(2/\delta)}+\tfrac23\log(2/\delta)$ with
probability at least $1-\delta$: the deterministic horizon is the expected work,
and the realized work fluctuates by order $\sqrt T$.

\section{Proof of the projected bound}\label{app:projected-bound}

The reserve left by $\theta<1$ is spent on the bias, which is why no gradient
integral remains on the left of~\eqref{eq:biased-localized-rate}.

\begin{proof}[Proof of Theorem~\ref{thm:master-inexact}(b)]
\medskip\noindent\textit{Step 1: localized drift.}\quad
Projection is nonexpansive relative to $x^\star\in\mathcal Q$, while
$z_t\in\mathcal Q$ implies
$|\langle z_t-x^\star,e_t\rangle|
\leq R_{\mathcal Q}\|e_t\|$.  Lemma~\ref{lem:predictable-drift} therefore
applies with this linear contribution.

\medskip\noindent\textit{Step 2: reserve accounting.}\quad
After inserting the multiplicative centered variance, the coefficient of
$\|\nabla F(x_{t-})\|^2$ is
$-(1-\theta^2)/2+V_\theta\eta$.  Completing the square gives
\[
-\left(\frac{1-\theta^2}{2}-V_\theta\eta\right)\|\nabla F\|^2
 +\theta^2\langle\nabla F,e_t\rangle
\leq
\frac{\theta^4\|e_t\|^2}
 {2(1-\theta^2-2V_\theta\eta)}.
\]
This uses all remaining gradient dissipation; hence no gradient-norm integral
appears on the left of \eqref{eq:biased-localized-rate}.

\medskip\noindent\textit{Step 3: compensation and normalization.}\quad
Apply the stopped compensator identity and remove the stopping exactly as in
the proof of Theorem~\ref{thm:master-inexact}(a), with Fatou's lemma for
$\Phi_{T\wedge\tau_R}$ and dominated convergence for the integrals.  Finally use
$A_T^{-1}\int_0^Tc_tdt=1/\gamma$ and
$A_T^{-1}\int_0^TA_tdt=T/3$ to obtain
\eqref{eq:biased-localized-rate}.
\end{proof}

\section{Variance reduction for shared-landscape finite sums}\label{app:vr}

This appendix proves the statements summarized in Remark~\ref{rem:vr}.  The
additive variance in Theorem~\ref{thm:static-network} grows with the horizon.
Quadratic growth permits fixed-length restarts, and an exact SVRG snapshot
turns the remaining stochastic error into a term proportional to the current
and snapshot gaps.

\paragraph{Finite sums and the control-variate row.}
Let $f_i(x)=n_i^{-1}\sum_{j=1}^{n_i}f_{ij}(x)$ and
$N_\Sigma:=\sum_{i=1}^m n_i$.  With a snapshot $y$ and the exact network
gradient $\nabla F(y)$, node $i$ uses
\begin{equation}\label{eq:svrg-row}
\widehat u_i(x;y)=\frac1b\sum_{r=1}^b
\bigl[\nabla f_{i,j_{i,r}}(x)-\nabla f_{i,j_{i,r}}(y)\bigr]+\nabla F(y).
\end{equation}
Conditionally on the stage-start sigma-field and the pre-event history, the
indices $j_{i,r}$ are i.i.d.\ uniform on $\{1,\ldots,n_i\}$ across $r$, the
same index is used in both gradients of each difference, and the batches are
conditionally independent across nodes and independent of the fresh clock and
root.  Cross-node independence may be replaced by the direct block-trace
domination~\eqref{eq:row-noise-trace-domination}.  An individual row mean need
not equal $\nabla f_i(x)$; only their average must equal $\nabla F(x)$, which
Lemma~\ref{lem:uniform-row} allows.

\begin{algorithm}[t]
\caption{Restarted exact-snapshot CQA-Net with variance reduction}
\label{alg:restarted-vr}
\begin{algorithmic}[1]
\Require stages $S$, stage horizon $T_{\rm st}$, inner depth $q$, batch $b$, common $x^{(0)}$
\For{$s=0,1,\ldots,S-1$}
  \State set the common snapshot $y^{(s)}\gets x^{(s)}$; every node computes $\nabla f_i(y^{(s)})$
  \State exactly reduce and broadcast $\nabla F(y^{(s)})$
  \State run CQA-Net from $x_0=z_0=y^{(s)}$ for horizon $T_{\rm st}$ with rows~\eqref{eq:svrg-row}
  \State set $x^{(s+1)}$ to the common stage output
\EndFor
\State \Return $x^{(S)}$
\end{algorithmic}
\end{algorithm}

\paragraph{Stage horizon and similarity.}
Let $\mathfrak F_s$ be the sigma-field at the beginning of stage $s$, after the
exact snapshot has been computed and before fresh inner clocks, component
samples, and roots are generated.  Exact-gradient acceleration and quadratic
growth give the deterministic stage factor $4L/(\mu\gamma^2\theta^2T^2)$, which
equals $1/4$ for
\begin{equation}\label{eq:stage-horizon}
T_{\rm st}=\frac4{\gamma\theta}\sqrt{\frac L\mu}.
\end{equation}
For comparison, quasar-convexity and quadratic growth imply
$\|\nabla F(x)\|^2\geq\gamma^2\mu[F(x)-F(x^\star)]/2$, so gradient descent with
step $1/L$ needs $\mathcal O(L/(\gamma^2\mu)\log(1/\varepsilon))$ iterations,
whereas restarted exact-gradient CQA needs
$\mathcal O(\sqrt{L/\mu}/(\gamma\theta)\log(1/\varepsilon))$ expected events.
The control variate is analyzed under two global \emph{gap-relative
similarity} conditions with constants $\beta,\Lambda_{\rm comp}\geq0$: for all
$x,y\in\mathbb R^d$ and every node $i$,
\begin{align}
\frac1m\sum_{i=1}^m
\|[\nabla f_i(x)-\nabla f_i(y)]-[\nabla F(x)-\nabla F(y)]\|^2
&\leq2\beta\bigl(F(x)-F(x^\star)+F(y)-F(x^\star)\bigr),
\label{eq:functional-similarity}\\
\mathbb E_j\|[\nabla f_{ij}(x)-\nabla f_{ij}(y)]
 -[\nabla f_i(x)-\nabla f_i(y)]\|^2
&\leq2\Lambda_{\rm comp}
\bigl(F(x)-F(x^\star)+F(y)-F(x^\star)\bigr).
\label{eq:component-similarity}
\end{align}

\begin{lemma}[Stage contraction under relative variance]\label{lem:relative-stage}
Suppose $F$ is $L$-smooth, $\gamma$-quasar-convex, and satisfies quadratic
growth with constant $\mu$.  Conditionally on $\mathfrak F_s$, let the
inner-stage oracle be unbiased and satisfy, for
$\mathfrak F_s$-measurable $0\leq a_s,c_s<\infty$,
$\varsigma_t^2\leq a_s(F(x_{t-})-F(x^\star))+c_s(F(x^{(s)})-F(x^\star))$.
Then, for a stage of length $T>0$,
\begin{equation}\label{eq:relative-stage-rate}
\mathbb E\!\left[F(x^{(s+1)})-F(x^\star)\mid\mathfrak F_s\right]
\leq e^{V_\theta a_sT/L}
\left[\frac{4L}{\mu\gamma^2\theta^2T^2}
 +\frac{V_\theta c_sT}{3L}\right]
\bigl(F(x^{(s)})-F(x^\star)\bigr).
\end{equation}
\end{lemma}

\begin{proof}
Condition on $\mathfrak F_s$ and reset inner time to zero.  The stopped
compensator inequality obtained from Lemma~\ref{lem:predictable-drift}, after
dropping nonpositive terms and using $A_t[F(x_t)-F(x^\star)]\leq\Phi_t$, reads
\[
\mathbb E[\Phi_{t\wedge\tau_R}\mid\mathfrak F_s]
\leq\Phi_0+\frac{V_\theta a_s}{L}
\int_0^t\mathbb E[\Phi_{r\wedge\tau_R}\mid\mathfrak F_s]dr+\frac{V_\theta c_s}{L}
\bigl(F(x^{(s)})-F(x^\star)\bigr)\int_0^tA_rdr .
\]
Because $\nabla F(x^\star)=0$, smoothness gives
$\|\nabla F(x)\|\leq L\|x-x^\star\|$ and
$F(x)-F(x^\star)\leq L\|x-x^\star\|^2/2$; with unbiasedness and the relative
variance bound this is a conditional linear second-moment bound, so
Lemma~\ref{lem:finite-moments} supplies the finite-horizon integrability and
nonexplosion.  Apply Lemma~\ref{lem:gronwall}(i), let $R\to\infty$ as in
Theorem~\ref{thm:master-inexact}, and use
$\Phi_0=\|x^{(s)}-x^\star\|^2/2\leq[F(x^{(s)})-F(x^\star)]/\mu$,
$\Phi_T\geq A_T[F(x^{(s+1)})-F(x^\star)]$, and $A_T^{-1}\int_0^TA_tdt=T/3$.
Division by $A_T$ gives~\eqref{eq:relative-stage-rate}.
\end{proof}

\begin{lemma}[Similarity-aware network oracle]\label{lem:similarity-oracle}
Under~\eqref{eq:functional-similarity} and~\eqref{eq:component-similarity}, let
$x,y$ be predictable and $W_q$ mean preserving, and assume conditional
cross-node independence of the mini-batch differences or the corresponding
block-trace domination.  The uniform-row estimate built from~\eqref{eq:svrg-row}
is conditionally unbiased for $\nabla F(x)$ and satisfies
\begin{equation}\label{eq:similarity-oracle-rate}
\mathbb E[\|\widehat g^{(k)}-\nabla F(x)\|^2\mid\mathcal F_{\tau_k-}]
\leq2K_{q,b}\bigl(F(x)-F(x^\star)+F(y)-F(x^\star)\bigr),
\end{equation}
where $K_{q,b}=\beta\|W_q-J\|_2^2
 +(\Lambda_{\rm comp}/b)\,m^{-1}\|W_q\|_{\mathrm F}^2$.
\end{lemma}

\begin{proof}
Condition on $\mathcal F_{\tau_k-}$.  The mean of row $i$ is
$\nabla f_i(x)-\nabla f_i(y)+\nabla F(y)$, whose network average is
$\nabla F(x)$, so Lemma~\ref{lem:uniform-row}(ii) with a uniform root and
$\omega\equiv1$ gives unbiasedness and the variance bound $\mathcal B^{(k)}$.
Its first term involves the centered mean rows
$[\nabla f_i(x)-\nabla f_i(y)]-[\nabla F(x)-\nabla F(y)]$, and
\eqref{eq:functional-similarity} bounds it by $2\beta\|W_q-J\|_2^2$ times the sum
of the two gaps.  In the noise term the same component is used at $x$ and $y$,
so~\eqref{eq:component-similarity} bounds one centered difference, averaging $b$
independent copies divides its variance by $b$, and the block-trace condition
contributes the factor $m^{-1}\|W_q\|_{\mathrm F}^2$.
\end{proof}

\begin{theorem}[Exact-snapshot restarted variance reduction]\label{thm:vr-main}
Suppose $F$ is $L$-smooth, $\gamma$-quasar-convex with respect to $x^\star$, and
satisfies quadratic growth with constant $\mu$ relative to the same minimizer.
Assume~\eqref{eq:functional-similarity}, \eqref{eq:component-similarity}, the
sampling conditions stated after~\eqref{eq:svrg-row}, and an exact snapshot
all-reduce on every stage.  Let $\theta\in(0,1]$ and use $T_{\rm st}$
from~\eqref{eq:stage-horizon}.  If
\begin{equation}\label{eq:vr-threshold}
K_{q,b}\leq K^\star
:=\frac{L\log(4/3)}{2V_\theta T_{\rm st}}
=\frac{\gamma\theta\log(4/3)}{8V_\theta}\sqrt{L\mu},
\end{equation}
then every stage of Algorithm~\ref{alg:restarted-vr} contracts the conditional
expected gap by
$r_{\rm st}=\frac43\bigl(\frac14+\frac13\log\frac43\bigr)\approx0.4612<\frac12$,
and
\begin{equation}\label{eq:vr-geometric-rate}
\mathbb E\bigl[F(x^{(S)})-F(x^\star)\bigr]
\leq r_{\rm st}^{S}\bigl(F(x^{(0)})-F(x^\star)\bigr),
\end{equation}
so $S=\lceil\log([F(x^{(0)})-F(x^\star)]/\varepsilon)/\log(1/r_{\rm st})\rceil$
stages suffice for accuracy $\varepsilon$.  On a fixed graph and for $q>0$, the
expected per-stage costs are at most $T_{\rm st}$ events, $2D+T_{\rm st}(q+D)$
critical-path rounds, $n_i+2bT_{\rm st}$ component gradients at node $i$,
$N_\Sigma+2mbT_{\rm st}$ in aggregate, and $2(m-1)+T_{\rm st}(2q|E|+M_{\rm br})$
transmitted vectors; for $q=0$, root-first sampling gives $n_i+2bT_{\rm st}/m$
component gradients at node $i$ and $N_\Sigma+2bT_{\rm st}$ in aggregate.
\end{theorem}

\begin{proof}
Fix $s$ and condition on $\mathfrak F_s$; the exact snapshot makes $x^{(s)}$
and $\nabla F(x^{(s)})$ fixed.  Lemma~\ref{lem:similarity-oracle} with
$y=x^{(s)}$ verifies Lemma~\ref{lem:relative-stage} with $a_s=c_s=2K_{q,b}$.  At
$T=T_{\rm st}$ the deterministic coefficient in~\eqref{eq:relative-stage-rate}
is $1/4$, and with $u=2V_\theta K_{q,b}T_{\rm st}/L$ the stage factor is
$e^u(1/4+u/3)$.  It increases in $u\geq0$, and~\eqref{eq:vr-threshold} gives
$u\leq\log(4/3)$, so the factor is at most $r_{\rm st}$; $\log(1+x)<x$ also gives
$r_{\rm st}<13/27<1/2$.  The tower property yields~\eqref{eq:vr-geometric-rate},
and solving $r_{\rm st}^S[F(x^{(0)})-F(x^\star)]\leq\varepsilon$ gives the stage
count.  For the costs, each stage begins with one exact local full pass, one
exact reduction, and one broadcast.  For $q>0$ every inner event evaluates the
same sampled components at the current point and at the snapshot on every node
($2b$ evaluations per node) and uses $q$ mixing rounds and one broadcast; for
$q=0$ only the root, revealed first, evaluates the two batches, and conditional
uniformity spreads its expected work equally over the nodes.
\end{proof}

The similarity conditions are strong.  A distance-based certificate, for
instance a uniform bound $\delta^2$ on
$m^{-1}\sum_i\|\nabla^2f_i(u)-\nabla^2F(u)\|_2^2$, gives
$\beta=2\delta^2/\mu$ through
$\|x-y\|^2\leq4[F(x)-F(x^\star)+F(y)-F(x^\star)]/\mu$, and for a generic deviation
$\delta^2\asymp L^2$ the threshold~\eqref{eq:vr-threshold} then requires
$\|W_q-J\|_2^2\lesssim(\mu/L)^{3/2}$ up to constants; on unstructured nonlinear
data such certificates can be vacuous.  The shared-landscape class avoids the
factor $1/\mu$.

\begin{lemma}[Shared-landscape certificate]\label{lem:shared-landscape}
Let $F$ be $L$-smooth and attain its global minimum at $x^\star$, and suppose
$f_{ij}(x)=a_{ij}F(x)+\langle c_{ij},x\rangle+d_{ij}$.  Write
$a_i=n_i^{-1}\sum_ja_{ij}$ and assume $m^{-1}\sum_ia_i=1$ and that the averaged
affine gradient vanishes.  Then~\eqref{eq:functional-similarity}
and~\eqref{eq:component-similarity} hold with
\begin{equation}\label{eq:shared-landscape-constants}
\beta=2Lv_{\rm node},
\qquad
\Lambda_{\rm comp}=2Lv_{\rm comp},
\end{equation}
where $v_{\rm node}:=m^{-1}\sum_i(a_i-1)^2$ and
$v_{\rm comp}:=\max_in_i^{-1}\sum_j(a_{ij}-a_i)^2$.
\end{lemma}

\begin{proof}
Affine terms cancel in gradient differences, so the node- and component-level
differences are $(a_i-1)[\nabla F(x)-\nabla F(y)]$ and
$(a_{ij}-a_i)[\nabla F(x)-\nabla F(y)]$.  Since $x^\star$ is a global minimizer,
the descent lemma gives $\|\nabla F(u)\|^2\leq2L[F(u)-F(x^\star)]$, hence
$\|\nabla F(x)-\nabla F(y)\|^2\leq4L[F(x)-F(x^\star)+F(y)-F(x^\star)]$, and
comparing with the factor $2$ in~\eqref{eq:functional-similarity}
and~\eqref{eq:component-similarity} gives~\eqref{eq:shared-landscape-constants}.
\end{proof}

\begin{corollary}[Shared-landscape finite sums]\label{cor:vr-shared}
Under the assumptions of Theorem~\ref{thm:vr-main} with components of the form
in Lemma~\ref{lem:shared-landscape}, the geometric rate~\eqref{eq:vr-geometric-rate}
holds whenever
\begin{equation}\label{eq:vr-shared-threshold}
v_{\rm node}\|W_q-J\|_2^2+\frac{v_{\rm comp}}{b}\frac{\|W_q\|_{\mathrm F}^2}{m}
\leq\frac{\gamma\theta\log(4/3)}{16V_\theta}\sqrt{\frac\mu L}.
\end{equation}
With the Chebyshev protocol of Lemma~\ref{lem:chebyshev-filter}, the inner depth
$q$ and the batch $b$ that satisfy~\eqref{eq:vr-shared-threshold} do not depend
on the target accuracy.
\end{corollary}

\begin{proof}
Insert~\eqref{eq:shared-landscape-constants} into~\eqref{eq:similarity-oracle-rate}:
$K_{q,b}=2L[v_{\rm node}\|W_q-J\|_2^2+(v_{\rm comp}/b)m^{-1}\|W_q\|_{\mathrm F}^2]$,
and~\eqref{eq:vr-threshold} becomes~\eqref{eq:vr-shared-threshold}.  Both sides
of~\eqref{eq:vr-shared-threshold} are independent of $\varepsilon$, and
Lemma~\ref{lem:chebyshev-filter} makes $\|W_q-J\|_2$ and $m^{-1}\|W_q\|_{\mathrm F}^2$
as small as needed by a fixed $q$.
\end{proof}

\begin{figure}[!htbp]
\centering
\includegraphics[width=.6\textwidth]{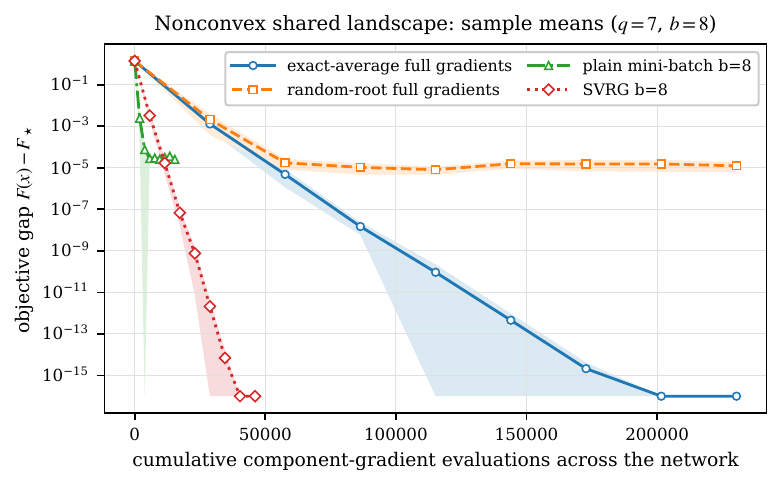}
\caption{Exact-snapshot variance reduction on the shared-landscape nonconvex
finite sum: objective gap versus component-gradient evaluations, including the
snapshot cost.  Curves are sample means with pointwise normal-approximation
bands.}
\label{fig:vr}
\end{figure}

\paragraph{Numerical check.}
Take the coordinatewise leaky-$\tanh$ landscape in dimension $12$ with
$\gamma=0.90$, $L=0.19219$, and $\mu=0.016875$, and components
$f_{ij}(x)=a_{ij}F(x)+\langle c_{ij},x\rangle$ on the $16$-node ring.  The
Hessian of $F$ has the negative eigenvalue $-3.61\cdot10^{-2}$ at an explicit
probe.  Lemma~\ref{lem:shared-landscape} gives $\beta=2.3543\cdot10^{-2}$ and
$\Lambda_{\rm comp}=1.2012\cdot10^{-2}$; with $q=7$ and $b=8$,
$K_{7,8}=4.7307\cdot10^{-4}<K^\star=1.8431\cdot10^{-3}$, so the threshold of
Theorem~\ref{thm:vr-main} holds before the run ($K_{7,8}/K^\star=0.25667$).
Figure~\ref{fig:vr} shows geometric progress after including the repeated
exact-snapshot cost: at stage eight the gaps are numerical zero
($<10^{-14}$) for full gradients and for SVRG with $b=8$, whereas the
restarted plain $b=8$ method has gap $2.55\cdot10^{-5}$.  The plain and SVRG
stage costs are $1{,}920$ and $5{,}760$ component-gradient evaluations, the
latter including the exact snapshot.  On the unstructured nonlinear data of
Section~\ref{sec:experiments}, by contrast, the generic Hessian certificate
misses the threshold by the factor $K/K^\star\approx2.44\cdot10^4$ even at the
largest tested batch.

\section{Proofs of the quasar-convex model certificates}
\label{app:quasar-certificates}

This appendix proves Proposition~\ref{prop:model-certificate}.

\subsection{Strongly monotone links}

Let $u_r=a_r^\top x$, $v_r=a_r^\top x^\star$, and
$d_r=\varphi(u_r)-\varphi(v_r)$.  The factors $d_r$ and $u_r-v_r$ have the same
sign.  Since $|d_r|\leq M|u_r-v_r|$ and $\varphi'(u_r)\geq\alpha$,
\[
\langle\nabla F(x),x-x^\star\rangle
=\frac1N\sum_rd_r\varphi'(u_r)(u_r-v_r)
\geq\frac{2\alpha}{M}F(x).
\]
Together with $F(x^\star)=0$ and the convention $\gamma\leq1$, this proves the
quasar certificate.  The lower derivative bound also gives
$|d_r|\geq\alpha|u_r-v_r|$, and therefore
\[
F(x)\geq\frac{\alpha^2}{2}(x-x^\star)^\top
\left(\frac1N\sum_ra_ra_r^\top\right)(x-x^\star),
\]
which proves the quadratic-growth constant.

For an arbitrary data partition, use the scaled local objective
\[
f_i(x)=\frac{m}{2N}\sum_{r\in\mathcal I_i}
\bigl[\varphi(a_r^\top x)-\varphi(a_r^\top x^\star)\bigr]^2.
\]
Then $m^{-1}\sum_i f_i=F$ even when the local sample counts are unequal.  A
noisy or misspecified response leaves a nonzero residual at $x^\star$ and is
not covered by this realizable certificate.

\subsection{The leaky-\texorpdfstring{$\tanh$}{tanh} smoothness envelope}

For $\alpha\in(0,1]$ and
$\varphi_\alpha(u)=\alpha u+(1-\alpha)\tanh u$, the Hessian is
$N^{-1}\sum_rc_r(x)a_ra_r^\top$, where
\[
c_r(x)=\varphi_\alpha'(u)^2+
[\varphi_\alpha(u)-\varphi_\alpha(v)]\varphi_\alpha''(u),
\qquad u=a_r^\top x,\qquad v=a_r^\top x^\star.
\]
Since
$\varphi_\alpha'(u)=\alpha+(1-\alpha)\mathop{\rm sech}\nolimits^2u$,
the squared derivative lies in $[\alpha^2,1]$.  Put $t=|\tanh u|$.  Then
\[
 |\varphi_\alpha''(u)|
 =2(1-\alpha)(1-t^2)t
 \leq \frac{4(1-\alpha)}{3\sqrt3}.
\]
The residual satisfies
\[
 |\varphi_\alpha(u)-\varphi_\alpha(v)|
 \leq \alpha|u|+\alpha|v|+2(1-\alpha).
\]
Also,
\[
 |u|\mathop{\rm sech}\nolimits^2u|\tanh u|
 \leq 4|u|e^{-2|u|}\leq\frac2e.
\]
Splitting the residual into the three terms above gives the explicit bound
\begin{align*}
 |[\varphi_\alpha(u)-\varphi_\alpha(v)]\varphi_\alpha''(u)|
 &\leq \frac{4\alpha(1-\alpha)}e
  +\frac{4\alpha(1-\alpha)}{3\sqrt3}|v|
  +\frac{8(1-\alpha)^2}{3\sqrt3}.
\end{align*}
Hence $|c_r(x)|\leq B_\alpha(v_r)$.  Finally,
\[
 \|\nabla^2F(x)\|_2
 \leq\frac1N\sum_{r=1}^N |c_r(x)|\,\|a_r\|^2
 \leq\frac1N\sum_{r=1}^N B_\alpha(v_r)\|a_r\|^2,
\]
which proves the global Lipschitz constant in
Proposition~\ref{prop:model-certificate}.

When the feature vectors are the coordinate basis, the Hessian is diagonal.
In that special case the generic sum can be sharpened to
$\|\nabla^2F(x)\|_2\leq\max_rB_\alpha(v_r)\|a_r\|^2/N$.

\subsection{Negative curvature within the certified class}

For $N=d=1$, $a_1=1$, $\alpha=1/2$ and $x^\star=-2$, the objective
$F(x)=[\varphi_{1/2}(x)-\varphi_{1/2}(-2)]^2/2$ is globally quasar-convex, and
outward-rounded interval evaluation of
$F''=\varphi_{1/2}'^2+[\varphi_{1/2}(x)-\varphi_{1/2}(-2)]\varphi_{1/2}''$ gives
$-0.253<F''(1)<-0.250$.

\section{Reproducibility statement}\label{app:reproducibility}

This appendix summarizes the setup behind the numerical figures and tables in
Section~\ref{sec:experiments}; it is independent of the mathematical proofs.

All experiments are CPU-only Python simulations using NumPy, pandas, and
Matplotlib.  Every stochastic stream has a distinct mixed-radix seed indexed
by experiment block, method, and replication, so all runs are
exactly repeatable.  The main certified monotone-link instance uses $m=16$
nodes, dimension $d=18$, a realizable leaky-$\tanh$ link, and the rounded
global constants reported in Section~\ref{sec:quasar-overview}; the
finite-mixing and damping examples use the explicit phase-fixed mode
$v_i=\sqrt{2}\cos(\pi i/2+\pi/2)$ rather than a numerically selected
eigenvector.
Roots are sampled after the rows for $q>0$ and before the sole local row
computation for $q=0$, matching Lemma~\ref{lem:uniform-row}.

All mixing operators are verified against explicit mean-preservation
assertions, and an automated validation pass recomputes the spectral
profiles, similarity constants, event counts, and every summary and
first-passage value quoted in the paper directly from the raw runwise data.
The manuscript rounds all constants to a readable number of significant
digits; full-precision values are kept in machine-readable metadata.  The
simulation code, raw data, and validation scripts are available from the
authors upon reasonable request.

\end{document}